\documentclass[11pt,a4paper,reqno]{amsart}
\usepackage{a4wide}
\usepackage{amsaddr}
\usepackage{amsmath}
\usepackage{amsfonts, amssymb, amsthm}
\usepackage{paralist}
\usepackage[colorlinks=true]{hyperref}
\hypersetup{urlcolor=blue, citecolor=red}

\makeatletter
\renewcommand{\email}[2][]{%
	\ifx\emails\@empty\relax\else{\g@addto@macro\emails{,\space}}\fi%
	\@ifnotempty{#1}{\g@addto@macro\emails{\textrm{(#1)}\space}}%
	\g@addto@macro\emails{#2}%
}
\makeatother

\theoremstyle{definition}

\numberwithin{equation}{section}

\usepackage[scrtime]{prelim2e}

\usepackage[normalem]{ulem}
\usepackage{color}
\newcommand{\R}{{\mathbb R}}

\begin{document}
\title[Nonlinear Choquard  equation with combined nonlinearity]
{Asymptotic profiles  of ground state solutions  for Choquard equations with a  general local perturbation}
\author{Shiwang Ma}\email{shiwangm@nankai.edu.cn}
\address{School of Mathematical Sciences and LPMC, Nankai University\\ 
	Tianjin 300071, China}
\author{Vitaly Moroz}\email{v.moroz@swansea.ac.uk}
\address{Department of Mathematics, Swansea University\\ 
	Fabian Way, Swansea SA1 8EN, 	Wales, UK}


\keywords{Nonlinear Choquard equation; groundstate solution; normalized solution; concentration compactness; asymptotic behavior.}

\subjclass[2010]{Primary 35J60, 35Q55; Secondary 35B25, 35B40, 35R09, 35J91}

\date{}

\begin{abstract}
In this paper, we study the asymptotic behavior of ground state solutions  for  the nonlinear  Choquard equation with a general local perturbation
$$
-\Delta u+\varepsilon u=(I_\alpha \ast |u|^{p})|u|^{p-2}u+ g(u), 
\quad {\rm in} \ \mathbb R^N,
 \eqno(P_\varepsilon)
 $$
 where $N\ge 3$ is an integer, $p=\frac{N+\alpha}{N}$, or $\frac{N+\alpha}{N-2}$, $I_\alpha$ is the Riesz potential and $\varepsilon>0$ is a parameter. Under some mild conditions on  $g(u)$, we show that as  $\varepsilon\to \infty$, after {\em a suitable rescaling} the ground state solutions of $(P_\varepsilon)$ converge  to a particular solution of some limit equations, and establish a sharp asymptotic characterisation of such a rescaling, which depend in a non-trivial way on the asymptotic behavior of the function $g(s)$ at infinity  and the space dimension $N$. Based on this study, we also present some results on the existence and asymptotic behaviors of positive normalized solutions  of $(P_\varepsilon)$  with the normalization constraint $\int_{\mathbb R^N}|u|^2=a^2$.  Particularly,  we obtain the asymptotic behavior of positive normalized solutions of such a problem as  $a\to 0$ and $a\to \infty$.

\end{abstract}

\maketitle

\newpage

\section{Introduction}

We study standing--wave solutions of the nonlinear Schr\"odinger equation with attractive  combined  nonlinearity
$$
i\psi_t=\Delta \psi+(I_\alpha\ast |\psi|^{p})|\psi|^{p-2}\psi+f(|\psi|)\psi\quad\text{in $\R^N\times\R$,}
\eqno(1.1)
$$
where $N\ge 3$ is an integer, $\psi: \mathbb R^N\times \mathbb R\to \mathbb C$,  $p\in [\frac{N+\alpha}{N},\frac{N+\alpha}{N-2}]$, and $I_\alpha$ is the Riesz potential defined for every $x\in \mathbb R^N\setminus \{0\}$ by 
$$
I_\alpha (x)=\frac{A_\alpha(N)}{|x|^{N-\alpha}}, \qquad A_\alpha(N)=\frac{\Gamma(\frac{N-\alpha}{2})}{\Gamma(\frac{\alpha}{2})\pi^{N/2}2^\alpha},
$$
where $\Gamma$ denotes the Gamma function, $f$ is a continuous function  to be specified later.

A theory of NLS with combined powers nonlinearity was developed by Tao, Visan and Zhang \cite{Tao} and attracted a lot attention during the past decade (cf. \cite{Akahori-2,Akahori-3, Cazenave-1, Coles, Jeanjean-5,   Li-2, Li-1, Li-3, Li-4, Liu-1, Moroz-1, Sun-1} and further references therein).  
The nonlocal convolution term in (1.1) represents the Newtonian gravitational attraction between bosonic particles. The local term takes into account the short--range self--interaction between bosons. For mare background of Choquard equation, we refer the readers to \cite{MV2017} and references therein.

A standing--wave solutions of (1.1) with a frequency $\varepsilon >0$ is a finite energy solution in the form
$$
\psi(t,x)=e^{-i\varepsilon t}u(x).
$$
Put $g(u)=f(|u|)u$, then this ansatz yields   the equation for $u$ in the form
$$
-\Delta u+\varepsilon u=(I_\alpha \ast |u|^{p})|u|^{p-2}u+g(u),
\quad {\rm in} \ \mathbb R^N.
 \eqno(P_\varepsilon)
$$

A solution $u_\varepsilon\in H^1(\mathbb R^N)$ of $(P_\varepsilon)$ is a critical point of the Action functional defined by
$$
I_\varepsilon(u)=\frac{1}{2}\int_{\mathbb R^N}|\nabla u|^2+\frac{\varepsilon}{2}\int_{\mathbb R^N}|u|^2-\frac{1}{2p}\int_{\mathbb R^N}(I_\alpha\ast |u|^p)|u|^p-\int_{\mathbb R^N}G(u),
\eqno(1.2)
$$
where $G(u)=\int_0^ug(s)ds$. The existence and properties of ground state $u_\varepsilon\in H^1(\mathbb R^N)$ to $(P_\varepsilon)$ can  been obtained by using the argument as in \cite{Li-2,Li-1}.
In the present paper, we are interested in the limit asymptotic profile of the groundsates $u_\varepsilon$ of the  problem $(P_\varepsilon)$, and in the asymptotic behaviors of different norms of $u_\varepsilon$, as $\varepsilon\to\infty$.

\smallskip

Many authors  consider the following Schr\"odinger equation with two focusing  exponents $p$ and $q$:
$$
-\Delta u+\varepsilon u=|u|^{p-2}u+|u|^{q-2}u, \quad {\rm in} \ \mathbb R^N,
\eqno(1.3)
$$
where $N\ge 3$, $2<q<p\le 2^*:=\frac{2N}{N-2}$ and $\varepsilon>0$ is a parameter.  It is well-known \cite{ASM2012, LLT2017, ZZ2012} that when $p=2^*$ and 
$q\in (2, 2^*)$ if $N\ge 4$, or $q\in (4,6)$ if $N=3$,
 the equation (1.3) admits a ground state solution $u_\varepsilon$ for every $\varepsilon>0$. Recently,  T. Akahori et al. \cite{Akahori-4}, and Wei and Wu \cite{Wei-2}  proved that (1.3) admits a ground state for all $N\ge 3$ and $\varepsilon>0$ small, while, for $N=3$ and $q\in (2,4]$, (1.3) has no ground state for large $\varepsilon>0$. When  $p=2^*$ and $q\in (2,2^*)$,  T. Akahori et al. \cite{Akahori-3}   proved that for small $\varepsilon>0$ the ground state is unique and as $\varepsilon\to 0$, the unique ground state $u_\varepsilon$ tends to the unique positive solution of the equation $-\Delta u+u=u^{q-1}$. After a suitable rescaling,  T. Akahori et al. \cite{Akahori-2}  established a uniform decay estimate for the ground state $u_\varepsilon$, and then   proved  the uniqueness and nondegeneracy in $H^1_r(\mathbb R^N)$  of ground states $u_\varepsilon$  for $N\ge 5$ and large $\varepsilon>0$,  and show that for $N\ge 3$, $q\in (2, 2^*)$ if $N\ge 4$, or $q\in (4,6)$ if $N=3$,  as $\varepsilon\to\infty$, $u_\varepsilon$ tends to a particular solution of the critical Emden--Fowler equation.
More recently, by using a global bifurcation argument, Jeanjean, Zhang and Zhong \cite{Jeanjean-4} also study the asymptotic behavior of solutions as $\varepsilon\to 0$ and $\varepsilon\to \infty$  for the equation  (1.3) with a general {\em subcritical} nonlinearity and discuss the connection  to the existence, non-existence and multiplicity of  prescribed mass positive solutions to (1.3) with  the  associated $L^2$ constraint condition $\int_{\mathbb R^N}|u|^2=a^2$. But the precise asymptotic profiles of positive solutions are not addressed there.  For  other related papers, we refer the reader to \cite{Ma-1} and the references therein. 
 
Li, Ma and Zhang \cite{Li-2}  and Li and Ma \cite{Li-1} considered $(P_\varepsilon)$ with $g(u)=|u|^{q-2}u$,
where $q\in (2,2^*]$, and  established the existence, regularity and symmetry of the ground state solutions of $(P_\varepsilon)$. Especially, for any $\varepsilon>0$, the equation $(P_\varepsilon)$ admits a ground state solution if $p=\frac{N+\alpha}{N}$ and $q\in (2,2+\frac{4}{N})$,  or  $p=\frac{N+\alpha}{N-2}$ and $q\in (2,2^*)$ for $N\ge 4$, $q\in (4,6)$ for $N=3$. 
 More recently, in \cite{Ma-2}, the authors consider again the equation $(P_\varepsilon)$ with $g(u)=|u|^{q-2}u$, and establish the precisely asymptotic behaviors of positive ground states as $\varepsilon\to 0$ and $\varepsilon\to \infty$. However, to the best of our knowledge, so far nothing has been achieved as for  the precise asymptotic behaviors of ground states of $(P_\varepsilon)$ with a general perturbation  $g(u)$.

In \cite{Moroz-1}, the second author and  Muratov first  study the asymptotic properties of ground states for  a class of scalar field equations   with a defocusing large exponent $p$ and a focusing smaller exponent $q$. More precisely, the following equation 
$$
-\Delta u+\varepsilon u=|u|^{p-2}u-|u|^{q-2}u, \quad {\rm in} \ \mathbb R^N,
\eqno(1.4)
$$
is discussed, where $N\ge 3$, $q>p>2$.  Later, in \cite{Lewin-1},  M. Lewin and S. Nodari prove a general result about the uniqueness and non-degeneracy of positive radial solutions to the above equation, and then the non-degeneracy of the unique solution $u_\varepsilon$ allows them to derive its behavior in the two limits $\varepsilon\to 0$ and $\varepsilon\to \varepsilon_*$, where $\varepsilon_*$ is a  threshold of existence.  Amongst other things, a precise asymptotic  expression of $M(\varepsilon)=\|u_\varepsilon\|_2^2$ is obtained. This gives the uniqueness of energy minimizers at fixed mass in certain regimes. 
 
  In \cite{Liu-1}, Zeng Liu and the second author extend the results in \cite{Moroz-1} to a class of Choquard type equation
$$
-\Delta u+\varepsilon u=(I_\alpha \ast |u|^{p})|u|^{p-2}u- |u|^{q-2}u, 
\quad {\rm in} \ \mathbb R^N,
\eqno(1.5)
 $$
where $N\ge 3$ is an integer. Under some assumptions on the exponents $p$ and $q$, the limit profiles of the ground states are discussed in  the two cases $\varepsilon\to 0$ and $\varepsilon\to \infty$. But the precisely 
asymptotic behaviors of ground states was not studied in \cite{Liu-1}. 

\smallskip

Alternatively, one can search for solutions to $(P_\varepsilon)$ having prescribed mass, and in this case $\varepsilon\in \mathbb R$ is part of the unknown.  That is, for a fixed $a>0$,  search for $u\in H^1(\mathbb R^N)$ and $\lambda\in \mathbb R$ satisfying
$$
\left\{\begin{array}{rl}
&-\Delta u=\lambda u+(I_\alpha\ast |u|^p)|u|^{p-2}u+g(u), \  \ in  \  \mathbb R^N,\\
&u\in H^1(\mathbb R^N),  \   \    \   \int_{\mathbb R^N}|u|^2=a^2.
\end{array}\right.
\eqno(1.6)
$$
The solution of (1.6) is usually denoted by a pair $(u,\lambda)$ and  called a normalized solution.

This approach seems to be particularly meaningful 
from the physical point of view, and often offers a good insight of the dynamical properties of the standing--wave solutions for (1.1), such as stability or instability.

  It is standard to check that the Energy functional 
$$
E(u)=\frac{1}{2}\int_{\mathbb R^N}|\nabla u|^2-\frac{1}{2p}\int_{\mathbb R^N}(I_\alpha\ast |u|^p)|u|^p-\int_{\mathbb R^N}G(u)
\eqno(1.7)
$$
is of class $C^1$ and that a critical point of $E$ restricted to the (mass) constraint 
$$
S(c)=\{ u\in H^1(\mathbb R^N): \ \|u\|_2^2=a^2\}
\eqno(1.8)
$$
gives a solution to (1.6).  Here   $\lambda:=-\varepsilon$ arises as a Lagrange multiplier. We refer the readers to \cite{Li-3,Li-4, Sun-1} and the  references therein.

\smallskip

 For quite a long time the paper \cite{Jeanjean-2} was the only one dealing with existence of normalized solutions in cases when the energy is unbounded from below on the $L^2$-constraint. More recently, however, problems of this type received much attention. In particular, many authors have made intensive studies to the following nonlinear Schr\"odinger equation with combined powers nonlinearity  
$$
\begin{cases}
&-\Delta u=\lambda u+|u|^{p-2}u+\mu |u|^{q-2}u, \quad  in \  \   \mathbb R^N,\\
& u\in H^1(\mathbb R^N), \   \  \|u\|_2^2=a^2, \end{cases}
\eqno(1.9)
$$
 where $N\ge 3$, $2<q<p\le 2^*$, $\mu>0$ is a new parameter and $\lambda\in \mathbb R$ appears as a Lagrange multiplier. We refer the readers to  \cite{ Jeanjean-3, Jeanjean-4, Soave-1,  Soave-2, Wei-1, Wei-2} and references therein  for the existence and multiplicity of  normalized solutions to the  probelm (1.9).
 In \cite{Wei-1}, Wei and Wu study the existence and asymptotic behaviors of normalized solutions  for (1.9) with $p=2^*$, and  obtain some precise asymptotic behaviors of ground states and mountain-pass solutions as $\mu\to 0$ and $\mu$ goes to its upper bound. It is worth mentioned that some ODE technique plays a crucial role in the studies in \cite{Wei-1}. We refer the readers to 
 \cite{ Jeanjean-3, Jeanjean-4, Soave-1,  Soave-2, Wei-2} for  the  asymptotic behaviors of normalized solutions when the parameter $\mu$ varies in its range.  
 
  Roughly speaking, the parameter $\mu$ in (1.9) plays a role of changing some thresholds for the existence. But it does not change the qualitative properties of solutions. In a sense,  changing the parameter $\mu$ is equavelent to changing the mass $a>0$.  More precisely,  it follows from the reduction given in \cite{Wei-2} that finding a normalized solution of (1.9) with $p=2^*$ and mass constrained  condition $\int_{\mathbb R^N}|u|^2=a^2$
is equivalent to finding a nomalized solution of the problem
   $$
\left\{\begin{array}{rl}
&-\Delta u=\lambda u+|u|^{2^*-2}u+|u|^{q-2}u, \  \ in  \  \   \mathbb R^N,\\
&u\in H^1(\mathbb R^N),  \   \    \   \int_{\mathbb R^N}|u|^2=a^2\mu^{\frac{4}{2N-q(N-2)}}.
\end{array}\right.
\eqno(1.10)
$$ 
 In particular,  sending $\mu\to 0$ (resp. $\mu\to \infty$) is equavelent to sending the mass $a\to 0$ (resp. $a\to \infty$).

A very complete theory has been established as to the existence  and multiplicity of normalized solutions to (1.9). However, to the best of our knowledge, little has been done for the Schr\"odinger equation with combined general nonlinearity such as 
 $$
 \begin{cases}
&-\Delta u=\lambda u+|u|^{2^*-2}u+ g(u), \quad  in \  \   \mathbb R^N,\\
& u\in H^1(\mathbb R^N), \   \  \|u\|_2^2=a^2. \end{cases}
\eqno(1.11)
$$
The main reason is  that the so called fibreing map is complicated and the technique to handle the Schr\"odinger equation with combined powers nonlinearity does not work any more in this case. To investigate the existence and asymptotic behavior of normalized solutions to (1.11), it seems that  some new approach must be employed.

As for the problem (1.6) with $g(u)=\mu |u|^{q-2}u$, by adopting the technique used in \cite{ Jeanjean-3, Jeanjean-4, Soave-1,  Soave-2}, some similar existence results on normalized solutions  are already obtained in \cite{ Li-3, Li-4, Li-5, Li-6, Sun-1, Sun-2} in the whole possible range of  parameters, but the precisely asymptotic behavior of the normalized solutions are not addressed there. 

\smallskip 

The purpose of this paper is two-fold. Firstly, by using some new technique, we  give some precise asymptotic behaviors of ground state solutions of $(P_\varepsilon)$ under some  mild subcritical growth conditions on the function $g(u)$. Secondly,    as an application of our main results, we  present  some existence results of positive normalized solutions to (1.6) and their precise asymptotic behaviors.  To explore the existence and precise asymptotic behaviors of normalized solutions of (1.6), 
it seems that the classical method developed recently in \cite{ Jeanjean-3, Jeanjean-4, Soave-1,  Soave-2}  is hard to be adopted since the associated fibreing   map is more complicated in our setting. Besides, due to the appearance 
of  the nonlocal term, the ODE technique does not work well any more, this prevent us from  adopting  the technique used in \cite{Wei-1} to derivate  the precise asymptotic behaviors of normalized solutions.

\smallskip

\noindent
\textbf{Organization of the paper}. In Section 2, we state  the main results in this paper. In Section 3, we give  some preliminary results which are needed in the proof of our main results. Sections 4--5 are devoted to the proofs of Theorems 2.1--2.2, respectively.  Finally, in the last section, we present some existence results of positive normalized solutions and their  precise asymptotic behaviors.

\smallskip

\noindent
\textbf{Basic notations}. Throughout this paper, we assume $N\geq
3$. $ C_c^{\infty}(\mathbb{R}^N)$ denotes the space of the functions
infinitely differentiable with compact support in $\mathbb{R}^N$.
$L^p(\mathbb{R}^N)$ with $1\leq p<\infty$ denotes the Lebesgue space
with the norms
$\|u\|_p=\left(\int_{\mathbb{R}^N}|u|^p\right)^{1/p}$.
 $ H^1(\mathbb{R}^N)$ is the usual Sobolev space with norm
$\|u\|_{H^1(\mathbb{R}^N)}=\left(\int_{\mathbb{R}^N}|\nabla
u|^2+|u|^2\right)^{1/2}$. $ D^{1,2}(\mathbb{R}^N)=\{u\in
L^{2^*}(\mathbb{R}^N): |\nabla u|\in
L^2(\mathbb{R}^N)\}$. $H_r^1(\mathbb{R}^N)=\{u\in H^1(\mathbb{R}^N):
u\  \mathrm{is\ radially \ symmetric}\}$.  $B_r$ denotes the ball in $\mathbb R^N$ with radius $r>0$ and centered at the origin,  $|B_r|$ and $B_r^c$ denote its Lebesgue measure  and its complement in $\mathbb R^N$, respectively.  As usual, $C$, $c$, etc., denote generic positive constants.
 For any large  $\epsilon>0$ and two nonnegative functions $f(\epsilon)$ and  $g(\epsilon)$, we write

(1)  $f(\epsilon)\lesssim g(\epsilon)$ or $g(\epsilon)\gtrsim f(\epsilon)$ if there exists a positive constant $C$ independent of $\epsilon$ such that $f(\epsilon)\le Cg(\epsilon)$.

(2) $f(\epsilon)\sim g(\epsilon)$ if $f(\epsilon)\lesssim g(\epsilon)$ and $f(\epsilon)\gtrsim g(\epsilon)$.

If $|f(\epsilon)|\lesssim |g(\epsilon)|$, we write $f(\epsilon)=O((g(\epsilon))$. We also denote by $\Theta=\Theta(\epsilon)$ a  generic  positive function satisfying
$
C_1\epsilon \le \Theta(\epsilon)\le C_2\epsilon
$
for some positive numbers $C_1,C_2>0$, which are independent of $\epsilon$. Finally, if $\lim f(\epsilon)/g(\epsilon)=1$ as $\epsilon\to \infty$, then we write $f(\epsilon)\simeq g(\epsilon)$ as $\epsilon\to \infty$.

\section{ Main Results}

In the present paper, we study the asymptotic behavior of positive groundstates of  $(P_\varepsilon)$ involving  a critical exponent in the sense of Hardy-Littlewood-Sobolev inequality. More precisely, if  $p=\frac{N+\alpha}{N}$, or $\frac{N+\alpha}{N-2}$, we  show that as $\varepsilon\to \infty$, after {\em a suitable rescaling} the ground state solutions $u_\varepsilon$ of $(P_\varepsilon)$ converge in $H^1(\mathbb R^N)$ to a particular solution of some limit equations. We also  establish a sharp asymptotic characteritation of such a rescaling, and the precise asymptotic behaviors of $u_\varepsilon(0), \|\nabla u_\varepsilon\|_2^2, \|u_\varepsilon\|_2^2$ and the least energy of the groundstates. 

\vskip 3mm 

To state our main results, we make the following assumptions:

\smallskip 

{\bf (H1)} $g\in C([0,\infty), [0,\infty) )$ satisfies $g(s)=o(s)$ as $s\to 0$;

\smallskip

{\bf (H2)}  there exists $q\in (2,2^*)$ such that
$
\lim_{s\to +\infty}g(s)s^{1-q}=A>0.
$

\vskip 5mm 

When $p=\frac{N+\alpha}{N}$, $q\in (2,2+\frac{4\alpha}{N(2+\alpha)})$ and $\varepsilon\to \infty$, we are going to show that after a suitable rescaling the limit equation for $(P_\varepsilon)$ is given by the critical Hardy-Littilewood-Sobolev  equation
$$
 U=(I_\alpha\ast |U|^{\frac{N+\alpha}{N}})U^{\frac{\alpha}{N}} \quad {\rm in } \   \ \mathbb R^N.
\eqno(2.1)
$$
It is well-known that the radial ground states of $(2.1)$ are given by the function 
$$
U_1(x):=\left(\frac{A_0}{1+|x|^2}\right)^{\frac{N}{2}}
\eqno(2.2)
$$
with a suitable constant $A_0 > 0$, and the family of its rescalings 
$$
U_\rho(x):=\rho^{-\frac{N}{2}}U_1(x/\rho),  \quad \rho>0.
\eqno(2.3)
$$

If $p=\frac{N+\alpha}{N}$, $q\in (2+\frac{4\alpha}{N(2+\alpha)}, 2+\frac{4}{N})$ and $\varepsilon\to \infty$, then it turns out that after some suitable rescaling the limit equation for $(P_\varepsilon)$ is given by
$$
-\Delta W+W=AW^{q-1}.
\eqno(2.4)
$$
For convenience, we set
$$
S_1:=\inf_{u\in H^1(\mathbb R^N)\setminus \{0\}}\frac{\int_{\mathbb R^N}|u|^2}{(\int_{\mathbb R^N}(I_\alpha\ast |u|^{\frac{N+\alpha}{N}}|u|^{\frac{N+\alpha}{N}})^{\frac{N}{N+\alpha}}},
\eqno(2.5)
$$ 
and 
$$
S_q:=\inf_{u\in H^1(\mathbb R^N)\setminus \{0\}}\frac{\int_{\mathbb R^N}|\nabla u|^2+|u|^2}{(\int_{\mathbb R^N}|u|^q)^{\frac{2}{q}}}.
\eqno(2.6)
$$

\vskip 5mm 

\noindent{\bf Theorem 2.1.} {\it Assume  (H1), (H2) hold, $p=\frac{N+\alpha}{N}$ and $q\in (2,2+\frac{4}{N})$,  then  for large $\varepsilon>0$, the problem  $(P_\varepsilon)$ admits a positive ground state $u_\varepsilon\in H^1(\mathbb R^N)$,  which is  radially symmetric and radially nonincreasing. Moreover,   if 
$q\in (2,2+\frac{4\alpha}{N(2+\alpha)})$ and $\varepsilon\to \infty$, then  there exists $\zeta_\varepsilon\in (0,+\infty)$ verifying 
$$
 \zeta_\varepsilon\sim \varepsilon^{-\frac{N(q-2)}{\alpha[4-N(q-2)]}}
 $$
such that  as $\varepsilon\to \infty$, the rescaled family of ground states
$
 \tilde w_\varepsilon(x)=\varepsilon^{-\frac{N}{2\alpha}}\zeta_\varepsilon^{\frac{N}{2}}u_\varepsilon(\zeta_\varepsilon x) 
$
 converges in $H^1(\mathbb R^N)$ to the extremal function  $U_{\rho_0}$ with
$$
\rho_0=\left(\frac{2q\int_{\mathbb R^N}|\nabla U_1|^2}{NA(q-2)\int_{\mathbb R^N}|U_1|^q}\right)^{\frac{2}{4-N(q-2)}}.
\eqno(2.7)
$$
If $q=2+\frac{4\alpha}{N(2+\alpha)}$ and $\varepsilon\to\infty$, then up to a subsequence,  the rescaled family of ground states
$
v_\varepsilon(x)=\varepsilon^{-\frac{1}{q-2}}u_\varepsilon(\varepsilon^{-\frac{1}{2}}x)
$
converges in $H^1(\mathbb R^N)$ to a positive solution of the equation
$$
-\Delta v+v=(I_\alpha\ast |v|^{\frac{N+\alpha}{N}})|v|^{\frac{N+\alpha}{N}-2}v+A|v|^{\frac{4\alpha}{N(2+\alpha)}}v.
$$
If $q\in (2+\frac{4\alpha}{N(2+\alpha)}, 2+\frac{4}{N})$ and $\varepsilon\to \infty$, then  the rescaled family of ground states
$
v_\varepsilon(x)=\varepsilon^{-\frac{1}{q-2}}u_\varepsilon(\varepsilon^{-\frac{1}{2}}x)
$
converges in $H^1(\mathbb R^N)$ to the unique positive solution $W$ of the equation (2.4). 

Furthermore, as $\varepsilon\to \infty$,  the least energy $m_\varepsilon$ of the ground state satisfies 
$$
m_\varepsilon=\left\{\begin{array}{rcl}
\varepsilon^{\frac{N+\alpha}{\alpha}}\left[\frac{\alpha}{2(N+\alpha)}S_1^{\frac{N+\alpha}{\alpha}}-\Theta(\varepsilon^{-\frac{4\alpha-N(2+\alpha)(q-2)}{\alpha[4-N(q-2)]}})\right],\ &if& q\in (2,2+\frac{4\alpha}{N(2+\alpha)}),\\
\varepsilon^{\frac{2N-q(N-2)}{2(q-2)}}\left[ \frac{q-2}{2q}A^{-\frac{2}{q-2}}S_q^{\frac{q}{q-2}}+o_\varepsilon(1)\right], \qquad  &if& q\in (2+\frac{4\alpha}{N(2+\alpha)}, 2+\frac{4}{N}).
\end{array}\right.
$$
}

\noindent{\bf Corollary 2.1.} {\it Assume  (H1), (H2) hold, $p=\frac{N+\alpha}{N}$ and $q\in (2,2+\frac{4}{N})$. Let $u_\varepsilon$ be the ground state solution  of  $(P_\varepsilon)$, then as $\varepsilon\to \infty$
$$
u_\varepsilon(0)\sim \left\{\begin{array}{rcl}
\varepsilon^{\frac{2N}{\alpha[4-N(q-2)]}},  \  &if& q\in (2,2+\frac{4\alpha}{N(2+\alpha)}),\\
\varepsilon^{\frac{1}{q-2}},\quad\quad \quad  &if&  q\in [2+\frac{4\alpha}{N(2+\alpha)}, 2+\frac{4}{N}),
\end{array}\right.
$$
$$
\|u_\varepsilon\|_2^2\sim \left\{\begin{array}{rcl}
\varepsilon^{\frac{N}{\alpha}}, \quad\qquad  &if& q\in (2,2+\frac{4\alpha}{N(2+\alpha)}),\\
\varepsilon^{\frac{4-N(q-2)}{2(q-2)}},\quad &if & q\in [2+\frac{4\alpha}{N(2+\alpha)}, 2+\frac{4}{N}),
\end{array}\right.
$$
$$
\|\nabla u_\varepsilon\|_2^2\sim \|u_\varepsilon\|_q^q\sim \left\{\begin{array}{rcl}
\varepsilon^{\frac{N[2N-q(N-2)]}{\alpha[4-N(q-2)]}}, \quad &if&q\in (2,2+\frac{4\alpha}{N(2+\alpha)}),\\
\varepsilon^{\frac{2N-q(N-2)}{2(q-2)}}, \quad &if & q\in [2+\frac{4\alpha}{N(2+\alpha)}, 2+\frac{4}{N}),
\end{array}\right.
$$
$$
\int_{\mathbb R^N}(I_\alpha\ast |u_\varepsilon|^{\frac{N+\alpha}{N}})|u_\varepsilon|^{\frac{N+\alpha}{N}}\sim \left\{\begin{array}{rcl}
\varepsilon^{\frac{N+\alpha}{\alpha}}, \quad \quad \qquad &if&q\in (2,2+\frac{4\alpha}{N(2+\alpha)}),\\
\varepsilon^{\frac{(N+\alpha)[4-N(q-2)]}{2N(q-2)}}, \  &if& q\in [2+\frac{4\alpha}{N(2+\alpha)}, 2+\frac{4}{N}).
\end{array}\right.
$$ }

\vskip 3mm

\noindent{\bf Remark 2.1.}  Corollary 2.1 is a direct consequence of Theorem 2.1. Surprisingly, in Theorem 2.1,  if  $q\in (2,2+\frac{4\alpha}{N(2+\alpha)})$, then as $\varepsilon \to \infty$,  the rescaled family of ground states
$ \tilde w_\varepsilon$  converges  to the extremal function  $U_{\rho_0}$ in $H^1(\mathbb R^N)$ other than in $L^2(\mathbb  R^N)$.
Besides,  from the main body of the paper, the following  more precise estimates are valid:

If 
$q\in (2,2+\frac{4\alpha}{N(2+\alpha)})$ and $\varepsilon\to \infty$, then the ground state solutions satisfy 
$$
\|u_\varepsilon\|^2_2=\varepsilon^{\frac{N}{\alpha}}\left[S_1^{\frac{N+\alpha}{\alpha}}+ O(\varepsilon^{-\frac{4\alpha-N(2+\alpha)(q-2)}{\alpha[4-N(q-2)]}})\right],
$$
$$
\|\nabla u_\varepsilon\|^2_2=\frac{NA(q-2)}{2q}\|u_\varepsilon\|_q^q+o(\varepsilon^{\frac{N[2N-q(N-2)]}{\alpha[4-N(q-2)]}})\sim \varepsilon^{\frac{N[2N-q(N-2)]}{\alpha[4-N(q-2)]}}, 
$$
$$
\int_{\mathbb R^N}(I_\alpha\ast |u_\varepsilon|^{\frac{N+\alpha}{N}})|u_\varepsilon|^{\frac{N+\alpha}{N}}=\varepsilon^{\frac{N+\alpha}{\alpha}}\left[S_1^{\frac{N+\alpha}{\alpha}}+ O(\varepsilon^{-\frac{4\alpha-N(2+\alpha)(q-2)}{\alpha[4-N(q-2)]}})\right].
$$ 
If $q\in (2+\frac{4\alpha}{N(2+\alpha)}, 2+\frac{4}{N})$ and $\varepsilon\to \infty$, then the ground state solutions satisfy
$$
\|u_\varepsilon\|_2^2=\varepsilon^{\frac{4-N(q-2)}{2(q-2)}}\left[\frac{2N-q(N-2)}{2q}A^{-\frac{2}{q-2}}S_q^{\frac{q}{q-2}}
+o_\varepsilon(1)\right],
$$
$$
\|\nabla u_\varepsilon\|_2^2=\varepsilon^{\frac{2N-q(N-2)}{2(q-2)}}\left[\frac{N(q-2)}{2q}A^{-\frac{2}{q-2}}S_q^{\frac{q}{q-2}}
+o_\varepsilon(1)\right],
$$
$$
\|u_\varepsilon\|_q^q=\varepsilon^{\frac{2N-q(N-2)}{2(q-2)}}\left[A^{-\frac{q}{q-2}}S_q^{\frac{q}{q-2}}
+o_\varepsilon(1)\right].
$$
\vskip 3mm

When $p=\frac{N+\alpha}{N-2}$ and $\varepsilon\to \infty$, we are going to show that after a suitable rescaling the limit equation for $(P_\varepsilon)$ is given by the critical 
Choquard equation
$$
-\Delta V=(I_\alpha\ast |V|^{\frac{N+\alpha}{N-2}})V^{\frac{2+\alpha}{N-2}} \quad {\rm in } \   \ \mathbb R^N.
\eqno(2.8)
$$
It is well-known that the radial ground states of $(2.8)$ are given by the  Talenti function 
$$
V_1(x):=[N(N-2)]^{\frac{N-2}{4}}\left(\frac{1}{1+|x|^2}\right)^{\frac{N-2}{2}}
\eqno(2.9)
$$
and the family of its rescalings 
$$
V_\rho(x):=\rho^{-\frac{N-2}{2}}V_1(x/\rho),  \quad \rho>0.
\eqno(2.10)
$$
For convenience, we set 
$$
S_\alpha:=\inf_{v\in D^{1,2}(\mathbb R^N)\setminus\{0\}}\frac{\int_{\mathbb R^N}|\nabla v|^2}{\left(\int_{\mathbb R^N}(I_\alpha\ast |v|^{\frac{N+\alpha}{N-2}})|v|^{\frac{N+\alpha}{N-2}}\right)^{\frac{N-2}{N+\alpha}}}.
\eqno(2.11)
$$

\noindent{\bf Theorem 2.2.} {\it Assume (H1) and (H2) hold,  $p=\frac{N+\alpha}{N-2}$,  $q\in (2,2^*)$ if $N\ge 4$ and $q\in (4,6)$ if $N=3$, then for large $\varepsilon>0$, the problem  $(P_\varepsilon)$ admits a positive ground state $u_\varepsilon\in H^1(\mathbb R^N)$,  which is  radially symmetric and radially nonincreasing. Moreovver, 
 if $N\ge 5$, then for large $\varepsilon>0$, there exists $\zeta_\varepsilon\in (0,+\infty)$
verifying  
$$\zeta_\varepsilon\sim \varepsilon^{-\frac{2}{(N-2)(q-2)}} $$
such that  as $\varepsilon\to \infty$, the rescaled family of ground states
$
w_\varepsilon(x)=\zeta_\varepsilon^{\frac{N-2}{2}}u_\varepsilon(\zeta_\varepsilon x)
$
converges in $H^1(\mathbb R^N)$ to  the extremal function $V_{\rho_0}$ with 
 $$
 \rho_0=\left(\frac{A[2N-q(N-2)]\int_{\mathbb R^N}|V_1|^q}{2q\int_{\mathbb R^N}|V_1|^2}\right)^ {\frac{2}{(N-2)(q-2)}}.
 \eqno(2.12)
 $$
 If  $N=4,3$,  
then  for large  $\varepsilon>0$, there exists $\zeta_\varepsilon\in (0,+\infty)$
verifying  
$$
\zeta_\varepsilon \sim\left\{\begin{array}{rcl}
(\varepsilon\ln\varepsilon)^{-\frac{1}{q-2}},   \  \quad if \   \ N=4,\\
\varepsilon^{-\frac{1}{q-4}},  \quad \qquad if  \    \  N=3,
\end{array}\right.
$$
such that 
the rescaled family of ground states
$
\tilde w_\varepsilon(x)=\zeta_\varepsilon^{\frac{N-2}{2}}u_\varepsilon(\zeta_\varepsilon x)
$
satisfies 
$$
\|\nabla \tilde w_\varepsilon\|_2^2\sim \|\tilde w_\varepsilon\|_{q}^{q}\sim \int_{\R^N}(I_\alpha\ast |\tilde w_\varepsilon|^{\frac{N+\alpha}{N-2}})|\tilde w_\varepsilon|^{\frac{N+\alpha}{N-2}}\sim 1, \quad 
  \|\tilde w_\varepsilon\|_2^2\sim \left\{\begin{array}{rcl}
\ln\varepsilon,   \quad  \quad if \   \ N=4,\\
\varepsilon^{\frac{6-q}{2(q-4)}},  \quad if  \    \  N=3,
\end{array}\right.
$$
and as $\varepsilon\to \infty$, $\tilde w_\varepsilon$ converges in $D^{1,2}(\mathbb R^N)$ and $L^{s}(\mathbb R^N)$  to $V_1$ for any $s>\frac{N}{N-2}$.

Furthermore, as $\varepsilon\to \infty$, the least energy $m_\varepsilon$ of the ground state satisfies 
$$
\frac{2+\alpha}{2(N+\alpha)}S_\alpha^{\frac{N+\alpha}{2+\alpha}}-m_\varepsilon\sim \left\{\begin{array}{rcl} \varepsilon^{-\frac{2N-q(N-2)}{(N-2)(q-2)}},  \quad  if \   \ N\ge 5,\\
(\varepsilon\ln\varepsilon)^{-\frac{4-q}{q-2}},   \quad  if \   \ N=4,\\
 \varepsilon^{-\frac{6-q}{2(q-4)}},    \qquad  if  \    \  N=3.
 \end{array}\right.
 $$}

\noindent{\bf Corollary 2.2.} {\it Assume (H1) and (H2) hold,  $p=\frac{N+\alpha}{N-2}$,  $q\in (2,2^*)$ if $N\ge 4$ and $q\in (4,6)$ if $N=3$.  
Let $u_\varepsilon$ be the ground state solution of  $(P_\varepsilon)$, then
 as $\varepsilon\to \infty$
 $$
u_\varepsilon(0)  \sim\left\{\begin{array}{rcl}
\varepsilon^{\frac{1}{q-2}}, \qquad \quad if \ \ N\ge 5,\\
(\varepsilon\ln\varepsilon)^{\frac{1}{q-2}},   \  \quad if \   \ N=4,\\
\varepsilon^{\frac{1}{2(q-4)}}, \  \qquad if  \    \  N=3,
\end{array}\right.
$$
$$
 \|u_\varepsilon\|_q^q\sim\left\{\begin{array}{rcl}
 \varepsilon^{-\frac{4}{(N-2)(q-2)}},\quad if \  \ N\ge 5,\\ 
(\varepsilon\ln\varepsilon)^{-\frac{4-q}{q-2}},   \  \quad if \   \ N=4,\\
\varepsilon^{-\frac{6-q}{2(q-4)}},   \qquad if  \    \  N=3,
\end{array}\right.
\quad
\|u_\varepsilon\|_2^2\sim\left\{\begin{array}{rcl}
\varepsilon^{-\frac{4}{(N-2)(q-2)}}, \quad \qquad  if \  \  N\ge 5,\\
\varepsilon^{-\frac{2}{q-2}}(\ln\varepsilon)^{-\frac{4-q}{q-2}},   \  \quad if \   \ N=4,\\
\varepsilon^{-\frac{q-2}{2(q-4)}},  \qquad \qquad if  \    \  N=3,
\end{array}\right.
$$
$$
\|\nabla u_\varepsilon\|^2_2=S_\alpha^{\frac{N+\alpha}{2+\alpha}}+ \left\{\begin{array}{rcl} 
O(\varepsilon^{-\frac{2N-q(N-2)}{(N-2)(q-2)}}), \quad if \  \ N\ge 5,\\
O((\varepsilon\ln\varepsilon)^{-\frac{4-q}{q-2}}),   \quad  if  \   \ N=4,\\
O( \varepsilon^{-\frac{6-q}{2(q-4)}}),   \qquad  if \    \  N=3,
 \end{array}\right.
$$
$$ 
\int_{\mathbb R^N}(I_\alpha\ast |u_\varepsilon|^{\frac{N+\alpha}{N-2}})|u_\varepsilon|^{\frac{N+\alpha}{N-2}}=S_\alpha^{\frac{N+\alpha}{2+\alpha}}+ \left\{\begin{array}{rcl} 
O(\varepsilon^{-\frac{2N-q(N-2)}{(N-2)(q-2)}}), \quad &if&  N\ge 5,\\
O((\varepsilon\ln\varepsilon)^{-\frac{4-q}{q-2}}),   \quad & if&  N=4,\\
O( \varepsilon^{-\frac{6-q}{2(q-4)}}),   \quad \quad   &if&  N=3.
 \end{array}\right.
$$}

\noindent{\bf Remark 2.2.} Corollary 2.2 is a direct consequence of Theorem 2.2.  In Theorem 2.2, we need not the assumption that $\alpha>N-4$  in the case $N\ge 5$, and if $\alpha>N-4$ holds true, then we can choose 
$
\zeta_\varepsilon= \varepsilon^{-\frac{2}{(N-2)(q-2)}}.
$
\vskip 5mm

\section{Preliminaries } 

In this section, we present some preliminary results which are needed in the proof of our main results. Firstly,  we consider the following Choquard type  equation with combined nonlinearites:
$$
-\Delta u+\varepsilon u=(I_{\alpha}\ast |u|^p)|u|^p+ g(u),  \quad \mathrm{in}\
\mathbb{R}^N,
\eqno(P_{\varepsilon})
$$
where $N\geq 3,\ \alpha\in(0,N)$, $p\in [\frac{N+\alpha}{N}, \frac{N+\alpha}{N-2}]$,
 $\varepsilon>0$ is a large parameter, and $g(u)$ is a function satisfying the assumptions (H1) and (H2) in the last section.

It has been proved in \cite[Theorem 2.1]{Li-1} that any weak solution of $(P_{\varepsilon})$ in
$H^1(\mathbb{R}^N)$ has additional regularity properties, which
allows us to establish the Poho\v{z}aev identity for all finite
energy solutions.

\smallskip

\noindent{\bf Lemma 3.1.}
{\it  Assume (H1) and (H2) hold.  If
$u\in H^1(\mathbb{R}^N)$ is a solution of $(P_{\varepsilon})$, then $u\in
W_{\mathrm{loc}}^{2,r}(\mathbb{R}^N)$ for every $r>1$. Moreover, $u$
satisfies the Poho\v{z}aev identity
$$
P_{\varepsilon}(u):=\frac{N-2}{2}\int_{\mathbb{R}^N}|\nabla
u|^2+\frac{N\varepsilon}{2}\int_{\mathbb{R}^N}|u|^2-\frac{N+\alpha}{2p}\int_{\mathbb{R}^N}(I_\alpha\ast
|u|^p)|u|^p-N\int_{\R^N}G(u)=0.
\eqno(3.1)
$$}

It is well known that any weak solution of  $(P_\varepsilon)$ corresponds to a critical point of the action  functionals $I_{\varepsilon}$ defined in (1.2), 
which is well defined and is of $C^1$ in $H^1(\R^N)$.  A nontrivial solution $u_\varepsilon\in H^1(\R^N)$  is called a ground-state if 
$$
I_\varepsilon(u_{\varepsilon})=m_{\varepsilon}:=\inf\{ I_{\varepsilon}(u): \ u\in H^1(\R^N)\setminus \{0\} \ {\rm and}  \  I'_{\varepsilon}(u)=0\}.
\eqno(3.2)
$$
In \cite{Li-2, Li-1} (see also the proof of the main results in \cite{Li-1}),  it has been shown that 
$$
m_{\varepsilon}=\inf_{u\in \mathcal M_{\varepsilon}}I_{\varepsilon}(u)=\inf_{u\in \mathcal P_{\varepsilon}}I_{\varepsilon}(u),
\eqno(3.3)
$$
where $\mathcal M_{\varepsilon}$ and $\mathcal P_{\varepsilon}$ are the correspoding Nehari  and Poho\v{z}aev manifolds defined by
$$
\mathcal M_{\varepsilon}:=\left\{ u\in H^1(\mathbb R^N)\setminus\{0\}  \ \left | \ \int_{\mathbb R^N}|\nabla u|^2+\varepsilon |u|^2=\int_{\mathbb R^N}(I_\alpha\ast |u|^p)|u|^p+\int_{\mathbb R^N}g(u)u \right. \right\}
$$
and 
$$
\mathcal P_{\varepsilon}:=\left\{ u\in H^1(\mathbb R^N)\setminus\{0\}  \ \left | \ P_{\varepsilon}(u)=0 \right. \right\},
$$
respectively.
Moreover, the following min-max descriptions are valid:

\smallskip

\noindent{\bf Lemma 3.2.}  {\it  Let
 $$
u_t(x)=\left\{\begin{array}{rcl} u(\frac{x}{t}), \quad if \  t>0,\\
0, \quad \  if \  \ t=0,
\end{array}\right.
$$
then 
$$
m_{\varepsilon}=\inf_{u\in H^1(\mathbb R^N)\setminus\{0\}}\sup_{t\ge 0}I_{\varepsilon}(tu)=\inf_{u\in H^1(\mathbb R^N)\setminus\{0\}}\sup_{t\ge 0}I_{\varepsilon}(u_t).
\eqno(3.4)
$$
In particular, we have $m_{\varepsilon}=I_{\varepsilon}(u_{\varepsilon})=\sup_{t>0}I_{\varepsilon}(tu_{\varepsilon})=\sup_{t>0}I_{\varepsilon}((u_{\varepsilon})_t)$.}
\smallskip

The following well
known Hardy-Littlewood-Sobolev inequality can be found in
\cite{Lieb-Loss 2001}.
\smallskip

\noindent{\bf Lemma 3.3.} {\it 
Let $p, r>1$ and $0<\alpha<N$ with $1/p+(N-\alpha)/N+1/r=2$. Let
$u\in L^p(\mathbb{R}^N)$ and $v\in L^r(\mathbb{R}^N)$. Then there
exists a sharp constant $C(N,\alpha,p)$, independent of $u$ and $v$,
such that
$$
\left|\int_{\mathbb{R}^N}\int_{\mathbb{R}^N}\frac{u(x)v(y)}{|x-y|^{N-\alpha}}\right|\leq
C(N,\alpha,p)\|u\|_p\|v\|_r.
$$
If $p=r=\frac{2N}{N+\alpha}$, then
$$
C(N,\alpha,p)=C_\alpha(N)=\pi^{\frac{N-\alpha}{2}}\frac{\Gamma(\frac{\alpha}{2})}{\Gamma(\frac{N+\alpha}{2})}\left\{\frac{\Gamma(\frac{N}{2})}{\Gamma(N)}\right\}^{-\frac{\alpha}{N}}.
$$}

\noindent{\bf Remark 3.1. }  By the Hardy-Littlewood-Sobolev inequality, for any $v\in L^s(\mathbb R^N)$ with $s\in (1,\frac{N}{\alpha})$, $I_\alpha\ast v\in L^{\frac{Ns}{N-\alpha s}}(\mathbb R^N)$ and 
$$
\|I_\alpha\ast v\|_{\frac{Ns}{N-\alpha s}}\le A_\alpha(N)C(N,\alpha, s)\|v\|_{s}.
\eqno(3.5)
$$

\smallskip

\noindent {\bf Lemma 3.4.} ( P. L. Lions \cite{Lions-1} )
{\it Let $r>0$ and $2\leq q\leq 2^{*}$. If $(u_{n})$ is bounded in $H^{1}(\mathbb{R}^N)$ and if
$$\sup_{y\in\mathbb{R}^N}\int_{B_{r}(y)}|u_{n}|^{q}\to0,\,\,\textrm{as\ }n\to\infty,$$
then $u_{n}\to0$ in $L^{s}(\mathbb{R}^N)$ for $2<s<2^*$. Moreover, if $q=2^*$, then $u_{n}\to0$ in $L^{2^{*}}(\mathbb{R}^N)$.
}

\smallskip

\noindent{\bf Lemma 3.5.}  ( Radial Lemma A.II, H. Berestycki and P. L. Lions \cite{Berestycki-1}  )
{\it  Let $N\ge 2$, then every radial function $u\in H^1(\mathbb R^N)$ is almost everywhere equal to a function $\tilde u(x)$, continuous for $x\not=0$, such that
$$
 |\tilde u(x)|\le  C_N |x|^{(1-N)/2}\|u\|_{H^1(\mathbb R^N)} \qquad for  \ |x|\ge \alpha_N,
\eqno(3.6)
$$
 where $C_N$ and $\alpha_N$ depend only on the dimension $N$.}
 \smallskip

\noindent {\bf Lemma 3.6.} ( Radial Lemma A.III. H. Berestycki and P. L. Lions \cite{Berestycki-1} ) 
 {\it Let $N\ge 3$, then every radial function $u$ in $D^{1,2}(\mathbb R^N)$ is almost everywhere equal to a function $\tilde u(x)$, continuous for $x\not=0$, such that
 $$
  |\tilde u(x)|\le C_N |x|^{(2-N)/2}\|u\|_{D^{1,2}(\mathbb R^N)} \qquad for  \ |x|\ge 1, 
\eqno(3.7)
  $$
  where $C_N$ only depends on $N$.} 
  \smallskip
  
  The following lemma is proved in \cite{Ma-2}.
  \smallskip
  
\noindent{\bf Lemma 3.7.} {\it Let $0<\alpha<N$, $0\le f(x)\in L^1(\mathbb R^N)$ be a radially symmetric function such that
$$
\lim_{|x|\to +\infty}f(|x|)|x|^N=0.
\eqno(3.8)
$$
If $\alpha\le 1$, we additionally assume that $f$ is monotone non-increasing. Then as $|x|\to +\infty$, we have 
$$
\int_{\mathbb R^N}\frac{f(y)}{|x-y|^{N-\alpha}}dy=\frac{\|f\|_{L^1}}{|x|^{N-\alpha}}+o\left(\frac{1}{|x|^{N-\alpha}}\right).
\eqno(3.9)
$$}
  
  \smallskip

The following Moser iteration lemma  is given in \cite[Proposition B.1]{Akahori-2}. See also  \cite{LiuXQ}  and \cite{GT}.

\smallskip

\noindent{\bf Lemma 3.8.} {\it Assume $N\ge 3$. Let $a(x)$ and $b(x)$ be functions on $B_4$, and let $u\in H^1(B_4)$ be a weak solution to 
$$
-\Delta u+a(x)u=b(x)u \qquad  in \  \ B_4.
\eqno(3.10)
$$
Suppose that $a(x)$ and $u$ satisfy that 
$$
a(x)\ge 0 \quad for \ a. e. \ x\in B_4, 
\eqno(3.11)
$$
and 
$$
 \int_{B_4}a(x)|u(x)v(x)|dx<\infty \quad for \ each \ v\in H_0^1(B_4).
\eqno(3.12)
$$
(i) Assume that for any $\varepsilon\in (0,1)$, there exists $t_\varepsilon>0$ such  that
$$
\|\chi_{[|b|>t_\varepsilon]}b\|_{L^{N/2}(B_4)}\le \varepsilon,
$$
where $[|b|>t]:=\{x\in B_4: \ |b(x)|>t\},$ and $\chi_A(x)$ denotes the characteristic function of $A\subset \mathbb R^N$. Then for any $r\in (0,\infty)$, there exists a constant $C(N,r,t_\varepsilon)$ such that 
$$
\||u|^{r+1}\|_{H^1(B_1)}\le C(N, r,t_\varepsilon)\|u\|_{L^{2^*}(B_4)}.
$$
(ii) \ Let $s>N/2$ and assume that $b\in L^s(B_4)$. Then there exits a constant $C(N,s,\|b\|_{L^s(B_4)})$ such that
$$
\|u\|_{L^\infty(B_1)}\le C(N,s,\|b\|_{L^s(B_4)})\|u\|_{L^{2^*}(B_4)}.
$$
Here, the constants $C(N,r, t_\varepsilon)$ and $C(N,s,\|b\|_{L^s(B_4)})$ in (i) and (ii) remain bounded as long as $r, t_\varepsilon$ and $\|b\|_{L^s(B_4)}$ are bounded.}




  \vskip 5mm

\section{Proof of Theorem 2.1}

In this section, we always assume that $2_\alpha:=\frac{N+\alpha}{N}$, $q\in (2, 2+\frac{4}{N})$ and $\varepsilon>0$ is a large  parameter. 
The existence  and symmetry of ground states in this case follows from the upper bound estimate of the least energy given by Lemma 4.3 below and the argument used in \cite{Ma-1}. See also \cite{VX-1, MV2017} and the references therein. 

\subsection{ The Case $\bf q\in (2,2+\frac{4\alpha}{N(2+\alpha)})$}  In this case, we set
$$
w(x)=\varepsilon^{-\frac{2N}{\alpha[4-N(q-2)]}}u(\varepsilon^{-\frac{N(q-2)}{\alpha[4-N(q-2)]}}x). 
\eqno(4.1)
$$
Then $(P_\varepsilon)$ is reduced to
$$
-\varepsilon^{-\sigma}\Delta w+w=(I_\alpha\ast |w|^{2_\alpha})|w|^{2_\alpha-2}w+\varepsilon^{-\sigma-\tau(q-1)}g(\varepsilon^\tau w),
\eqno(4.2)
$$
 where
$$
\sigma:=\frac{4\alpha-N(2+\alpha)(q-2)}{\alpha[4-N(q-2)]}>0,\quad 
\tau:=\frac{2N}{\alpha[4-N(q-2)]}.
$$ 
 The corresponding  functional is defined by 
$$
\begin{array}{rcl}
 J_\varepsilon(w):&=&\frac{\varepsilon^{-\sigma}}{2}\int_{\mathbb R^N}|\nabla w|^2+\frac{1}{2}\int_{\mathbb R^N}|w|^2-\frac{1}{22_\alpha}\int_{\mathbb R^N}(I_\alpha\ast |w|^{2_\alpha})|w|^{2_\alpha}\\
 &\mbox{}&-
\varepsilon^{-\sigma-\tau q}\int_{\mathbb R^N}G(\varepsilon^{\tau}w).
\end{array}
$$
Let $u_\varepsilon\in H^1(\mathbb R^N)$ be the ground state for $(P_\varepsilon)$ and 
$$
w_\varepsilon(x)=\varepsilon^{-\frac{2N}{\alpha[4-N(q-2)]}}u_\varepsilon(\varepsilon^{-\frac{N(q-2)}{\alpha[4-N(q-2)]}}x).
\eqno(4.3)
$$ 
Then  $w_\varepsilon$ is a ground state of (4.2).

To show that $w_\varepsilon$ is bounded in $H^1(\mathbb R^N)$, we also make the following 
$$
w_\varepsilon(x)=\varepsilon^{\frac{N[4\alpha-N(2+\alpha)(q-2)]}{4\alpha[4-N(q-2)]}}v_\varepsilon(\varepsilon^{\frac{4\alpha-N(2+\alpha)(q-2)}{2\alpha[4-N(q-2)]}}x).
\eqno(4.4)
$$
Then $v=v_\varepsilon$ is a ground state of  the following equation 
$$
-\Delta v+v=(I_\alpha\ast |v|^{2_\alpha})|v|^{2_\alpha-2}v+\varepsilon^{-\frac{4-N(q-2)}{4}\sigma}\varepsilon^{-\frac{N(2+\alpha)(q-1)}{4\alpha}}g(\varepsilon^{\frac{N(2+\alpha)}{4\alpha}}v).
\eqno(4.5)
$$
The corresponding  functional  of (4.5) is defined by
$$
\begin{array}{rcl}
\tilde J_\varepsilon(v):&=&\frac{1}{2}\int_{\mathbb R^N}|\nabla v|^2+|v|^2-\frac{1}{22_\alpha}\int_{\mathbb R^N}(I_\alpha\ast |v|^{2_\alpha})|v|^{2_\alpha}\\
&\mbox{}& -\varepsilon^{-\frac{4-N(q-2)}{4}\sigma}\varepsilon^{-\frac{N(2+\alpha)q}{4\alpha}}\int_{\mathbb R^N}G(\varepsilon^{\frac{N(2+\alpha)}{4\alpha}}v).
\end{array}
$$

Clearly, the following lemma holds true.

\noindent{\bf Lemma 4.1.}  {\it  Let $\varepsilon>0$, and $w_\varepsilon$, $v_\varepsilon$ are given in (4.3) and (4.4) respectively. Then

(1)  \ $ \ \|w_\varepsilon\|_2^2=\|v_\varepsilon\|_2^2, \quad \int_{\mathbb R^N}(I_\alpha\ast |w_\varepsilon|^{2_\alpha})|w_\varepsilon|^{2_\alpha}= \int_{\mathbb R^N}(I_\alpha\ast |v_\varepsilon|^{2_\alpha})|v_\varepsilon|^{2_\alpha},$

(2) \  $\varepsilon^{-\sigma}\|\nabla w_\varepsilon\|_2^2=\|\nabla v_\varepsilon\|_2^2, \quad   \varepsilon^{-\frac{N\sigma(q-2)}{4}}\|w_\varepsilon\|_q^q=\|v_\varepsilon\|_q^q$,

(3) \   $J_\varepsilon(w_\varepsilon)=\tilde J_\varepsilon(v_\varepsilon)$. }

\smallskip

As $\varepsilon\to \infty$, it follows from (H2) that the formal limit equation of (4.5) is as follows
$$
-\Delta v+v=(I_\alpha\ast |v|^{2_\alpha})|v|^{2_\alpha-2}v.
\eqno(4.6)
$$
Let 
$$
\tilde J_\infty(v):=\frac{1}{2}\int_{\mathbb R^N}|\nabla v|^2+|v|^2-\frac{1}{22_\alpha}\int_{\mathbb R^N}(I_\alpha\ast |v|^{2_\alpha})|v|^{2_\alpha},
$$
and 
$$
\tilde{\mathcal M}_\infty:=\left\{v\in H^1(\mathbb R^N)\setminus \{0\} \ \left | \  \int_{\mathbb R^N}|\nabla v|^2+|v|^2=\int_{\mathbb R^N}|v|^{2_\alpha}\right.\right\}.
$$
Then, from a minimax characterizations for $\tilde J_\infty$ similar to (4.9) below, it is easy to see    that
$$
\tilde J_\varepsilon(v_\varepsilon)\le \inf_{v\in H^1(\mathbb R^N)\setminus \{0\}}\sup_{t>0}\tilde J_\infty(tv)= \inf_{v\in \tilde{\mathcal M}_\infty}\tilde J_\infty(v)<\infty.
$$
Since $v_\varepsilon$ satisfies the Poho\v zaev identity, we obtain
$$
\begin{array}{rcl}
\frac{1}{2^*}\int_{\mathbb R^N}|\nabla v_\varepsilon|^2+\frac{1}{2}\int_{\mathbb R^N}|v_\varepsilon|^2&=&\frac{1}{2}\int_{\mathbb R^N}(I_\alpha\ast |v_\varepsilon|^{2_\alpha})|v_\varepsilon|^{2_\alpha}\\
&\mbox{}&+\varepsilon^{-\frac{4-N(q-2)}{4}\sigma}\varepsilon^{-\frac{N(2+\alpha)q}{4\alpha}}\int_{\mathbb R^N}G(\varepsilon^{\frac{N(2+\alpha)}{4\alpha}}v).
\end{array}
\eqno(4.7)
$$
Therefore, we have 
$$
\tilde J_\varepsilon(v_\varepsilon)=\frac{1}{N}\int_{\mathbb R^N}|\nabla v_\varepsilon|^2+\frac{2_\alpha-1}{22_\alpha}\int_{\mathbb R^N}(I_\alpha\ast |v_\varepsilon|^{2_\alpha})|v_\varepsilon|^{2_\alpha},
$$
and hence $\|\nabla v_\varepsilon\|_2$ and $\int_{\mathbb R^N}(I_\alpha\ast |v_\varepsilon|^{2_\alpha})|v_\varepsilon|^{2_\alpha}$ are bounded. 
Thus, by (4.7), (H2) and  the Gagliardo-Nirenberg inequality, we obtain
$$
\begin{array}{rcl}
\|v_\varepsilon\|_2^2&\le &\frac{2}{2^*}\|\nabla v_\varepsilon\|_2^2+\|v_\varepsilon\|_2^2\\
&=&\int_{\mathbb R^N}(I_\alpha\ast |v_\varepsilon|^{2_\alpha})|v_\varepsilon|^{2_\alpha}+2\varepsilon^{-\frac{4-N(q-2)}{4}\sigma}\varepsilon^{-\frac{N(2+\alpha)q}{4\alpha}}\int_{\mathbb R^N}G(\varepsilon^{\frac{N(2+\alpha)}{4\alpha}}v)\\
&\le &C+\varepsilon^{-\frac{4-N(q-2)}{4}\sigma}(\frac{A}{q}\int_{\mathbb R^N}|v_\varepsilon|^q+o_\varepsilon(1))\\
&\le &C+C\varepsilon^{-\frac{4-N(q-2)}{4}\sigma}\|\nabla v_\varepsilon\|_2^{\frac{N(q-2)}{2}}\|v_\varepsilon\|_2^{\frac{2N-q(N-2)}{2}},
\end{array}
$$
which together with the fact that $\|\nabla v_\varepsilon\|_2$ is bounded and $\frac{2N-q(N-2)}{2}<2$ yields $\|v_\varepsilon\|_2$ is bounded.
Thus, by Lemma 4.1, it follows that $w_\varepsilon$ is bounded in $L^2(\mathbb R^N)$.

It is also easy to see from Lemma 4.1 that  $I_\varepsilon(u)=\varepsilon^{\frac{N+\alpha}{\alpha}}J_\varepsilon(w)=\varepsilon^{\frac{N+\alpha}{\alpha}}\tilde J_\varepsilon(v)$ and  the energy of the ground states given by
$$
m_\varepsilon=\inf_{u\in \mathcal N_\varepsilon}I_\varepsilon(u)=\varepsilon^{\frac{N+\alpha}{\alpha}}\inf_{w\in \mathcal M_\varepsilon}J_\varepsilon(w)=\varepsilon^{\frac{N+\alpha}{\alpha}}\inf_{v\in \tilde{\mathcal M}_\varepsilon}\tilde J_\varepsilon(v)
$$
are well-defined and positive, where $\mathcal M_\varepsilon$ and $\tilde{\mathcal M}_\varepsilon$ denote the correspoding Nehari manifolds 
$$
\mathcal M_\varepsilon:=\left\{ w\in H^1(\mathbb R^N)\setminus\{0\}  \ \left | \ \begin{array}{rl}
&\int_{\mathbb R^N}\varepsilon^{-\sigma}|\nabla w|^2+|w|^2=\int_{\mathbb R^N}(I_\alpha\ast |w|^{2_\alpha})|w|^{2_\alpha}\\
&\mbox{}+\varepsilon^{-\sigma-\tau(q-1)}\int_{\mathbb R^N}g(\varepsilon^{\tau}w)w\end{array} \right. \right\},
$$
$$
\tilde{\mathcal M}_\varepsilon:=\left\{ v\in H^1(\mathbb R^N)\setminus\{0\}  \ \left | \ \begin{array}{rl}
&\int_{\mathbb R^N}|\nabla v|^2+|v|^2=\int_{\mathbb R^N}(I_\alpha\ast |v|^{2_\alpha})|v|^{2_\alpha}\\
&\mbox{}+\varepsilon^{-\frac{4-N(q-2)}{4}\sigma}\varepsilon^{-\frac{N(2+\alpha)(q-1)}{4\alpha}}\int_{\mathbb R^N}g(\varepsilon^{\frac{N(2+\alpha)}{4\alpha}}v)v \end{array}\right. \right\}.
$$

Define the Poho\v{z}aev manifold as follows
$$
\mathcal P_\varepsilon:=\{w\in H^1(\mathbb R^N)\setminus\{0\} \   | \  P_\varepsilon(w)=0 \},
$$
where 
$$
\begin{array}{rcl}
P_\varepsilon(w):&=&\frac{\varepsilon^{-\sigma}(N-2)}{2}\int_{\mathbb R^N}|\nabla w|^2+\frac{ N}{2}\int_{\mathbb R^N}|w|^2\\ \\
&\quad &-\frac{N+\alpha}{22_\alpha}\int_{\mathbb R^N}(I_\alpha\ast |w|^{2_\alpha})|w|^{2_\alpha}-N\varepsilon^{-\sigma -\tau q}\int_{\mathbb R^N}G(\varepsilon^\tau w).
\end{array}
\eqno(4.8)
$$
Then by Lemma 3.1, $w_\varepsilon\in \mathcal P_\varepsilon$. Moreover,  we have the following minimax characterizations for the least energy level $m_\varepsilon$.
$$
\varepsilon^{-\frac{N+\alpha}{\alpha}}m_\varepsilon=\tilde m_\varepsilon:=\inf_{w\in H^1(\mathbb R^N)\setminus\{0\}}\sup_{t\ge 0}J_\varepsilon(tw)=\inf_{w\in H^1(\mathbb R^N)\setminus\{0\}}\sup_{t\ge 0}J_\varepsilon(w_t).
\eqno(4.9)
$$
In particular, we have $\varepsilon^{-\frac{N+\alpha}{\alpha}}m_\varepsilon=\tilde m_\varepsilon=J_\varepsilon(w_\varepsilon)=\sup_{t>0}J_\lambda(tw_\varepsilon)=\sup_{t>0}J_\varepsilon((w_\varepsilon)_t)$. 

\smallskip

\noindent{\bf Lemma 4.2.} {\it The rescaled family of solutions $\{w_\varepsilon\}$ is bounded in $H^1(\mathbb R^N)$, and 
$$
\|\nabla w_\varepsilon\|_2^2=\frac{NA(q-2)}{2q}\|w_\varepsilon\|_q^q+o_\varepsilon(1),
\eqno(4.10)
$$
where and in what follows, $\lim_{\varepsilon\to\infty}o_\varepsilon(1)=0$.}
\begin{proof}   Since $\{w_\varepsilon\}$ is bounded in $L^2(\mathbb R^N)$, it suffices to show that $\{w_\varepsilon\}$ is also bounded 
in $D^{1,2}(\mathbb R^N)$.  By $w_\varepsilon\in \mathcal M_\varepsilon\cap\mathcal P_\varepsilon$,  we obtain
$$
\varepsilon^{-\sigma}\int_{\mathbb R^N}|\nabla w_\varepsilon|^2+\int_{\mathbb R^N}|w_\varepsilon|^2=\int_{\mathbb R^N}(I_\alpha\ast |w_\varepsilon|^{2_\alpha})|w_\varepsilon|^{2_\alpha}+\varepsilon^{-\sigma-\tau(q-1)}\int_{\mathbb R^N}g(\varepsilon^\tau w_\varepsilon)w_\varepsilon,
$$
$$
\frac{\varepsilon^{-\sigma}}{2^*}\int_{\mathbb R^N}|\nabla w_\varepsilon|^2+\frac{1}{2}\int_{\mathbb R^N}|w_\varepsilon|^2=\frac{N+\alpha}{2Np}\int_{\mathbb R^N}(I_\alpha\ast |w_\varepsilon|^{2_\alpha})|w_\varepsilon|^{2_\alpha}+\varepsilon^{-\sigma-\tau q}\int_{\mathbb R^N}G(\varepsilon^\tau w_\varepsilon).
$$
Therefore, we have 
$$
\varepsilon^{-\sigma}\|\nabla w_\varepsilon\|_2^2=N\varepsilon^{-\sigma-\tau q}\int_{\mathbb R^N}\left[\frac{1}{2}g(\varepsilon^\tau w_\varepsilon)\varepsilon^\tau w_\varepsilon-G(\varepsilon^\tau w_\varepsilon)\right].
\eqno(4.11)
$$
Set 
$g(s)=As^{q-1}+\tilde g(s)$ for $s\ge 0$, then $G(s)=\frac{A}{q}s^q+\tilde G(s)$, and 
$$
\lim_{s\to\infty}\tilde g(s)s^{1-q}=\lim_{s\to\infty}\tilde G(s)s^{-q}=0.
$$
Therefore, for any $\delta>0$ there exists a constant $C_\delta>0$ such that
$$
\left|\frac{1}{2}\tilde g(u)u-\tilde G(u)\right|\le \delta |u|^q+C_\delta |u|^2,
$$
and hence
$$
\begin{array}{rcl}
\varepsilon^{-\tau q}\int_{\mathbb R^N}\left|\frac{1}{2}\tilde g(\varepsilon^\tau w_\varepsilon)\varepsilon^\tau w_\varepsilon-\tilde G(\varepsilon^\tau w_\varepsilon)\right|
&\le &\delta \int_{\mathbb R^N}|w_\varepsilon|^q+C_\delta \varepsilon^{-\tau(q-2)}\int_{\mathbb R^N}|w_\varepsilon|^2\\
&\le &\delta \|w_\varepsilon\|_q^q+CC_\delta \varepsilon^{-\tau(q-2)}.
\end{array}
\eqno(4.12)
$$
By the Gagliardo-Nirenberg Inequality, we obtain 
$$
\|w_\varepsilon\|_q^q\le C\|\nabla w_\varepsilon\|_2^{\frac{N(q-2)}{2}}\|w_\varepsilon\|_2^{\frac{2N-q(N-2)}{2}}\le \tilde C\|\nabla w_\varepsilon\|_2^{\frac{N(q-2)}{2}}.
$$
Therefore, we get
$$
\begin{array}{rcl}
\|\nabla w_\varepsilon\|_2^2&=&N\varepsilon^{-\tau q}\int_{\mathbb R^N}\left[\frac{1}{2}g(\varepsilon^\tau w_\varepsilon)\varepsilon^\tau w_\varepsilon-G(\varepsilon^\tau w_\varepsilon)\right]\\
&\le &\left(\frac{NA(q-2)}{2q}+\delta\right)\int_{\mathbb R^N}|w_\varepsilon|^q+C\\
&\le &\tilde C\left(\frac{NA(q-2)}{2q}+\delta\right)\|\nabla w_\varepsilon\|_2^{\frac{N(q-2)}{2}}+C,
\end{array}
$$
which together with  $2>\frac{N(q-2)}{2}$ implies that $\{w_\varepsilon\}$ is bounded in $D^{1,2}(\mathbb R^N)$.
 
 By virtue of (4.12), for any $\delta>0$, we have 
$$
\limsup_{\varepsilon\to\infty} \varepsilon^{-\tau q}\int_{\mathbb R^N}\left|\frac{1}{2}\tilde g(\varepsilon^\tau w_\varepsilon)\varepsilon^\tau w_\varepsilon-\tilde G(\varepsilon^\tau w_\varepsilon)\right|
\le C\delta,
$$
and hence 
$$
\lim_{\varepsilon\to\infty} \varepsilon^{-\tau q}\int_{\mathbb R^N}\left[\frac{1}{2}\tilde g(\varepsilon^\tau w_\varepsilon)\varepsilon^\tau w_\varepsilon-\tilde G(\varepsilon^\tau w_\varepsilon)\right]=0.
\eqno(4.13)
$$
Clearly,  (4.11) implies that
$$
\|\nabla w_\varepsilon\|_2^2=\frac{NA(q-2)}{2q}\|w_\varepsilon\|_q^q+N\varepsilon^{-\tau q}\int_{\mathbb R^N}\left[\frac{1}{2}\tilde g(\varepsilon^\tau w_\varepsilon)\varepsilon^\tau w_\varepsilon-\tilde G(\varepsilon^\tau w_\varepsilon)\right],
$$
the equality (4.10)  then follows from (4.13).
The proof is complete.
\end{proof}

\smallskip

As $\varepsilon\to \infty$,  the limit equation of (4.2) is  the Hardy-Littlewood-Sobolev critical equation $w=(I_\alpha\ast |w|^{2_\alpha})|w|^{2_\alpha-2}w$, the corresponding  functional is defined by
$$
J_\infty(w):=\int_{\mathbb R^N}|w|^2-\frac{1}{22_\alpha}\int_{\mathbb R^N}(I_\alpha\ast |w|^{2_\alpha})|w|^{2_\alpha}.
$$
We define the Nehari manifolds as follows.
$$
\mathcal{M}_\infty=
\left\{w\in H^1(\mathbb R^N)\setminus\{0\} \ \left | \ \int_{\mathbb R^N}|w|^2=\int_{\mathbb R^N}(I_\alpha\ast |w|^{2_\alpha})|w|^{2_\alpha}\  \right. \right\}. 
$$
Then 
$$
\tilde m_\infty:=\inf_{w\in \mathcal {M}_\infty}J_\infty(w)
$$
is  well-defined and positive. Moreover,  as in (4.9), a similar minimax characterizations for the least energy level $\tilde m_\infty$ also holds, and 
$J_\infty$ is attained on $\mathcal M_\infty$ and 
$\tilde m_\infty=\frac{\alpha}{2(N+\alpha)}S_1^{\frac{N+\alpha}{\alpha}}$.

 For $w\in H^1(\mathbb R^N)\setminus \{0\}$, we set
$$
\tau_1(w)=\frac{\int_{\mathbb R^N}|w|^2}{\int_{\mathbb R^N}(I_\alpha\ast |w|^{2_\alpha})|w|^{2_\alpha}}.
\eqno(4.14)
$$
Then $(\tau_1(w))^{\frac{N}{2\alpha}}w\in \mathcal M_\infty$ for any $w\in H^1(\mathbb R^N)\setminus\{0\}$,  and $w\in \mathcal M_\infty$  if and only if $\tau_1(w)=1$.

\smallskip

Now, we give  the following estimation on $\tau_1(w_\varepsilon)$ and  the least energy $\tilde m_\varepsilon$:

\smallskip

\noindent{\bf Lemma 4.3.} {\it Let $N\ge 3$ and $q\in (2,2+\frac{4}{N})$. Then $\tau_1(w_\varepsilon)=1+O(\varepsilon^{-\sigma})$ and 
$$
\frac{\alpha}{2(N+\alpha)}S_1^{\frac{N+\alpha}{\alpha}}-\tilde m_\varepsilon=\tilde m_\infty-\tilde m_\varepsilon\sim \varepsilon^{-\sigma},
\eqno(4.15)
$$
 as $\varepsilon \to \infty$}.
\begin{proof} 
For any $\delta>0$ there exists a constant $C_\delta>0$ such that 
$$
g(s)s\le C_\delta |s|^q+\delta |s|^2, \quad for \ all \ s\ge 0.
$$
By  the Gagliardo-Nirenberg inequality and the fact that $\frac{N(q-2)}{2}<2$, it follows that 
$$
\begin{array}{cl}
&\frac{\varepsilon^{-\tau q}\int_{\mathbb R^N}g(\varepsilon^\tau w_\varepsilon)\varepsilon^\tau w_\varepsilon -\int_{\mathbb R^N}|\nabla w_\varepsilon|^2}{\int_{\mathbb R^N}|w_\varepsilon|^2}\\
&\le \frac{C_\delta \int_{\mathbb R^N}|w_\varepsilon|^q-\int_{\mathbb R^N}|\nabla w_\varepsilon|^2}{\int_{\mathbb R^N}|w_\varepsilon|^2}+\delta\varepsilon^{-\tau (q-2)}\\
&\le C_\delta\left(\frac{\|\nabla w_\varepsilon\|_2}{\|w_\varepsilon\|_2}\right)^{\frac{N(q-2)}{2}}\|w_\varepsilon\|_2^{q-2}-\left(\frac{\|\nabla w_\varepsilon\|_2}{\|w_\varepsilon\|_2}\right)^2
+\delta\varepsilon^{-\tau (q-2)}\\
&\le CC_\delta\left(\frac{\|\nabla w_\varepsilon\|_2}{\|w_\varepsilon\|_2}\right)^{\frac{N(q-2)}{2}}-\left(\frac{\|\nabla w_\varepsilon\|_2}{\|w_\varepsilon\|_2}\right)^2
+\delta\varepsilon^{-\tau (q-2)}\\
&\le C<\infty.
\end{array}
$$
First, since $v_\varepsilon\in \mathcal{M}_\varepsilon$,  by (4.14), it follows that 
$$
\begin{array}{rcl}
\tau_1(w_\varepsilon)&=&\frac{\int_{\mathbb R^N}|w_\varepsilon|^2}{\int_{\mathbb R^N}(I_\alpha\ast |w_\varepsilon|^p)|w_{\varepsilon}|^p}\\
&
=&\frac{\int_{\mathbb R^N}|w_\varepsilon|^2}{\int_{\mathbb R^N}|w_\varepsilon|^2-\varepsilon^{-\sigma}\left\{\varepsilon^{-\tau q}\int_{\mathbb R^N}g(\varepsilon^\tau w_\varepsilon)\varepsilon^\tau w_\varepsilon-\int_{\mathbb R^N} |\nabla w_\varepsilon|^2\right\}}\\
&\le &1+C\varepsilon^{-\sigma},
\end{array}
$$
from which  we get
$$
\int_{\mathbb R^N}|w_\varepsilon|^2=\tau_1(w_\varepsilon)\int_{\mathbb R^N}(I_\alpha\ast |w_\varepsilon|^{2_\alpha})|w_\varepsilon|^{2_\alpha}\le C\tau_1(w_\varepsilon)(\int_{\mathbb R^N}|w_\varepsilon|^2)^{2_\alpha},
$$
and hence
$$
\int_{\mathbb R^N}|w_\varepsilon|^2\ge (C\tau_1(w_\varepsilon))^{-\frac{N}{\alpha}}\ge C>0.
\eqno(4.16)
$$
By Lemma 4.1 and the Sobolev  inequality, we have
$$
\begin{array}{rcl}
\tau_1(w_\varepsilon)&=&\frac{\int_{\mathbb R^N}|w_\varepsilon|^2}{\int_{\mathbb R^N}(I_\alpha\ast |w_\varepsilon|^{2_\alpha})|w_{\varepsilon}|^{2_\alpha}}\\
&
=&\frac{\int_{\mathbb R^N}|w_\varepsilon|^2}{\int_{\mathbb R^N}|w_\varepsilon|^2-\varepsilon^{-\sigma}\left\{\varepsilon^{-\tau q}\int_{\mathbb R^N}g(\varepsilon^\tau w_\varepsilon)\varepsilon^\tau w_\varepsilon-\int_{\mathbb R^N} |\nabla w_\varepsilon|^2\right\}}\\
&
\ge &\frac{\int_{\mathbb R^N}|w_\varepsilon|^2}{\int_{\mathbb R^N}|w_\varepsilon|^2+\varepsilon^{-\sigma}\int_{\mathbb R^N} |\nabla w_\varepsilon|^2}\\
&\ge  &1-\frac{1}{2}\varepsilon^{-\sigma}\frac{\|\nabla w_\varepsilon\|_2^2}{\|w_\varepsilon\|_2^2}\\
&\ge &1-C\varepsilon^{-\sigma},
\end{array}
$$
thus we conclude that $\tau_1(w_\varepsilon)=1+O(\varepsilon^{-\sigma})$ as $\varepsilon\to \infty$.

For $w\in H^1(\mathbb R^N)$, let $w_t$ be defined as in Lemma 3.2. 
Then by Lemma 3.2 and Poho\v{z}aev's identity, it is easy to show that $\sup_{t\ge 0}J_\varepsilon((w_\varepsilon)_t)=J_\varepsilon(w_\varepsilon)=\tilde m_\varepsilon$. Therefore, we get 
$$
\begin{array}{rcl}
\tilde m_\infty&\le& \sup_{t\ge 0}J_\infty((w_\varepsilon)_t)=J_\infty((w_\varepsilon)_t)|_{t=\tau_1((w_\varepsilon))^{1/\alpha}}\\
&\le & \sup_{t\ge 0}J_\varepsilon((w_\varepsilon)_t)+\varepsilon^{-\sigma-\tau q}(\tau_1(w_\varepsilon))^{N/\alpha}\int_{\mathbb R^N}G(\varepsilon^\tau w_\varepsilon)\\
&\le& \tilde m_\varepsilon+\varepsilon^{-\sigma}(1+O(\varepsilon^{-\sigma}))^{N/\alpha}(\frac{A}{q}\|v_\varepsilon\|_q^q+o_\varepsilon(1))\\
&=& \tilde m_\varepsilon+O(\varepsilon^{-\sigma}).
\end{array}
$$

For each $\rho>0$, the family  $U_\rho(x):=\rho^{-\frac{N}{2}}U_1(x/\rho)$ are radial ground states of  $ v=(I_\alpha\ast |v|^{2_\alpha})v^{2_\alpha-1}$, and verify that
 $$
 \|\nabla U_\rho\|_2^2=\rho^{-2}\|\nabla U_1\|_2^2, \qquad \int_{\mathbb R^N}|U_\rho|^q=\rho^{N-\frac{N}{2}q}\int_{\mathbb R^N}|U_1|^q.
 \eqno(4.17)
 $$ 
 Let $g_0(\rho)=\frac{A}{q}\int_{\mathbb R^N}|U_\rho|^q-\frac{1}{2}\int_{\mathbb R^N}|\nabla U_\rho|^2$. Then
 $$
\begin{array}{rcl}
 g_0(\rho)&=&\frac{A}{q}\rho^{N-\frac{N}{2}q}\int_{\mathbb R^N}|U_1|^q-\frac{1}{2}\rho^{-2}\int_{\mathbb R^N}|\nabla U_1|^2\\
& =&\rho^{-2}\left\{\frac{A}{q}\rho^{\frac{4-N(q-2)}{2}}\int_{\mathbb R^N}|U_1|^q-\frac{1}{2}\int_{\mathbb R^N}|\nabla U_1|^2\right\}.
\end{array}
 $$
Since $q\in (2,2+\frac{4}{N})$, we have $\lim_{\rho\to 0}g_0(\rho)=-\infty, \ \lim_{\rho\to +\infty}g_0(\rho)=0$.  
 It is easy to see that  there exists $\rho_0=\rho(q)\in (0,+\infty)$ with
 $$
 \rho_0=\left(\frac{2q\int_{\mathbb R^N}|\nabla U_1|^2}{NA(q-2)\int_{\mathbb R^N}|U_1|^q}\right)^{\frac{2}{4-N(q-2)}}
 $$
such that  $g'_0(\rho_0)=0$ and hence
 $$
 g_0(\rho_0)=\sup_{\rho>0}g_0(\rho)=\frac{4-N(q-2)}{2N(q-2)}\left(\frac{NA(q-2)\int_{\mathbb R^N}|U_1|^q}{2q\int_{\mathbb R^N}|\nabla U_1|^2}\right)^{\frac{4}{4-N(q-2)}}\int_{\mathbb R^N}|\nabla U_1|^2.
 $$ 
 Let $U_0=U_{\rho_0}$, then there exists $t_\varepsilon\in (0,+\infty)$ such that
 $$
\begin{array}{rcl}
\tilde m_\varepsilon&\le &\sup_{t\ge 0}J_\varepsilon((U_0)_t)=J_\varepsilon( (U_0)_{t_\varepsilon})\\
&\le &\sup_{t>0}\left\{\frac{t^N}{2}\int_{\mathbb R^N}| U_0|^2-\frac{t^{N+\alpha}}{22_\alpha}\int_{\R^N}(I_\alpha\ast |U_0|^{2_\alpha})|U_0|^{2_\alpha}\right\}\\
&\mbox{}& -\varepsilon^{-\sigma}\left\{\varepsilon^{-\tau q}t^N_\varepsilon\int_{\mathbb R^N}G(\varepsilon^\tau U_0)-\frac{1}{2}t^{N-2}_\varepsilon\int_{\mathbb R^N}|\nabla U_0|^2\right\}\\
&\le & \tilde m_\infty-\varepsilon^{-\sigma}\left\{t^N_\varepsilon\int_{\mathbb R^N}\left[\frac{A}{q}|U_0|^q+\varepsilon^{-\tau q}\tilde G(\varepsilon^\tau U_0)\right]-\frac{1}{2}t^{N-2}_\varepsilon\int_{\mathbb R^N}|\nabla U_0|^2\right\}.
\end{array}
\eqno(4.18)
$$
where $t_\varepsilon\in (0,+\infty)$ satisfies
$$
\begin{array}{rcl}
\frac{1}{2^*}t_\varepsilon^{N-2}\varepsilon^{-\sigma}\int_{\mathbb R^N}|\nabla U_0|^2+\frac{1}{2}t_\varepsilon^N\int_{\mathbb R^N}|U_0|^2
&=&\frac{1}{2}t_\varepsilon^{N+\alpha}\int_{\mathbb R^N}(I_\alpha\ast |U_0|^{2_\alpha})|U_0|^{2_\alpha}\\
&\mbox{}&\quad +t_\varepsilon^N\varepsilon^{-\sigma-\tau q}\int_{\mathbb R^N}G(\varepsilon^\tau U_0).
\end{array}
$$

If $t_\varepsilon\ge 1$, then
$$
\frac{1}{2^*}\varepsilon^{-\sigma}\int_{\mathbb R^N}|\nabla U_0|^2+\frac{1}{2}\int_{\mathbb R^N}|U_0|^2\ge \frac{1}{2}t_\varepsilon^\alpha \int_{\mathbb R^N}(I_\alpha\ast |U_0|^{2_\alpha})|U_0|^{2_\alpha}+\varepsilon^{-\sigma-\tau q}\int_{\mathbb R^N}G(\varepsilon^\tau U_0).
$$
Hence
$$
t_\varepsilon^\alpha\le \frac{\frac{1}{2^*}\varepsilon^{-\sigma}\int_{\mathbb R^N}|\nabla U_0|^2+\frac{1}{2}\int_{\mathbb R^N}| U_0|^2-\varepsilon^{-\sigma-\tau q}\int_{\mathbb R^N}G(\varepsilon^\tau U_0)}{\frac{1}{2}\int_{\mathbb R^N}(I_\alpha\ast |U_0|^{2_\alpha})|U_0|^{2_\alpha}},
$$
which implies that $\limsup_{\varepsilon\to \infty}t_\varepsilon^{\alpha}\le 1.$

If $t_\varepsilon\le  1$, then a similar argument yields that
$$
t_\varepsilon^\alpha\ge \frac{\frac{1}{2^*}\varepsilon^{-\sigma}\int_{\mathbb R^N}|\nabla U_0|^2+\frac{1}{2}\int_{\mathbb R^N}| U_0|^2-\varepsilon^{-\sigma-\tau q}\int_{\mathbb R^N}G(\varepsilon^\tau U_0)}{\frac{1}{2}\int_{\mathbb R^N}(I_\alpha\ast |U_0|^{2_\alpha})|U_0|^{2_\alpha}},
$$
which implies that $\liminf_{\varepsilon\to \infty}t_\varepsilon^{\alpha}\ge 1.$

Therefore, $\lim_{\lambda\to 0}t_\lambda=1$ and hence
$$
\begin{array}{rl}
&\lim_{\varepsilon\to \infty}\left\{t^N_\varepsilon\int_{\mathbb R^N}\left[\frac{A}{q}|U_0|^q+\varepsilon^{-\tau q}\tilde G(\varepsilon^\tau U_0)\right]-\frac{1}{2}t^{N-2}_\varepsilon\int_{\mathbb R^N}|\nabla U_0|^2\right\}\\
&\quad =\frac{A}{q}\int_{\mathbb R^N}|U_0|^q-\frac{1}{2}\int_{\mathbb R^N}|\nabla U_0|^2=g_0(\rho_0).
\end{array}
$$
Thus, there exists a constant $C>0$ such that
$$
\tilde m_\varepsilon\le \tilde m_\infty-C\varepsilon^{-\sigma},
$$
for large $\varepsilon>0$. The proof is complete. 
\end{proof}

Let 
$$
\mathbb D(u):=\int_{\mathbb R^N}(I_\alpha\ast |u|^{\frac{N+\alpha}{N}})|u|^{\frac{N+\alpha}{N}}.
$$
The following result is  a special case of the classical Brezis-Lieb lemma \cite{Brezis-1} for Riesz potentials, for a proof, we refer the reader to \cite[Lemma 2.4]{Moroz-2}.

\smallskip

\noindent{\bf Lemma 4.4.} {\it Let $N \in\mathbb  N, \alpha\in  (0,N)$, and $(w_n)_{n\in \mathbb N}$ be a bounded sequence in $L^2(\mathbb R^N)$. If $w_n \to  w$ almost everywhere on $\mathbb R^N$ as $n\to \infty$, then
$$
\lim_{n\to\infty}\mathbb{D}(w_n)-\mathbb D(w_n-w_0) = \mathbb D(w_0).
$$}

\noindent{\bf Lemma 4.5.} {\it  $\|w_\varepsilon\|_2^2\sim \mathbb D(w_\varepsilon)\sim \|\nabla w_\varepsilon \|_2^2\sim \|w_\varepsilon\|_q^q\sim 1$ as $\varepsilon\to \infty$.}
\begin{proof}  
By the definition of $\tau_1(w_\varepsilon)$, Lemma 4.3  and the  Hardy-Sobolev-Littlewood inequality, for large $\varepsilon>0$, we have 
$$
\|w_\varepsilon\|_2^2=\tau_1(w_\varepsilon)\int_{\mathbb R^N}(I_\alpha\ast |w_\varepsilon|^{2_\alpha})|w_\varepsilon|^{2_\alpha}\le 2S_1^{-2_\alpha}\|w_\varepsilon\|_2^{22_\alpha},
$$
and thus it follows that 
$$
\|w_\varepsilon\|_2^2\ge 2^{-\frac{N}{\alpha}}S_1^{\frac{N+\alpha}{\alpha}},
$$
which together with the boundedness of $w_\varepsilon$ implies that $\|w_\varepsilon\|_2^2\sim1$ as $\varepsilon\to \infty$.

It   follows from Lemma 4.3  that
$$
\begin{array}{rcl}
\tilde m_\infty&\le &J_\infty((\tau_1(w_\varepsilon))^{\frac{N}{2\alpha}}w_\varepsilon)\\
&=&\frac{1}{2}(\tau_1(w_\varepsilon))^{\frac{N}{\alpha}}\|w_\varepsilon\|_2^2-\frac{1}{22_\alpha}(\tau_1(w_\varepsilon))^{\frac{N+\alpha}{\alpha}}\int_{\mathbb R^N}(I_\alpha\ast |w_\varepsilon|^{2_\alpha})|w_\varepsilon|^{2_\alpha}\\
&\le& (\tau_1(w_\varepsilon))^{\frac{N}{\alpha}}\left[\frac{1}{2}\|w_\varepsilon\|_2^2-\frac{1}{22_\alpha}\int_{\mathbb R^N}(I_\alpha\ast |w_\varepsilon|^{2_\alpha})|w_\varepsilon|^{2_\alpha}\right].
\end{array}
\eqno(4.19)
$$
Since $w_\varepsilon\in \mathcal N_\varepsilon$,  by Lemma 4.4, we obtain
$$
\begin{array}{rcl}
\tau_1(w_\varepsilon)&=&\frac{\int_{\mathbb R^N}|w_\varepsilon|^2}{\int_{\mathbb R^N}(I_\alpha\ast |w_\varepsilon|^{2_\alpha})|w_\varepsilon|^{2_\alpha}}\\
&=&\frac{\int_{\mathbb R^N}|w_\varepsilon|^2}{\int_{\mathbb R^N}|w_\varepsilon|^2+\varepsilon^{-\sigma}(\int_{\mathbb R^N}|\nabla w_\varepsilon|^2-\varepsilon^{-\tau q}\int_{\mathbb R^N}g(\varepsilon^\tau w_\varepsilon)\varepsilon^\tau w_\varepsilon)}\\
&=&\frac{\int_{\mathbb R^N}|w_\varepsilon|^2}{\int_{\mathbb R^N}|w_\varepsilon|^2-\varepsilon^{-\sigma}\frac{2N-q(N-2)}{2q}(\int_{\mathbb R^N}|w_\varepsilon|^q+o_\varepsilon(1))}\\
&\le & 1+C_1\varepsilon^{-\sigma} \|w_\varepsilon\|_q^q.
\end{array}
\eqno(4.20)
$$
Therefore, by (4.10), (4.19) and  (4.20),  we obtain 
$$
\begin{array}{rcl}
\tilde m_\varepsilon&=&\frac{1}{2}\int_{\mathbb R^N}\varepsilon^{-\sigma}|\nabla w_\varepsilon|_2^2+|w_\varepsilon|_2^2-\frac{1}{22_\alpha}\int_{\mathbb R^N}(I_\alpha\ast |w_\varepsilon|^{2_\alpha})|w_\varepsilon|^{2_\alpha}-\varepsilon^{-\sigma-\tau q}\int_{\mathbb R^N}G(\varepsilon^\tau w_\varepsilon)\\
&\ge& \varepsilon^{-\sigma}(\frac{1}{2}\|\nabla w_\varepsilon\|_2^2-\frac{A}{q}\|w_\varepsilon\|_q^q+o_\varepsilon(1))+\tilde m_\infty(\tau_1(w_\varepsilon))^{-\frac{N}{\alpha}}\\
&\ge &\varepsilon^{-\sigma}\frac{[N(q-2)-4]A}{4q}(\|w_\varepsilon\|_q^q+o_\varepsilon(1))+\tilde m_\infty-C_1\varepsilon^{-\sigma}\|w_\varepsilon\|_q^q.
\end{array}
$$
By Lemma 4.3, we have $\tilde m_\infty-\tilde m_\varepsilon\ge C_2\varepsilon^{-\sigma}$. Hence,  we get
$$
\frac{[4-N(q-2)]A}{4q}\|w_\varepsilon\|_q^q\ge C_2-C_1\|w_\varepsilon\|_q^q, 
$$
which yields 
$$
\|w_\varepsilon\|_q^q\ge \frac{4qC_2}{[4-N(q-2)]A+4qC_1}>0.
$$
Since $w_\varepsilon$ is bounded in $H^1(\mathbb R^N)$, it follows that $\|w_\varepsilon\|_q^q\sim 1$ as $\varepsilon\to \infty$.

Since $\|w_\varepsilon\|_2^2\sim 1$ as $\varepsilon\to \infty$,  the Gagliardo-Nirenberg inequality implies 
$$
\|w_\varepsilon\|_q^q\le C\|\nabla w_\varepsilon\|_2^{\frac{N(q-2)}{2}},
$$
which together with the boundedness of $w_\varepsilon$ in $H^1(\mathbb R^N)$  yields $\|\nabla w_\varepsilon \|_2^2\sim 1$ as $\varepsilon\to \infty$.

Finally, by the definition of $\tau_1(w_\varepsilon)$, Lemma 4.3 and Lemma 4.4,  it follows  that 
$$
\mathbb D(w_\varepsilon)=\int_{\mathbb R^N}(I_\alpha\ast |w_\varepsilon|^{\frac{N+\alpha}{N}})|w_\varepsilon|^{\frac{N+\alpha}{N}}=(\tau_1(w_\varepsilon))^{-1} \|w_\varepsilon\|_2^2\sim 1,
$$
as $\varepsilon\to \infty$. The proof is complete. 
\end{proof}

\noindent{\bf Lemma 4.6.} {\it Let $N\ge 3$ and $q\in (2, 2+\frac{4}{N})$,  then for any $\varepsilon_n\to \infty$, there exists $\rho\in [\rho_0, +\infty)$, such that, up to a subsequence, $w_{\varepsilon_n}\to U_{\rho}$ in $L^2(\mathbb R^N)$, where 
$$
  \rho_0=\rho_0(q):= \left(\frac{2q\int_{\mathbb R^N}|\nabla U_1|^2}{NA(q-2)\int_{\mathbb R^N}|U_1|^q}\right)^{\frac{2}{4-N(q-2)}}.
 $$
Moreover, $w_{\varepsilon_n}\to U_\rho$ in $D^{1,2}(\mathbb R^N)$ if and only if $\rho=\rho_0$. }
\begin{proof} The proof is similar to that in \cite[Lemma 4.8]{Ma-2}, and for the reader's convenience, we give the detail proof. Note that $w_{\varepsilon_n}$ is a positive radially symmetric function, and by Lemma 4.2, $\{w_{\varepsilon_n}\}$ is bounded in $H^1(\mathbb R^N)$. Then there exists $w_\infty\in H^1(\mathbb R^N)$ such that 
$$
w_{\varepsilon_n} \rightharpoonup w_\infty   \quad {\rm weakly \ in} \  H^1(\mathbb R^N), \quad w_{\varepsilon_n}\to w_\infty \quad {\rm in} \ L^p(\mathbb R^N) \quad {\rm for \ any} \ p\in (2,2^*),
\eqno(4.21)
$$
and 
$$
w_{\varepsilon_n}(x)\to w_\infty(x) \quad a. \ e. \  {\rm on} \ \mathbb R^N,  \qquad w_{\varepsilon_n}\to w_\infty \quad {\rm in} \   L^2_{loc}(\mathbb R^N).
\eqno(4.22)
$$
Observe that
$$
J_\infty(w_{\varepsilon_n})=J_{\varepsilon_n}(w_{\varepsilon_n})+\varepsilon_n^{-\sigma-\tau q}\int_{\mathbb R^N}G(\varepsilon^\tau w_{\varepsilon_n})-\frac{\varepsilon_n^{-\sigma}}{2}\int_{\mathbb R^N}|\nabla w_{\varepsilon_n}|^2=\tilde m_{\varepsilon_n}+o_n(1)=\tilde m_\infty+o_n(1),
$$
and 
$$
J'_\infty(w_{\varepsilon_n})w=J'_{\varepsilon_n}(w_{\varepsilon_n})w+\varepsilon_n^{-\sigma-\tau(q-1)}\int_{\mathbb R^N}g(\varepsilon^\tau w_{\varepsilon_n}) w-\varepsilon_n^{-\sigma}\int_{\mathbb R^N}\nabla w_{\varepsilon_n} \nabla w=o_n(1).
$$
Therefore, $\{w_{\varepsilon_n}\}$ is a PS sequence of $J_\infty$ at level $\tilde m_\infty=\frac{\alpha}{2(N+\alpha)}S_1^{\frac{N+\alpha}{\alpha}}$.

By Lemma 4.5, we have $w_\infty\not=0$.
By Lemma 4.3 and Lemma 4.4, we have
$$
\begin{array}{rcl}
o_n(1)&=&\tilde m_{\varepsilon_n}-\tilde m_\infty\\
&\ge &\frac{\varepsilon_n^{-\sigma}}{2}\|\nabla w_{\varepsilon_n}\|_2^2+\frac{1}{2}[\|w_{\varepsilon_n}\|_2^2-\|w_0\|_2^2]-\frac{N}{2(N+\alpha)}[\mathbb D(w_{\varepsilon_n})-\mathbb D(w_\infty)]\\
&\mbox{}&-\varepsilon_n^{-\sigma-\tau q}\int_{\mathbb R^N}G(\varepsilon_n^\tau w_{\varepsilon_n})\\
&=&\frac{1}{2}\|w_{\varepsilon_n}-w_\infty\|_2^2-\frac{N}{2(N+\alpha)}\mathbb D(w_{\varepsilon_n}-w_\infty)+o_n(1),
\end{array}
$$
$$
\begin{array}{rcl}
0&=&\langle J_{\varepsilon_n}'(w_{\varepsilon_n})-J_\infty'(w_\infty), w_{\varepsilon_n}-w_\infty\rangle\\
&=& \varepsilon_n^{-\sigma}\int_{\mathbb R^N}\nabla w_{\varepsilon_n}(\nabla w_{\varepsilon_n}-\nabla w_\infty)+\int_{\mathbb R^N}|w_{\varepsilon_n}-w_\infty|^2\\
&\mbox{}&\quad 
-\int_{\mathbb R^N}(I_\alpha\ast |w_{\varepsilon_n}|^{\frac{N+\alpha}{N}})w_{\varepsilon_n}^{\frac{\alpha}{N}}(w_{\varepsilon_n}-w_\infty)
+\int_{\mathbb R^N}(I_\alpha\ast |w_\infty|^{\frac{N+\alpha}{N}})w_\infty^{\frac{\alpha}{N}}(w_{\varepsilon_n}-w_\infty)\\
&\mbox{}&\quad
-\varepsilon_n^\sigma\int_{\mathbb R^N}|w_{\varepsilon_n}|^{q-2}w_{\varepsilon_n}(w_{\varepsilon_n}-w_\infty)\\
&=& \|w_{\varepsilon_n}-w_0\|_2^2-\mathbb D(w_{\varepsilon_n}-w_\infty)+o_n(1).
\end{array}
$$
Hence, it follows that
$$
\|w_{\varepsilon_n}-w_\infty\|_2^2\le \frac{N}{N+\alpha}\mathbb D(w_{\varepsilon_n}-w_\infty)+o_n(1)=\frac{N}{N+\alpha}\|w_{\varepsilon_n}-w_\infty\|_2^2+o_n(1),
$$
and hence
$$
\|w_{\varepsilon_n}-w_\infty\|_2\to 0, \qquad  {\rm as} \  \  \varepsilon_n\to \infty.
$$
By the Hardy-Littlewood-Sobolev inequality and Lemma 4.4, it follows that 
$$
\lim_{\varepsilon_n\to \infty}\mathbb D(w_{\varepsilon_n})=\mathbb D(w_\infty).
$$
Since $\tau_1(w_{\varepsilon_n})\to 1$ as $\varepsilon_n\to \infty$ by Lemma 4.3,  it follows that $\tau_1(w_\infty)=1$ and hence $w_\infty\in \mathcal M_\infty$.

On the other hand, by Lemma 4.1 and the boundedness of $v_{\varepsilon_n}$ in $H^1(\mathbb R^N)$, we have 
$$
\begin{array}{rcl}
\tilde m_{\varepsilon_n}&=&J_{\varepsilon_n}(w_{\varepsilon_n})\\
&=&\frac{1}{2}\int_{\mathbb R^N}\varepsilon_n^{-\sigma}|\nabla w_{\varepsilon_n}|^2+|w_{\varepsilon_n}|^2-\frac{1}{2p}\int_{\mathbb R^N}
(I_\alpha\ast |w_{\varepsilon_n}|^{\frac{N+\alpha}{N}})|w_{\varepsilon_n}|^{\frac{N+\alpha}{N}}\\
&\mbox{}& -\varepsilon_n^{-\sigma-\tau q}\int_{\mathbb R^N}G(\varepsilon^\tau w_{\varepsilon_n})\\
&\ge &  \frac{1}{2}\int_{\mathbb R^N}|w_{\varepsilon_n}|^2-\frac{1}{2p}\int_{\mathbb R^N}
(I_\alpha\ast |w_{\varepsilon_n}|^{\frac{N+\alpha}{N}})|w_{\varepsilon_n}|^{\frac{N+\alpha}{N}}\\
&\mbox{}& -\varepsilon_n^{-\sigma}\{C_\delta\int_{\mathbb R^N}|w_{\varepsilon_n}|^q+\delta\varepsilon_n^{-\tau(q-2)}\int_{\mathbb R^N}|w_{\varepsilon_n}|^2\}.
\end{array}
$$
Sending $\varepsilon_n\to \infty$, it then follows from Lemma 4.3 that
$$
\tilde m_\infty\ge \frac{1}{2}\int_{\mathbb R^N}|w_\infty|^2-\frac{1}{2p}\int_{\mathbb R^N}
(I_\alpha\ast |w_\infty|^{\frac{N+\alpha}{N}})|w_\infty|^{\frac{N+\alpha}{N}}
=J_\infty(w_\infty).
$$
Therefore, note that $w_\infty\in \mathcal M_\infty$, we obtain $J_\infty(w_\infty)=\tilde m_\infty$. 
Thus, $w_\infty=U_{\rho}$ for some $\rho\in (0,+\infty)$.

 Moreover, by (4.10), we obtain
$$
\|\nabla w_\infty\|_2^2\le \lim_{\varepsilon_n\to \infty}\|\nabla w_{\varepsilon_n}\|_2^2=\frac{NA(q-2)}{2q}\int_{\mathbb R^N}|w_\infty|^q,
$$
from which it follows that  
 $$
 \rho\ge\left(\frac{2q\int_{\mathbb R^N}|\nabla U_1|^2}{NA(q-2)\int_{\mathbb R^N}|U_1|^q}\right)^{\frac{2}{4-N(q-2)}}.
 $$
If $\rho=\rho_0$, then (4.10) implies that $\lim_{n\to\infty}\|\nabla w_{\varepsilon_n}\|_2^2=\|\nabla U_{\rho_0}\|_2^2$, and hence $w_{\varepsilon_n}\to U_{\rho_0}$ in $D^{1,2}(\mathbb R^N)$.
\end{proof}

\vskip 3mm


\begin{proof}[Proof of Theorem 2.1 in the case $ q\in (2,2+\frac{4\alpha}{N(2+\alpha)})$]   
Let 
$$
M_\varepsilon=w_\varepsilon(0),  \quad z_\varepsilon=M_\varepsilon [U_{\rho_0}(0)]^{-1},  
$$
where $\rho_0$ is given in Lemma 4.6.  Then
$$
z_\varepsilon=\frac{1}{U_{\rho_0}(0)}w_\varepsilon(0)=\frac{1}{U_{\rho_0}(0)}\varepsilon^{-\frac{2N}{\alpha[4-N(q-2)]}}u_\varepsilon(0).
$$
Hence
$$
\zeta_\varepsilon :=z_\varepsilon^{-\frac{2}{N}} \varepsilon^{-\frac{N(q-2)}{\alpha[4-N(q-2)]}}=\left(\frac{U_{\rho_0}(0)}{u_\varepsilon(0)}\right)^{\frac{2}{N}}\varepsilon^{\frac{1}{\alpha}},
\quad 
z_\varepsilon^{-1}\varepsilon^{-\frac{2N}{\alpha[4-N(q-2)]}}=\varepsilon^{-\frac{N}{2\alpha}}\zeta_\varepsilon^{\frac{N}{2}}.
$$
Arguing as in the proof of \cite[Theorem 2.1]{Ma-2}, we can show that 
$$
\zeta_\varepsilon \sim \varepsilon^{-\frac{N(q-2)}{\alpha[4-N(q-2)]}}
$$
and  for large $\varepsilon>0$, the rescaled family of ground states
$$
\tilde w_\varepsilon(x)=z_\varepsilon^{-1}w_\varepsilon(z_\varepsilon^{-\frac{2}{N}} x)=z_\varepsilon^{-1}\varepsilon^{-\frac{2N}{\alpha[4-N(q-2)]}}u_\varepsilon(z_\varepsilon^{-\frac{2}{N}}\varepsilon^{-\frac{N(q-2)}{\alpha[4-N(q-2)]}} x)=\varepsilon^{-\frac{N}{2\alpha}}\zeta_\varepsilon^{\frac{N}{2}}u_\varepsilon(\zeta_\varepsilon x)
$$
satisfies 
$$
\|\nabla \tilde w_\varepsilon\|^2_2\sim \|\tilde w_\varepsilon\|_{q}^{q} \sim \int_{\mathbb R^N}(I_\alpha \ast |\tilde w_\varepsilon|^{\frac{N+\alpha}{N}})|\tilde w_\varepsilon|^{\frac{N+\alpha}{N}}\sim \|\tilde w_\varepsilon\|_2^2\sim 1,
$$
and as $\varepsilon\to\infty$, $\tilde w_\varepsilon$ converges in $L^2(\mathbb R^N)$ to the extremal function  $U_{\rho_0}$. 
Then by Lemma 4.6, we also have $\tilde w_\varepsilon\to U_{\rho_0}$ in $D^{1,2}(\mathbb R^N)$. Thus we conclude that $\tilde w_\varepsilon\to U_{\rho_0}$ in $H^1(\mathbb R^N)$. 
 
 Since $w_\varepsilon\in \mathcal M_\varepsilon$ and $w_\varepsilon$ is bounded in $H^1(\mathbb R^N)$,  it follows that
$$
\begin{array}{lcl}
\tilde m_\varepsilon&=&\frac{2_\alpha-1}{22_\alpha}\varepsilon^{-\sigma}\int_{\R^N}|\nabla  w_\varepsilon|^2+\frac{2_\alpha-1}{22_\alpha}\int_{\R^N}|w_\varepsilon|^2-\varepsilon^{-\sigma-\tau q}\int_{\R^N}[G(\varepsilon^\tau w_\varepsilon)-g(\varepsilon^\tau w_\varepsilon)\varepsilon^\tau w_\varepsilon]\\
&=&\frac{\alpha}{2(N+\alpha)}\int_{\R^N}|w_\varepsilon|^2+O(\varepsilon^{-\sigma}).
\end{array}
$$
Similarly, we also have 
$$
\tilde m_\infty=\frac{\alpha}{2(N+\alpha)}\int_{\R^N}|U_1|^2.
$$
Then it follows  from Lemma 4.3 that 
$$
\int_{\R^N}| U_1|^2-\int_{\R^N}|w_\varepsilon|^2=\frac{2(N+\alpha)}{\alpha}(\tilde m_\infty-\tilde m_\varepsilon)+O(\varepsilon^{-\sigma})=O(\varepsilon^{-\sigma}).
$$
Since 
$
\| U_1\|_2^2=\int_{\R^N}(I_\alpha\ast |U_1|^p)|U_1|^p=S_1^{\frac{N+\alpha}{\alpha}},
$
we conclude  that 
$$
\|\tilde w_\varepsilon\|_2^2=\| w_\varepsilon\|_2^2=S_1^{\frac{N+\alpha}{\alpha}}+O(\varepsilon^{-\sigma}).
$$
Finally, by $w_\varepsilon\in \mathcal M_\varepsilon$ and the boundedness of $w_\varepsilon\in H^1(\mathbb R^N)$, we also have
$$
\int_{\R^N}(I_\alpha\ast |\tilde w_\varepsilon|^{\frac{N+\alpha}{N}})|\tilde w_\varepsilon|^{\frac{N+\alpha}{N}}=\int_{\R^N}(I_\alpha\ast |w_\varepsilon|^{\frac{N+\alpha}{N}})|w_\varepsilon|^{\frac{N+\alpha}{N}}=\|w_\varepsilon\|_2^2+O(\varepsilon^{-\sigma})=  S_1^{\frac{N+\alpha}{\alpha}}+O(\varepsilon^{-\sigma}).
$$
The statements on $u_\varepsilon$ follow from the corresponding results on $v_\varepsilon$ and $\tilde w_\varepsilon$, and 
the proof is complete.  \end{proof}

\subsection{ The Case $\bf q\in (2+\frac{4\alpha}{N(2+\alpha)}, 2+\frac{4}{N})$}   In this case, we set
$$
v(x)=\varepsilon^{-\frac{1}{q-2}}u(\varepsilon^{-\frac{1}{2}}x).
\eqno(4.23)
$$
 Then $(P_\varepsilon)$ is reduced to
$$
-\Delta v+v=\varepsilon^{-\frac{N(2+\alpha)(q-2)-4\alpha}{2N(q-2)}}(I_\alpha\ast |v|^{2_\alpha})|v|^{2_\alpha-2}v+\varepsilon^{\frac{1-q}{q-2}}g(\varepsilon^{\frac{1}{q-2}} v).
\eqno(4.24)
$$
Let $u_\varepsilon\in H^1(\mathbb R^N)$ be the ground state for $(P_\varepsilon)$ and 
$$
v_\varepsilon(x)=\varepsilon^{-\frac{1}{q-2}}u_\varepsilon(\varepsilon^{-\frac{1}{2}}x).
$$ 
Then  $v_\varepsilon$ is a ground state of  (4.24).  Moreover, we have 
$$
\|u_\varepsilon\|_2^2=\varepsilon^{\frac{4-N(q-2)}{2(q-2)}}\|v_\varepsilon\|_2^2,\quad 
\|\nabla u_\varepsilon\|_2^2=\varepsilon^{\frac{2N-q(N-2)}{2(q-2)}}\|\nabla v_\varepsilon\|_2^2,\quad 
\|u_\varepsilon\|_q^q=\varepsilon^{\frac{2N-q(N-2)}{2(q-2)}}\|v_\varepsilon\|_q^q.
\eqno(4.25)
$$
The limiting equation of (4.24) reads as follows
$$
-\Delta v+v=A|v|^{q-2}v.
\eqno(4.26)
$$
Associated to (4.24) and (4.26), the corresponding  functionals are  defined by 
$$
\begin{array}{rcl}
 J_\varepsilon(v):&=&\frac{1}{2}\int_{\mathbb R^N}|\nabla v|^2+|v|^2-\frac{1}{22_\alpha}\varepsilon^{-\frac{N(2+\alpha)(q-2)-4\alpha}{2N(q-2)}}\int_{\mathbb R^N}(I_\alpha\ast |v|^{2_\alpha})|v|^{2_\alpha}\\
 &\mbox{}&\qquad \qquad \qquad \qquad -
\varepsilon^{-\frac{q}{q-2}}\int_{\mathbb R^N}G(\varepsilon^{\frac{1}{q-2}}v),
\end{array}
$$
and 
$$
J_\infty(v):=\frac{1}{2}\int_{\mathbb R^N}|\nabla v|^2+|v|^2-\frac{A}{q}\int_{\mathbb R^N}|v|^q,
$$
respectively. The associated  Nehari manifolds are as follows
$$
 \mathcal N_\varepsilon:=\left\{v\in H^1(\mathbb R^N) \ \left| \begin{array}{rcl}
 \ \int_{\mathbb R^N}|\nabla v_\varepsilon|^2+|v_\varepsilon|^2&=&\varepsilon^{-\eta_q}\int_{\mathbb R^N}(I_\alpha\ast |v_\varepsilon|^{2_\alpha})|\varepsilon|^{2_\alpha}\\
 &\mbox{}&+\varepsilon^{\frac{1-q}{q-2}}\int_{\mathbb R^N}g(\varepsilon^{\frac{1}{q-2}}v_\varepsilon)v_\varepsilon\end{array} \right.\right\},
$$
$$
 \mathcal N_\infty:=\left\{v\in H^1(\mathbb R^N) \ \left| \ \int_{\mathbb R^N}|\nabla v_\varepsilon|^2+|v_\varepsilon|^2=A\int_{\mathbb R^N}|v|^q\right.\right\},
 $$
 respectively.
Denote the lest energies by $\tilde m_\varepsilon$ and $\tilde m_\infty$, then
$$
\tilde m_\varepsilon=\inf_{v\in \mathcal N_\varepsilon}J_\varepsilon(v), \qquad \tilde m_\infty=\inf_{v\in \mathcal N_\infty}J_\infty(v)
$$
are well defined and positive. Arguing as in Lemma 4.3, we show that $\tilde m_\varepsilon= \tilde m_\infty+o_\varepsilon(1)$. 

Let $\eta_q:=\frac{N(2+\alpha)(q-2)-4\alpha}{2N(q-2)}$, then
$$
\tilde m_\varepsilon=\frac{1}{2}\int_{\mathbb R^N}|\nabla v_\varepsilon|^2+|v_\varepsilon|^2-\frac{1}{22_\alpha}\varepsilon^{-\eta_q}\int_{\mathbb R^N}(I_\alpha\ast |v_\varepsilon|^{2_\alpha})|v_\varepsilon|^{2_\alpha}-
\varepsilon^{-\frac{q}{q-2}}\int_{\mathbb R^N}G(\varepsilon^{\frac{1}{q-2}}v_\varepsilon),
\eqno(4.27)
$$
$$
\int_{\mathbb R^N}|\nabla v_\varepsilon|^2+|v_\varepsilon|^2=\varepsilon^{-\eta_q}\int_{\mathbb R^N}(I_\alpha\ast |v_\varepsilon|^{2_\alpha})|v_\varepsilon|^{2_\alpha}+\varepsilon^{\frac{1-q}{q-2}}\int_{\mathbb R^N}g(\varepsilon^{\frac{1}{q-2}}v_\varepsilon)v_\varepsilon,
\eqno(4.28)
$$
and 
$$
\frac{1}{2^*}\int_{\mathbb R^N}|\nabla v_\varepsilon|^2+\frac{1}{2}\int_{\mathbb R^N}|v_\varepsilon|^2=\frac{1}{2}\varepsilon^{-\eta_q}\int_{\mathbb R^N}(I_\alpha\ast |v_\varepsilon|^{2_\alpha})|v_\varepsilon|^{2_\alpha}+\varepsilon^{-\frac{q}{q-2}}\int_{\mathbb R^N}G(\varepsilon^{\frac{1}{q-2}}v_\varepsilon).
\eqno(4.29)
$$
By (4.27) and (4.29), we get
$$
\tilde m_\infty+o_\varepsilon(1)\ge \tilde m_\varepsilon=\frac{1}{N}\int_{\mathbb R^N}|\nabla v_\varepsilon|^2+\frac{2_\alpha-1}{22_\alpha}\varepsilon^{-\eta_q}\int_{\mathbb R^N}(I_\alpha\ast |v_\varepsilon|^{2_\alpha})|v_\varepsilon|^{2_\alpha},
$$
from which, we conclude  that $v_\varepsilon$ is bounded in $D^{1,2}(\mathbb R^N)$, $\varepsilon^{-\eta_q}\int_{\mathbb R^N}(I_\alpha\ast |v_\varepsilon|^{2_\alpha})|v_\varepsilon|^{2_\alpha}$ is bounded. Then by (4.28), it is easy to show that $v_\varepsilon$ is also bounded in $L^2(\mathbb R^N)$. Therefore, $v_\varepsilon$ is bounded in $H^1(\mathbb R^N)$, and hence, by the Hardy-Litlewood-Sobolev indetity, $\varepsilon^{-\eta_q}\int_{\mathbb R^N}(I_\alpha\ast |v_\varepsilon|^{2_\alpha})|v_\varepsilon|^{2_\alpha}=O(\varepsilon^{-\eta_q})$ and  there holds
$$
\tilde m_\varepsilon=\frac{1}{2}\int_{\mathbb R^N}|\nabla v_\varepsilon|^2+|v_\varepsilon|^2
-\frac{A}{q}\int_{\mathbb R^N}|v_\varepsilon|^q+o_\varepsilon(1),
$$
$$
\int_{\mathbb R^N}|\nabla v_\varepsilon|^2+|v_\varepsilon|^2=
A\int_{\mathbb R^N}|v_\varepsilon|^q+o_\varepsilon(1),
$$
and
$$
\frac{1}{2^*}\int_{\mathbb R^N}|\nabla v_\varepsilon|^2+\frac{1}{2}\int_{\mathbb R^N}|v_\varepsilon|^2=
\frac{A}{q}\int_{\mathbb R^N}|v_\varepsilon|^q+o_\varepsilon(1).
$$
Therefore, we obtain
$$
\int_{\mathbb R^N}|\nabla v_\varepsilon|^2=\frac{NA(q-2)}{2q}\int_{\mathbb R^N}|v_\varepsilon|^q+o_\varepsilon(1).
$$
Put
$$
A_\varepsilon:=\int_{\mathbb R^N}|\nabla v_\varepsilon|^2, \quad  B_\varepsilon:=\int_{\mathbb R^N}|v_\varepsilon|^2, \quad C_\varepsilon:=\int_{\mathbb R^N}|v_\varepsilon|^q,
$$
then
$$
\left\{\begin{array}{rcl} \frac{1}{2}A_\varepsilon+\frac{1}{2}B_\varepsilon-\frac{A}{q}C_\varepsilon&=&\tilde m_\varepsilon+o_\varepsilon(1),\\
A_\varepsilon+B_\varepsilon-AC_\varepsilon&=&o_\varepsilon(1),\\
\frac{1}{2^*}A_\varepsilon+\frac{1}{2}B_\varepsilon-\frac{A}{q}C_\varepsilon&=&o_\varepsilon(1).\end{array}\right.
$$
Solving this system, we obtain
$$
A_\varepsilon=N\tilde m_\varepsilon+o_\varepsilon(1), \quad B_\varepsilon=\frac{2N-q(N-2)}{q-2}\tilde m_\varepsilon+o_\varepsilon(1),
\quad C_\varepsilon=\frac{2q}{A(q-2)}\tilde m_\varepsilon +o_\varepsilon(1).
$$
It is easy to show that
$$
\tilde m_\varepsilon=\tilde m_\infty+o_\varepsilon(1)=\frac{q-2}{2q}A^{-\frac{2}{q-2}}S_q^{\frac{q}{q-2}}+o_\varepsilon(1),
$$
$$
 m_\varepsilon=\varepsilon^{\frac{2N-q(N-2)}{2(q-2)}}\tilde m_\varepsilon
=\varepsilon^{\frac{2N-q(N-2)}{2(q-2)}}\left(\frac{q-2}{2q}A^{-\frac{2}{q-2}}S_q^{\frac{q}{q-2}}+o_\varepsilon(1)\right).
$$
Therefore, we conclude that
$$
A_\varepsilon=N\tilde m_\varepsilon+o_\varepsilon(1)=\frac{N(q-2)}{2q}A^{-\frac{2}{q-2}}S_q^{\frac{q}{q-2}}+o_\varepsilon(1),
$$
$$
 \quad B_\varepsilon=\frac{2N-q(N-2)}{q-2}\tilde m_\varepsilon+o_\varepsilon(1)=\frac{2N-q(N-2)}{2q}A^{-\frac{2}{q-2}}S_q^{\frac{q}{q-2}}+o_\varepsilon(1),
 $$
 $$
\quad C_\varepsilon=\frac{2q}{A(q-2)}\tilde m_\varepsilon +o_\varepsilon(1)=A^{-\frac{q}{q-2}}S_q^{\frac{q}{q-2}}+o_\varepsilon(1).
$$
Form which, we conclude the proof of Theorem 2.1 in the case $q\in (2+\frac{4\alpha}{N(2+\alpha)}, 2+\frac{4}{N})$.
\vskip 3mm 

\subsection{ The Case $\bf q=2+\frac{4\alpha}{N(2+\alpha)}$}   In this case, we set
$$
v(x)=\varepsilon^{-\frac{1}{q-2}}u(\varepsilon^{-\frac{1}{2}}x).
\eqno(4.30)
$$
 Then $(P_\varepsilon)$ is reduced to
$$
-\Delta v+v=(I_\alpha\ast |v|^{2_\alpha})|v|^{2_\alpha-2}v+\varepsilon^{\frac{1-q}{q-2}}g(\varepsilon^{\frac{1}{q-2}} v).
\eqno(4.31)
$$
Let $u_\varepsilon\in H^1(\mathbb R^N)$ be the ground state for $(P_\varepsilon)$ and 
$$
v_\varepsilon(x)=\varepsilon^{-\frac{1}{q-2}}u_\varepsilon(\varepsilon^{-\frac{1}{2}}x).
$$ 
Then  $v_\varepsilon$ is a ground state of (4.26).  Moreover, we have 
$$
\|u_\varepsilon\|_2^2=\varepsilon^{\frac{4-N(q-2)}{2(q-2)}}\|v_\varepsilon\|_2^2,\quad 
\|\nabla u_\varepsilon\|_2^2=\varepsilon^{\frac{2N-q(N-2)}{2(q-2)}}\|\nabla v_\varepsilon\|_2^2,\quad 
\|u_\varepsilon\|_q^q=\varepsilon^{\frac{2N-q(N-2)}{2(q-2)}}\|v_\varepsilon\|_q^q.
\eqno(4.32)
$$
The limiting equation of (4.31) reads as follows
$$
-\Delta v+v=(I_\alpha\ast |v|^{2_\alpha})|v|^{2_\alpha-2}v+A|v|^{q-2}v.
\eqno(4.33)
$$
The rest of proof is similar to the case $q\in  (2+\frac{4\alpha}{N(2+\alpha)}, 2+\frac{4}{N})$ and is omitted.

\smallskip

\vskip 5mm

\section{Proof of Theorem 2.2}

 In this section, we always assume that $2_\alpha^*:=\frac{N+\alpha}{N-2}$ and  $q\in (2, 2^*)$ if $N\ge 4$, and $q\in (4,6)$ if $N=3$. 
 For large $\varepsilon>0$, the existence  and symmetry of ground states in this case follows from the upper bound estimates of the least energy given by Lemma 5.3  and Lemma 5.4 below and the argument used in \cite{Ma-1}. See also \cite{MV2017} and references therein.

 \subsection{The first scaling}
  It is easy to see that under the rescaling 
$$
w(x)=\varepsilon^{-\frac{1}{q-2}}u(\varepsilon^{-\frac{2}{(N-2)(q-2)}}x), 
\eqno(5.1)
$$
the equation $(P_\varepsilon)$ is reduced to 
$$
-\Delta w+ \varepsilon^{-\sigma}w=(I_\alpha\ast |w|^{2_\alpha^*})|w|^{2_\alpha^*-2}w+\varepsilon^{-\sigma-\tau(q-1)} g(\varepsilon^\tau w), 
\eqno(5.2)
$$
where $\sigma:=\frac{2^*-q}{q-2}=\frac{2N-q(N-2)}{(N-2)(q-2)}>0$ and $\tau:=\frac{1}{q-2}.$ 

The associated functional is defined by 
$$
J_\varepsilon(w):=\frac{1}{2}\int_{\mathbb R^N}|\nabla w|^2+\varepsilon^{-\sigma}|w|^2-\frac{1}{22_\alpha^*}\int_{\mathbb R^N}(I_\alpha\ast |w|^{2_\alpha^*})|w|^{2_\alpha^*}-\varepsilon^{-\sigma-\tau q}\int_{\mathbb R^N}G(\varepsilon^\tau w).
$$

\noindent{\bf Lemma  5.1.}  {\it Let $\varepsilon>0, u\in H^1(\mathbb R^N)$ and $w$ is the rescaling  of $u$ given in (5.1).  Then

(1) $ \    \   \|\nabla w\|_2^2=\|\nabla u\|_2^2, \quad \int_{\mathbb R^N}(I_\alpha\ast |w|^{2_\alpha^*})|w|^{2_\alpha^*}= \int_{\mathbb R^N}(I_\alpha\ast |u|^{2_\alpha^*})|u|^{2_\alpha^*},$

(2)  $ \   \   \varepsilon^{-\sigma}\|w\|_2^2=\varepsilon \|u\|_2^2,  \quad   \varepsilon^{-\sigma}\|w\|_q^q=\|u\|_q^q$,

(3)  $\    \  I_\varepsilon(u)=J_\varepsilon(w)$.}

\smallskip

We define the Nehari manifolds as follows.
$$
\mathcal{N}_\varepsilon=
\left\{w\in H^1(\mathbb R^N)\setminus\{0\} \ \left | \ \begin{array}{rl}&\int_{\mathbb R^N}|\nabla w|^2+\varepsilon^{-\sigma}\int_{\mathbb R^N}|w|^2\\
&=\int_{\mathbb R^N}(I_\alpha\ast |w|^{2_\alpha^*})|w|^{2_\alpha^*}+\varepsilon^{-\sigma-\tau(q-1)}
\int_{\mathbb R^N}g(\varepsilon^\tau w)w \end{array} \right. \right\}
$$
and 
$$
\mathcal{N}_\infty=
\left\{w\in H^1(\mathbb R^N)\setminus\{0\} \ \left | \ \int_{\mathbb R^N}|\nabla w|^2=\int_{\mathbb R^N}(I_\alpha\ast |w|^{2_\alpha^*})|w|^{2_\alpha^*}\  \right. \right\}. 
$$
Then 
$$
m_\varepsilon:=\inf_{w\in \mathcal {N}_\varepsilon}J_\varepsilon(w) \quad {\rm and} \quad 
m_\infty:=\inf_{w\in \mathcal {N}_\infty}J_\infty(w)
$$
are well-defined and positive. Moreover, $J_\infty$ is attained on $\mathcal N_\infty$.

 For $w\in H^1(\mathbb R^N)\setminus \{0\}$, we set
$$
\tau_2(w)=\frac{\int_{\mathbb R^N}|\nabla w|^2}{\int_{\mathbb R^N}(I_\alpha\ast |w|^{2_\alpha^*})|w|^{2_\alpha^*}}.
\eqno(5.3)
$$
Then $(\tau_2(w))^{\frac{N-2}{2(2+\alpha)}}w\in \mathcal N_\infty$ for any $w\in H^1(\mathbb R^N)\setminus\{0\}$,  and $w\in \mathcal N_\infty$  if and only if $\tau_2(w)=1$.

Define the Poho\v zaev manifold as follows
$$
\mathcal P_\varepsilon:=\{w\in H^1(\mathbb R^N)\setminus\{0\} \   | \  P_\varepsilon(w)=0 \},
$$
where 
$$
\begin{array}{rcl}
P_\varepsilon(w):&=&\frac{N-2}{2}\int_{\mathbb R^N}|\nabla w|^2+\frac{ \varepsilon^{-\sigma} N}{2}\int_{\mathbb R^N}|w|^2\\ \\
&\quad &-\frac{N+\alpha}{22_\alpha^*}\int_{\mathbb R^N}(I_\alpha\ast |w|^{2_\alpha^*})|w|^{2_\alpha^*}-N\varepsilon^{-\sigma-\tau q}\int_{\mathbb R^N}G(\varepsilon^\tau w).
\end{array}
\eqno(5.4)
$$
Then by Lemma 3.1, $w_\varepsilon\in \mathcal P_\varepsilon$. Moreover,  we have a similar minimax characterizations for the least energy level $m_\varepsilon$ as in Lemma 3.2.
A similar result also holds for $m_\infty$ and $J_\infty$.

\smallskip

\noindent{\bf Lemma 5.2.} {\it The rescaled family of solutions $\{w_\varepsilon\}$ is bounded in $H^1(\mathbb R^N)$, and 
$$
\|w_\varepsilon\|_2^2=\frac{NA(2^*-q)}{2^*q}\|w_\varepsilon\|_q^q+o_\varepsilon(1), \quad as \  \ \varepsilon\to\infty.
\eqno(5.5)
$$}

\begin{proof}  To prove $w_\varepsilon$ is bounded in $D^{1,2}(\mathbb R^N)$, we let
$$
v(x)=\varepsilon^{-\frac{N-2}{4}}u(\varepsilon^{-\frac{1}{2}}x).
$$
Then $(P_\varepsilon)$ reduces to 
$$
-\Delta v+v=(I_\alpha\ast |v|^{2_\alpha^*})|v|^{2_\alpha^*-2}v+\varepsilon^{-\frac{N+2}{4}}g(\varepsilon^{\frac{N-2}{4}}v),
$$
and $v_\varepsilon(x)=\varepsilon^{-\frac{N-2}{4}}u_\varepsilon(\varepsilon^{-\frac{1}{2}}x)$  is a grond state solution, and hence
$$
\frac{1}{2^*}\int_{\mathbb R^N}|\nabla v_\varepsilon|^2+\frac{1}{2}\int_{\mathbb R^N}|v_\varepsilon|^2=\frac{N+\alpha}{22_\alpha^*N}\int_{\mathbb R^N}(I_\alpha\ast |v_\varepsilon|^{2_\alpha^*})|v_\varepsilon|^{2_\alpha^*}+\varepsilon^{-\frac{N}{2}}\int_{\mathbb R^N}G(\varepsilon^{\frac{N-2}{4}}v_\varepsilon),
$$
and 
$$
\begin{array}{rcl}
m_\infty\ge m_\varepsilon&=&\frac{1}{2}\int_{\mathbb R^N}|\nabla v_\varepsilon|^2+|v_\varepsilon|^2-\frac{1}{22_\alpha^*}\int_{\mathbb R^N}(I_\alpha\ast |v_\varepsilon|^{2_\alpha^*})|v_\varepsilon|^{2_\alpha^*}-\varepsilon^{-\frac{N}{2}}\int_{\mathbb R^N}G(\varepsilon^{\frac{N-2}{4}}v_\varepsilon)\\
&=&\frac{1}{N}\int_{\mathbb R^N}|\nabla v_\varepsilon|^2+\frac{\alpha}{22_\alpha^*N}\int_{\mathbb R^N}(I_\alpha\ast |v_\varepsilon|^{2_\alpha^*})|v_\varepsilon|^{2_\alpha^*}.
\end{array}
$$
Therefore, $v_\varepsilon$ is bounded in $D^{1,2}(\mathbb R^N)$. Clearly, we have $\|\nabla v_\varepsilon\|_2=\|\nabla u_\varepsilon\|_2=\|\nabla w_\varepsilon\|_2$, it follows that $\{w_\varepsilon\}$ is bounded in $D^{1,2}(\mathbb R^N)$. Thus, it suffices to show that it is also bounded 
in $L^2(\mathbb R^N)$.  By $w_\varepsilon\in \mathcal N_\varepsilon\cap \mathcal P_\varepsilon$, we obtain
$$
\int_{\mathbb R^N}|\nabla w_\varepsilon|^2+\varepsilon^{-\sigma} \int_{\mathbb R^N}|w_\varepsilon|^2=\int_{\mathbb R^N}(I_\alpha\ast |w_\varepsilon|^{2_\alpha^*})|w_\varepsilon|^{2_\alpha^*}+\varepsilon^{-\sigma-\tau(q-1)}\int_{\mathbb R^N}g(\varepsilon^\tau w_\varepsilon)w_\varepsilon,
$$
$$
\frac{N-2}{2}\int_{\mathbb R^N}|\nabla w_\varepsilon|^2+\frac{\varepsilon^{-\sigma} N}{2}\int_{\mathbb R^N}|w_\varepsilon|^2=\frac{N+\alpha}{22_\alpha^*}\int_{\mathbb R^N}(I_\alpha\ast |w_\varepsilon|^{2_\alpha^*})|w_\varepsilon|^{2_\alpha^*}+N\varepsilon^{-\sigma-\tau q}\int_{\mathbb R^N}G(\varepsilon^\tau w_\varepsilon).
$$
Let $g(s)=As^{q-1}+\tilde g(s)$ for $s\ge 0$, then $G(s)=\frac{A}{q}s^q+\tilde G(s)$, then
$$
\frac{1}{N}\int_{\mathbb R^N}|w_\varepsilon|^2=\frac{A(2^*-q)}{2^*q}\int_{\mathbb R^N}|w_\varepsilon|^q+\varepsilon^{-\tau q}\int_{\mathbb R^N}[\tilde G(\varepsilon^\tau w_\varepsilon)-\frac{1}{2^*}\tilde g(\varepsilon^\tau w_\varepsilon)\varepsilon^\tau w_\varepsilon].
\eqno(5.6)
$$
By (H1) and (H2), for any $\delta>0$ there exists a constant $C_\delta>0$ such that
$$
|\tilde G(s)-\frac{1}{2^*}\tilde g(s)s|\le \delta s^{q}+C_\delta s^2, \quad {\rm for \ all} \ s\ge 0.
$$
Therefore, we get 
$$
\varepsilon^{-\tau q}\int_{\mathbb R^N}|\tilde G(\varepsilon^\tau w_\varepsilon)-\frac{1}{2^*}\tilde g(\varepsilon^\tau w_\varepsilon)\varepsilon^\tau w_\varepsilon|
\le \delta \int_{\mathbb R^N}|w_\varepsilon|^q+C_\delta\varepsilon^{-\tau(q-2)} \int_{\mathbb R^N}|w_\varepsilon|^2. 
\eqno(5.7)
$$
By (5.5), we have
$$
[\frac{1}{N}-C_\delta \varepsilon^{-\tau(q-2)}]\int_{\mathbb R^N}|w_\varepsilon|^2\le [\frac{A(2^*-q)}{2^*q}+\delta]\int_{\mathbb R^N}|w_\varepsilon|^q.
\eqno(5.8)
$$
By the Sobolev embedding theorem  and the  interpolation inequality, we obtain 
$$
\begin{array}{rcl}
\int_{\mathbb R^N}|w_\varepsilon|^q&\le& \left(\int_{\mathbb R^N}|w_\varepsilon|^2\right)^{\frac{2^*-q}{2^*-2}}\left(\int_{\mathbb R^N}|w_\varepsilon|^{2^*}\right)^{\frac{q-2}{2^*-2}}\\
&\le& \left(\int_{\mathbb R^N}|w_\varepsilon|^2\right)^{\frac{2^*-q}{2^*-2}}\left(\frac{1}{S}\int_{\mathbb R^N}|\nabla w_\varepsilon|^2\right)^{\frac{2^*(q-2)}{2(2^*-2)}},
\end{array}
$$
where $S$ is the best Sobolev constant.  It then follows from (5.7) and Lemma 5.1 that
$$
\frac{1}{2N}\left(\int_{\mathbb R^N}|w_\varepsilon|^2\right)^{\frac{q-2}{2^*-2}}\le [\frac{A(2^*-q)}{2^*q}+\delta]\left(\frac{1}{S}\int_{\mathbb R^N}|\nabla w_\varepsilon|^2\right)^{\frac{2^*(q-2)}{2(2^*-2)}},
$$
which together with the boundedness of $w_\varepsilon$ in $D^{1,2}(\mathbb R^N)$ implies that $w_\varepsilon$ is bounded in $L^2(\mathbb R^N)$.  

Finally, it follows from (5.7) that
$$
\lim_{\varepsilon\to \infty}\varepsilon^{-\tau q}\int_{\mathbb R^N}[\tilde G(\varepsilon^\tau w_\varepsilon)-\frac{1}{2^*}\tilde g(\varepsilon^\tau w_\varepsilon)\varepsilon^\tau w_\varepsilon]=0,
$$
this together with (5.6) yields (5.5), and completes the proof.
\end{proof}

Now, we give  the following estimation on the least energy:

\smallskip

\noindent{\bf Lemma 5.3.} {\it  If $N\ge 5$ and $q\in (2,2^*)$, then
  $$
  m_\infty-m_\varepsilon \sim \varepsilon^{-\sigma}\qquad  as \quad \varepsilon\to \infty.
  $$
  If $N=4$ and $q\in (2,4)$, or $N=3$ and $q\in (4,6)$, then
  $$
  m_\infty-m_\varepsilon\lesssim \varepsilon^{-\sigma},  \qquad as \quad \varepsilon\to \infty.
  $$
  }
\begin{proof}  First, we claim that  
$$
\tau_2(w_\varepsilon)=1+O(\varepsilon^{-\sigma}),\quad as \  \ \varepsilon\to \infty.
\eqno(5.9)
$$
In fact, since  $w_\varepsilon\in \mathcal{N}_\varepsilon$,  we see that
$$
\tau_2(w_\varepsilon)=\frac{\int_{\mathbb R^N}|\nabla w_\varepsilon|^2}{\int_{\mathbb R^N}(I_\alpha \ast |w_\varepsilon|^{2_\alpha^*})|w_\varepsilon|^{2_\alpha^*}}
=1+\varepsilon^{-\sigma} \frac{\varepsilon^{-\tau(q-1)}\int_{\mathbb R^N}g(\varepsilon^\tau w_\varepsilon)w_\varepsilon-\int_{\mathbb R^N}|w_\varepsilon|^2}{\int_{\mathbb R^N}(I_\alpha\ast |w_\varepsilon|^{2_\alpha^*})|w_\varepsilon|^{2_\alpha^*}}.
$$
Since 
$$
\int_{\mathbb R^N}|w_\varepsilon|^q\le \left(\int_{\mathbb R^N}|w_\varepsilon|^2\right)^{\frac{2^*-q}{2^*-2}}\left(\int_{\mathbb R^N}|w_\varepsilon|^{2^*}\right)^{\frac{q-2}{2^*-2}}, 
$$
we see that 
$$
\begin{array}{rcl}
\frac{\varepsilon^{-\tau(q-1)}\int_{\mathbb R^N}g(\varepsilon^\tau w_\varepsilon)w_\varepsilon-\int_{\mathbb R^N}|w_\varepsilon|^2}{\int_{\mathbb R^N}|w_\varepsilon|^{2^*}}
&\le&\frac{C_\delta\int_{\mathbb R^N}|w_\varepsilon|^q-(1-\delta\varepsilon^{-\tau(q-2)})\int_{\mathbb R^N}|w_\varepsilon|^2}{\int_{\mathbb R^N}|w_\varepsilon|^{2^*}}\\
&\le &\zeta_\varepsilon^{\theta_q}(C_\delta-(1-\delta\varepsilon^{-\tau(q-2)})\zeta_\varepsilon^{1-\theta_q}).
\end{array}
 $$
where 
$$
\theta_q=\frac{2^*-q}{2^*-2}, \qquad \zeta_\varepsilon=\frac{\int_{\mathbb R^N}|w_\varepsilon|^2}{\int_{\mathbb R^N}|w_\varepsilon|^{2^*}}.
$$
Let $
g(x):=x^{\theta_q}(C_\delta-(1-\delta\varepsilon^{-\tau(q-2)})x^{1-\theta_q}).$ Then  
$$
\sup_{x\ge 0}g(x)=K_q:=C_\delta\theta_q(1-\theta_q)\left(\frac{C_\delta\theta_q}{1-\delta\varepsilon^{-\tau(q-2)}}\right)^{\frac{\theta_q}{1-\theta_q}}.
$$
Therefore, by the boundedness of $w_\varepsilon$ in $D^{1,2}(\R^N)$, we get
$$
\begin{array}{rcl}
\tau_2(w_\varepsilon)&\le& 1+\varepsilon^{-\sigma}K_q\frac{\int_{\mathbb R^N}|w_\varepsilon|^{2^*}}{\int_{\mathbb R^N}(I_\alpha\ast |w_\varepsilon|^{2_\alpha^*})|w_\varepsilon|^{2_\alpha^*}}\\
&\le & 1+\varepsilon^{-\sigma}K_qS^{-\frac{N}{N-2}}\frac{(\int_{\mathbb R^N}|\nabla w_\varepsilon|^{2})^{\frac{N}{N-2}}}{\int_{\mathbb R^N}(I_\alpha\ast |w_\varepsilon|^{2_\alpha^*})|w_\varepsilon|^{2_\alpha^*}}\\
&=&1+\varepsilon^{-\sigma }K_qS^{-\frac{N}{N-2}}\tau_2(w_\varepsilon)(\int_{\mathbb R^N}|\nabla w_\varepsilon|^{2})^{\frac{2}{N-2}}\\
&\le & 1+\varepsilon^{-\sigma} C\tau_2(w_\varepsilon),
\end{array}
$$
and hence for large $\varepsilon>0$, there holds
$$
\tau_2(w_\varepsilon)\le \frac{1}{1-\varepsilon^{-\sigma} C}=1+\varepsilon^{-\sigma}\frac{C}{1-\varepsilon^{-\sigma} C}\le 1+\frac{1}{2}C\varepsilon^{-\sigma}.
$$
On the other hand,  by (5.3), the Hardy-Littlewood-Sobolev and the Sobolev inequalities, we obtain that 
$$
\int_{\mathbb R^N}|\nabla w_\varepsilon|^2=\tau_2(w_\varepsilon)\int_{\mathbb R^N}(I_\alpha\ast |w_\varepsilon|^{2_\alpha^*})|w_\varepsilon|^{2_\alpha^*}
\le C\tau_2(w_\varepsilon)(\frac{1}{S}\int_{\mathbb R^N}|\nabla w_\varepsilon|^2)^{2_\alpha^*}.
$$
Therefore, we get
$$
\int_{\mathbb R^N}|\nabla w_\varepsilon|^2\ge (\frac{C\tau_2(w_\varepsilon)}{S^{2_\alpha^*}})^{-\frac{1}{2_\alpha^*-1}}\ge C>0.
$$
By (5.3), $w_\varepsilon\in \mathcal N_\varepsilon$ and the boundedness of $w_\varepsilon$ in $H^1(\mathbb R^N)$, we then  obtain
$$
\begin{array}{rcl}
\tau_2(w_\varepsilon)&=&\frac{\int_{\mathbb R^N}|\nabla w_\varepsilon|^2}{\int_{\mathbb R^N}(I_\alpha \ast |w_\varepsilon|^{2_\alpha^*})|w_\varepsilon|^{2_\alpha^*}}\\
&=&\frac{\int_{\mathbb R^N}|\nabla w_\varepsilon|^2}{\int_{\mathbb R^N}|\nabla w_\varepsilon|^2+\varepsilon^{-\sigma}\int_{\mathbb R^N}|w_\varepsilon|^2-\varepsilon^{-\sigma-\tau(q-1)}\int_{\mathbb R^N}g(\varepsilon^\tau w_\varepsilon)w_\varepsilon}\\
&\ge &\frac{\int_{\mathbb R^N}|\nabla w_\varepsilon|^2}{\int_{\mathbb R^N}|\nabla w_\varepsilon|^2+\varepsilon^{-\sigma}\int_{\mathbb R^N}|w_\varepsilon|^2}\\
&\ge &1-\frac{1}{2}\varepsilon^{-\sigma} \frac{\int_{\mathbb R^N}|w_\varepsilon|^2}{\int_{\mathbb R^N}|\nabla w_\varepsilon|^2}\\
&\ge &1-C\varepsilon^{-\sigma}.
\end{array}
$$
Thus  we get 
$\tau_2(w_\varepsilon)=1+O(\varepsilon^{-\sigma})$ as $\varepsilon\to \infty$. This proved the claim.

 If $N\ge 3$, by Lemma 3.2, (H2) and the boundedness of $\{w_\varepsilon\}$, we find
$$
\begin{array}{rcl}
m_\infty&\le& \sup_{t\ge 0} J_\varepsilon((w_\varepsilon)_t)+\varepsilon^{-\sigma} t^N_\varepsilon \left(\varepsilon^{-\tau q}\int_{\mathbb R^N}G(\varepsilon^\tau w_\varepsilon)-\frac{1}{2}\int_{\mathbb R^N}|w_\varepsilon|^{2}\right)\\
&\le &m_\varepsilon +C\varepsilon^{-\sigma},
\end{array}
\eqno(5.10)
$$
where 
$$
t_\varepsilon=\left(\frac{\int_{\mathbb R^N}|\nabla w_\varepsilon|^2}{\int_{\mathbb R^N}(I_\alpha\ast |w_\varepsilon|^{2_\alpha^*})|w_\varepsilon|^{2_\alpha^*}}\right)^{\frac{1}{2+\alpha}}=(\tau_2(w_\varepsilon))^{\frac{1}{2+\alpha}}.
$$

Next, we assume $N\ge 5$. 
For each $\rho>0$, the family  $V_\rho(x):=\rho^{-\frac{N-2}{2}}V_1(x/\rho)$ are radial ground states of  $-\Delta w=(I_\alpha\ast |w|^{2_\alpha^*})w^{2_\alpha^*-1}$, and verify that
 $$
 \|V_\rho\|_2^2=\rho^2\|V_1\|_2^2, \qquad \int_{\mathbb R^N}|V_\rho|^q=\rho^{N-\frac{N-2}{2}q}\int_{\mathbb R^N}|V_1|^q.
 \eqno(5.11)
 $$ 
 Let $g_0(\rho)=\frac{A}{q}\int_{\mathbb R^N}|V_\rho|^q-\frac{1}{2}\int_{\mathbb R^N}|V_\rho|^2$. Then there exists $\rho_0=\rho(q)\in (0,+\infty)$ with
 $$
 \rho_0=\left(\frac{A[2N-q(N-2)]\int_{\mathbb R^N}|V_1|^q}{2q\int_{\mathbb R^N}|V_1|^2}\right)^{\frac{2}{(N-2)(q-2)}}
 $$
such that 
 $$
 g_0(\rho_0)=\sup_{\rho>0}g_0(\rho)=\frac{A(N-2)(q-2)}{4q}\left(\frac{A[2N-q(N-2)]\int_{\mathbb R^N}|V_1|^q}{2q\int_{\mathbb R^N}|V_1|^2}\right)^{\frac{2N-q(N-2)}{(N-2)(q-2)}}\int_{\mathbb R^N}|V_1|^q.
 $$ 
 Let $V_0=V_{\rho_0}$, then there exists $t_\varepsilon\in (0,+\infty)$ such that
 $$
\begin{array}{rcl}
m_\varepsilon&\le &\sup_{t\ge 0}J_\varepsilon((V_0)_t)=J_\varepsilon( (V_0)_{t_\varepsilon})\\
&=&\frac{t^{N-2}_\varepsilon}{2}\int_{\mathbb R^N}|\nabla V_0|^2-\frac{t^{N+\alpha}_\varepsilon}{22_\alpha^*}(I_\alpha\ast |V_0|^{2_\alpha^*})|V_0|^{2_\alpha^*}\\
&\mbox{}& -\varepsilon^{-\sigma}t_\varepsilon^N\{\varepsilon^{-\tau q}\int_{\mathbb R^N}G(\varepsilon^\tau V_0)-\frac{1}{2}\int_{\mathbb R^N}|V_0|^2\}\\
&\le &\sup_{t\ge 0}\left(\frac{t^{N-2}}{2}-\frac{t^{N+\alpha}}{2p}\right)\int_{\mathbb R^N}|\nabla V_0|^2 \\
&\mbox{}&-\varepsilon^{-\sigma}t_\varepsilon^N\{\varepsilon^{-\tau q}\int_{\mathbb R^N}G(\varepsilon^\tau V_0)-\frac{1}{2}\int_{\mathbb R^N}|V_0|^2\}
\\
&=& m_\infty-\varepsilon^{-\sigma}t_\varepsilon^N\{\varepsilon^{-\tau q}\int_{\mathbb R^N}G(\varepsilon^\tau V_0)-\frac{1}{2}\int_{\mathbb R^N}|V_0|^2\}.
\end{array}
\eqno(5.12)
$$
Moreover, $t_\varepsilon$ satisfies 
$$
\frac{1}{2^*t^2}\int_{\mathbb R^N}|\nabla V_0|^2=\frac{t^{\alpha}}{2^*}\int_{\mathbb R^N}(I_\alpha\ast |V_0|^{2_\alpha^*})|V_0|^{2_\alpha^*}+\varepsilon^{-\sigma}\{\varepsilon^{-\tau q}\int_{\mathbb R^N}G(\varepsilon^\tau V_0)-\frac{1}{2}\int_{\mathbb R^N}|V_0|^2\}.
$$

If $t_\varepsilon\ge 1$, then
$$
\int_{\mathbb R^N}|\nabla V_0|^2\ge t_\varepsilon^2\left\{\int_{\mathbb R^N}(I_\alpha\ast |V_0|^{2_\alpha^*})|V_0|^{2_\alpha^*}+\varepsilon^{-\sigma-\tau q}\int_{\mathbb R^N}G(\varepsilon^\tau V_0)-\frac{1}{2}\varepsilon^{-\sigma}\int_{\mathbb R^N}|V_0|^2\right\}.
$$
Hence
$$
t_\varepsilon^2\le \frac{\int_{\mathbb R^N}|\nabla V_0|^2}{\int_{\mathbb R^N}(I_\alpha\ast |V_0|^{2_\alpha^*})|V_0|^{2_\alpha^*}+2^*\varepsilon^{-\sigma} \{\varepsilon^{-\tau q}\int_{\mathbb R^N}G(\varepsilon^\tau V_0)-\frac{1}{2}\int_{\mathbb R^N}|V_0|^2\}}.
$$
If $t_\varepsilon\le  1$, then
$$
\int_{\mathbb R^N}|\nabla V_0|^2\le t_\varepsilon^2\left\{\int_{\mathbb R^N}(I_\alpha\ast |V_0|^{2_\alpha^*})|V_0|^{2_\alpha^*}+\varepsilon^{-\sigma-\tau q}\int_{\mathbb R^N}G(\varepsilon^\tau V_0)-\frac{1}{2}\varepsilon^{-\sigma}\int_{\mathbb R^N}|V_0|^2\right\}.
$$
Hence
$$
t_\varepsilon^2\ge \frac{\int_{\mathbb R^N}|\nabla V_0|^2}{\int_{\mathbb R^N}(I_\alpha\ast |V_0|^{2_\alpha^*})|V_0|^{2_\alpha^*}+2^*\varepsilon^{-\sigma} \{\varepsilon^{-\tau q}\int_{\mathbb R^N}G(\varepsilon^\tau V_0)-\frac{1}{2}\int_{\mathbb R^N}|V_0|^2\}}.
$$
Since
 $$
 \int_{\mathbb R^N}(I_\alpha\ast |V_0|^{2_\alpha^*})|V_0|^{2_\alpha^*}=\int_{\mathbb R^N}|\nabla V_0|^2,
 $$
we conclude that $\lim_{\varepsilon\to \infty}t_\varepsilon=1$.  
By (H2), it follows that
$$
\lim_{\varepsilon\to \infty}\varepsilon^{-\tau q}\int_{\mathbb R^N}G(\varepsilon^\tau V_0)-\frac{1}{2}\int_{\mathbb R^N}|V_0|^2
=\frac{A}{q}\int_{\mathbb R^N}|V_0|^q-\frac{1}{2}\int_{\mathbb R^N}|V_0|^2=g_0(\rho_0),
$$
Thus, by (5.12), there exists a constant $C>0$ such that
$$
m_\varepsilon\le m_\infty-C\varepsilon^{-\sigma},
$$
for large  $\varepsilon>0$. The proof is complete. 
\end{proof}

\noindent{\bf Lemma 5.4.}  {\it There exists a constant $\varpi=\varpi(q)>0$ such that  for $\varepsilon>0$ large, 
$$
m_\varepsilon\le \left\{\begin{array}{rcl}
m_\infty-(\varepsilon\ln \varepsilon)^{-\frac{4-q}{q-2}}\varpi, \
&if&    N=4, \\
m_\infty-\varepsilon^{-\frac{6-q}{2(q-4)}}\varpi, \ \  \   \ &if&   N=3 \  and \  q\in (4,6).
\end{array} \right.
$$}
\begin{proof}  
Let  $\rho>0$ and $R$ be a large parameter, and $\eta_R\in C_0^\infty(\mathbb R)$ is a cut-off function such that $\eta_R(r)=1$ for $|r|<R$, $0<\eta_R(r)<1$ for $R<|r|<2R$, $\eta_R(r)=0$ for $|r|>2R$ and $|\eta'_R(r)|\le 2/R$.  

For $\ell\gg 1$, a straightforward computation shows that 
$$
\int_{\mathbb R^N}|\nabla (\eta_\ell V_1)|^2=\left\{ \begin{array}{rcl} \frac{2(N+\alpha)}{2+\alpha}m_\infty+O(\ell^{-2}), \quad \qquad &{\rm if} & \  \ N=4,\\
\frac{2(N+\alpha)}{2+\alpha}m_\infty+O(\ell^{-1}),  \qquad \  \    \ &{\rm if} & \   \ N=3.
\end{array}\right.
$$ 
$$
\int_{\mathbb R^N} (I_\alpha\ast |\eta_\ell V_1|^{2_\alpha^*})|\eta_\ell V_1|^{2_\alpha^*}=\frac{2(N+\alpha)}{2+\alpha}m_\infty+O(\ell^{-N}),
$$
$$
\int_{\mathbb R^N}|\eta_\ell V_1|^2=\left\{ \begin{array}{rcl} \ln \ell(1+o(1)), \quad &{\rm if} & \  \ N=4,\\
\ell(1+o(1)), \quad  \quad &{\rm if} & \   \ N=3.
\end{array}\right.
$$ 
By Lemma 3.2, we find
$$
\begin{array}{rcl}
m_\varepsilon&\le &\sup_{t\ge 0}J_\varepsilon((\eta_RV_\rho)_t)=J_\varepsilon((\eta_RV_\rho)_{t_\varepsilon})\\
&\le &\sup_{t\ge 0}\left(\frac{t^{N-2}}{2}\int_{\mathbb R^N}|\nabla(\eta_RV_\rho)|^2-\frac{t^{N+\alpha}}{2p}\int_{\mathbb R^N}(I_\alpha\ast |\eta_RV_\rho|^{2_\alpha^*})|\eta_RV_\rho|^{2_\alpha^*}\right)\\
&\quad & -\varepsilon^{-\sigma}
t_\varepsilon^N\left[\varepsilon^{-\tau q}\int_{\mathbb R^N}G(\varepsilon^\tau\eta_RV_\rho)-\frac{1}{2}\int_{\mathbb R^N}|\eta_RV_\rho|^2\right]\\
&= & (I)-\varepsilon^{-\sigma} (II),
\end{array}
\eqno(5.13)
$$
where $t_\varepsilon\in (0, +\infty)$ is the unique critical point of the function $g(t)$ defined by 
$$
\begin{array}{rcl}
g(t)&=&\frac{t^{N-2}}{2}\int_{\mathbb R^N}|\nabla(\eta_RV_\rho)|^2+ \frac{t^N}{2}\varepsilon^{-\sigma}\int_{\mathbb R^N}|\eta_RV_\rho|^2\\
&\quad &-\frac{t^{N+\alpha}}{22_\alpha^*}\int_{\mathbb R^N}(I_\alpha\ast |\eta_RV_\rho|^{2_\alpha^*})|\eta_RV_\rho|^{2_\alpha^*}-t^N\varepsilon^{-\sigma-\tau q}\int_{\mathbb R^N}G(\varepsilon^\tau \eta_RV_\rho).
\end{array}
$$
That is, $t=t_\varepsilon$ solves the equation $\ell_1(t)=\ell_2(t)$, where
$$
\ell_1(t):=\frac{1}{2^*t^2}\int_{\mathbb R^N}|\nabla(\eta_RV_\rho)|^2
$$
and 
$$
\ell_2(t):=\frac{1}{2^*}t^\alpha\int_{\mathbb R^N}(I_\alpha\ast |\eta_RV_\rho|^{2_\alpha^*})|\eta_RV_\rho|^{2_\alpha^*}
+\varepsilon^{-\sigma-\tau q}\int_{\mathbb R^N}G(\varepsilon^\tau \eta_RV_\rho)-\frac{1}{2}\varepsilon^{-\sigma}\int_{\mathbb R^N}|\eta_RV_\rho|^2.
$$

If $t_\varepsilon\ge 1$, then
$$
\begin{array}{rcl}
\frac{1}{2^*t_\varepsilon^{2}}\int_{\mathbb R^N}|\nabla(\eta_RV_\rho)|^2&\ge& \frac{1}{2^*}\int_{\mathbb R^N}(I_\alpha\ast |\eta_RV_\rho|^{2_\alpha^*})|\eta_RV_\rho|^{2_\alpha^*}\\
&\quad &
+\varepsilon^{-\sigma-\tau q}\int_{\mathbb R^N}G(\varepsilon^\tau \eta_RV_\rho)-\frac{1}{2}\varepsilon^{-\sigma}\int_{\mathbb R^N}|\eta_RV_\rho|^2,
\end{array}
$$
and hence 
$$
\begin{array}{rcl}
t_\varepsilon^2
&\le &\frac{\int_{\mathbb R^N}|\nabla(\eta_RV_\rho)|^2}{\int_{\mathbb R^N}(I_\alpha\ast |\eta_RV_\rho|^{2_\alpha^*})|\eta_RV_\rho|^{2_\alpha^*}+2^*\varepsilon^{-\sigma-\tau q}\int_{\mathbb R^N}G(\varepsilon^\tau \eta_RV_\rho)-\frac{1}{2}\varepsilon^{-\sigma}\int_{\mathbb R^N}|\eta_RV_\rho|^2}.
\end{array}
\eqno(5.14)
$$
If $t_\varepsilon\le 1$, then
$$
\begin{array}{rcl}
\frac{1}{2^*t_\varepsilon^2}\int_{\mathbb R^N}|\nabla(\eta_RV_\rho)|^2&\le& \frac{1}{2^*}\int_{\mathbb R^N}(I_\alpha\ast |\eta_RV_\rho|^{2_\alpha^*})|\eta_RV_\rho|^{2_\alpha^*}\\
&\quad &
+\varepsilon^{-\sigma-\tau q}\int_{\mathbb R^N}G(\varepsilon^\tau \eta_RV_\rho)-\frac{N}{2}\varepsilon^{-\sigma}\int_{\mathbb R^N}|\eta_RV_\rho|^2,
\end{array}
$$
and hence 
$$
\begin{array}{rcl}
t_\varepsilon^2
&\ge& \frac{\int_{\mathbb R^N}|\nabla(\eta_RV_\rho)|^2}{\int_{\mathbb R^N}(I_\alpha\ast |\eta_RV_\rho|^{2_\alpha^*})|\eta_RV_\rho|^{2_\alpha^*}+2^*\varepsilon^{-\sigma-\tau q}\int_{\mathbb R^N}G(\varepsilon^\tau \eta_RV_\rho)-\frac{1}{2}\varepsilon^{-\sigma}\int_{\mathbb R^N}|\eta_RV_\rho|^2}.
\end{array}
\eqno(5.15)
$$
Therefore, set $\ell=R/\rho$, then we obtain
$$
\begin{array}{rcl}
\limsup_{\varepsilon\to \infty}|t_\varepsilon^2-1|&\le &\left|\frac{\int_{\mathbb R^N}|\nabla(\eta_RV_\rho)|^2}{\int_{\mathbb R^N}(I_\alpha\ast |\eta_RV_\rho|^{2_\alpha^*})|\eta_RV_\rho|^{2_\alpha^*}}-1\right|=\left|\frac{\int_{\mathbb R^N}|\nabla(\eta_\ell V_1)|^2}{\int_{\mathbb R^N}(I_\alpha\ast |\eta_\ell V_1|^{2_\alpha^*})|\eta_\ell V_1|^{2_\alpha^*}}-1\right|\to 0,
\end{array}
$$
as $\ell \to \infty$. Therefore, we get 

$$
\begin{array}{rcl}
(I)&=&\frac{2+\alpha}{2(N+\alpha)}\frac{(\int_{\mathbb R^N}|\nabla(\eta_\ell V_1)|^2)^{\frac{N+\alpha}{2+\alpha}}}{(\int_{\mathbb R^N}(I_\alpha\ast |\eta_\ell V_1|^{2_\alpha^*})|\eta_\ell V_1|^{2_\alpha^*})^{\frac{N-2}{2+\alpha}}}\\
&=&\left\{\begin{array}{rcl} \frac{2+\alpha}{2(N+\alpha)}\frac{(\frac{2(N+\alpha)}{2+\alpha}m_\infty+O(\ell^{-2}))^{\frac{N+\alpha}{2+\alpha}}}{(\frac{2(N+\alpha)}{2+\alpha}m_\infty+O(\ell^{-4}))^{\frac{N-2}{2+\alpha}}}, \  \  {\rm if} \   \  N=4,\\
\frac{2+\alpha}{2(N+\alpha)}\frac{(\frac{2(N+\alpha)}{2+\alpha}m_\infty+O(\ell^{-1}))^{\frac{N+\alpha}{2+\alpha}}}{(\frac{2(N+\alpha)}{2+\alpha}m_\infty+O(\ell^{-3}))^{\frac{N-2}{2+\alpha}}}, \  \  {\rm if} \   \  N=3,
\end{array}\right.
\\
&=&\left\{\begin{array}{rcl} m_\infty+O(\ell^{-2}), \  \  {\rm if} \   \  N=4,\\
m_\infty+O(\ell^{-1}),   \   \   {\rm if}  \   \  N=3.\end{array}\right.
\end{array}
\eqno(5.16)
$$
Since 
$$
\begin{array}{rcl}
\phi(\rho)&=&\varepsilon^{-\tau q}\int_{\mathbb R^N}G(\varepsilon^\tau \eta_RV_\rho)-\frac{1}{2}\int_{\mathbb R^N}|\eta_RV_\rho|^2\\
&=&\frac{A}{q}\int_{\mathbb R^N}|\eta_RV_\rho|^q-\frac{1}{2}\int_{\mathbb R^N}|\eta_RV_\rho|^2\\
&=&\frac{A}{q}\rho^{N-\frac{N-2}{2}q}\int_{\mathbb R^N}|\eta_\ell V_1|^q-\frac{1}{2}\rho^2\int_{\mathbb R^N}|\eta_\ell V_1|^2
\end{array}
$$
take its maximum value $\varphi(\rho_\ell)$ at the unique point 
$$
\begin{array}{rcl}
\rho_\ell:&=&\left(\frac{A[2N-q(N-2)]\int_{\mathbb R^N}|\eta_\ell V_1|^q}{2q\int_{\mathbb R^N}|\eta_\ell V_1|^2}\right)^{\frac{2}{(N-2)(q-2)}}\\
&\sim &\left\{\begin{array}{rcl}
(\ln \ell)^{-\frac{2}{(N-2)(q-2)} }\quad &{\rm if}&  \   N=4,\\
\ell^{-\frac{2}{(N-2)(q-2)}} \quad &{\rm if}&  \   N=3,
\end{array}\right.
\end{array}
$$
we obtain 
$$
\begin{array}{rcl}
\phi(\rho_\ell)&=&\sup_{\rho\ge 0}\phi(\rho)\\
&=&\frac{4+q(N-2)-2N}{4q}\rho_\ell^{N-\frac{N-2}{2}q}\int_{\mathbb R^N}|\eta_\ell V_1|^q\\
&=&\frac{4+q(N-2)-2N}{4q}\left(\frac{2N-q(N-2)}{2q}\right)^{\frac{2N-q(N-2)}{(N-2)(q-2)}}\frac{(\int_{\mathbb R^N}|\eta_\ell V_1|^q)^{\frac{4}{(N-2)(q-2)}}}{(\int_{\mathbb R^N}|\eta_\ell V_1|^2)^{\frac{2N-q(N-2)}{(N-2)(q-2)}}}\\
&\le &\frac{4+q(N-2)-2N}{4q}\left(\frac{2N-q(N-2)}{2q}\right)^{\frac{2N-q(N-2)}{(N-2)(q-2)}}\int_{\mathbb R^N}|\eta_\ell V_1|^{2^*}\\
&\to &\frac{4+q(N-2)-2N}{4q}\left(\frac{2N-q(N-2)}{2q}\right)^{\frac{2N-q(N-2)}{(N-2)(q-2)}}\int_{\mathbb R^N}| V_1|^{2^*},
\end{array}
$$
as $\ell \to +\infty$, where we have used the interpolation inequality
$$
\int_{\mathbb R^N}|\eta_\ell V_1|^q\le \left(\int_{\mathbb R^N}|\eta_\ell V_1|^2\right)^{\frac{2^*-q}{2^*-2}}\left(\int_{\mathbb R^N}|\eta_\ell V_1|^{2^*}\right)^{\frac{q-2}{2^*-2}}.
$$
Since 
$
\int_{\mathbb R^N}|\eta_\ell V_1|^q\to \int_{\mathbb R^N}|V_1|^q,
$
as $\ell \to +\infty$, it follows that
$$
\begin{array}{rcl}
\phi(\rho_\ell)&=&\frac{(N-2)(q-2)}{4q}A^{\frac{4}{(N-2)(q-2)}}\left(\frac{[2N-q(N-2)]\int_{\mathbb R^N}|\eta_\ell V_1|^q}{2q\int_{\mathbb R^N}|\eta_\ell V_1|^2}\right)^{\frac{2N-q(N-2)}{(N-2)(q-2)}}\int_{\mathbb R^N}|\eta_\ell V_1|^q\\
&=& \left\{\begin{array}{rcl} C(\ln \ell(1+o(1))^{-\frac{2N-q(N-2)}{(N-2)(q-2)}}  \qquad &{\rm if}&  \   \  N=4,\\
C(\ell(1+o(1))^{-\frac{2N-q(N-2)}{(N-2)(q-2)}} \qquad &{\rm if }& \   \ N=3.
\end{array}\right.
\end{array}
$$
Thus, for large $\varepsilon>0$, we have 
$$
\begin{array}{rcl}
(II)&=&\phi(\rho_\ell)+(t_\varepsilon^N-1)\phi(\rho_\ell )\\
&\sim &\left\{\begin{array}{rcl} (\ln \ell)^{-\frac{2N-q(N-2)}{(N-2)(q-2)}},  \qquad &{\rm if}&  \   \  N=4,\\
\ell^{-\frac{2N-q(N-2)}{(N-2)(q-2)}}, \qquad &{\rm if }& \   \ N=3.
\end{array}\right.\end{array}
$$
It follows that if $N=4$, then 
$$
\begin{array}{rcl}
m_\varepsilon &\le &(I)-\varepsilon^{-\sigma} (II)\\
&\le &m_\infty+O(\ell^{-2})-C\varepsilon^{-\sigma} (\ln \ell)^{-\frac{2N-q(N-2)}{(N-2)(q-2)} }.
\end{array}
\eqno(5.17)
$$
Take $\ell=\varepsilon^M$. Then 
$$
m_\varepsilon\le m_\infty+O(\varepsilon^{-2M})-C\varepsilon^{-\sigma} M^{-\frac{2N-q(N-2)}{(N-2)(q-2)}}(\ln\varepsilon)^{-\frac{2N-q(N-2)}{(N-2)(q-2) }}.
$$
If $M>\frac{4-q}{2(q-2)}$, then $2M>\sigma$, and hence
$$
m_\varepsilon\le m_\infty-\varepsilon^{-\sigma}(\ln\varepsilon)^{-\frac{2N-q(N-2)}{(N-2)(q-2)} }\varpi=m_\infty-(\varepsilon\ln \varepsilon)^{-\frac{4-q}{q-2}}\varpi,
\eqno(5.18)
$$
for large  $\varepsilon>0$, where 
$$
\varpi=\frac{1}{2}CM^{-\frac{2N-q(N-2)}{(N-2)(q-2)}}.
$$

If  $N=3$, then 
$$
\begin{array}{rcl}
m_\lambda&\le &(I)-\lambda^\sigma (II)\\
&\le &m_\infty+O(\ell^{-1})-C\varepsilon^{-\sigma} \ell^{-\frac{2N-q(N-2)}{(N-2)(q-2)}}.
\end{array}
\eqno(5.19)
$$
Take $\ell=\delta^{-1}\varepsilon^\tau$. Then 
$$
m_\varepsilon\le m_\infty+\delta O(\varepsilon^{-\tau})-C\varepsilon^{-\sigma} \delta^{\frac{2N-q(N-2)}{(N-2)(q-2)}}\varepsilon^{-\tau\frac{2N-q(N-2)}{(N-2)(q-2)} }
$$
If $q\in (4,6)$ and  $\tau=\frac{6-q}{2(q-4)}$, then
$$
\tau=\frac{2(N-2)}{2+q(N-2)-2N}=\frac{2}{q-4},
$$
then
$$
m_\varepsilon\le m_\infty+(\delta O(1)-C\delta^{\frac{2N-q(N-2)}{(N-2)(q-2)}})\varepsilon^{-\frac{6-q}{2(q-4)}}.
$$
Since 
$$
1>\frac{2N-q(N-2)}{(N-2)(q-2)},
$$
it follows that for some small $\delta>0$, there exists $\varpi>0$ such that 
$$
m_\varepsilon\le m_\infty-\varepsilon^{-\frac{6-q}{2(q-4)}}\varpi.
$$
This completes the proof. \end{proof}

Combining Lemma 5.3 and Lemma 5.4, we get the following

\smallskip

\noindent{\bf Lemma 5.5.} {\it Let $\delta_\varepsilon:=m_\infty-m_\varepsilon$, then 
$$
\varepsilon^{-\sigma}\gtrsim \delta_\varepsilon\gtrsim \left\{\begin{array}{rcl} 
 \varepsilon^{-\sigma}, \  \  \     \quad  \qquad &{ if}&  \    N\ge 5,\\
(\varepsilon\ln\varepsilon)^{-\frac{4-q}{q-2}},  \quad  \ &{ if}& \   N=4,\\
\varepsilon^{-\frac{6-q}{2(q-4)}},   \qquad  \   &{ if}& \     N=3 \ and \ q\in (4, 6).
 \end{array}\right.
 $$   }


\noindent{\bf Lemma 5.6.} {\it 
Assume $N\ge 5$. Then $\|w_\varepsilon\|_q\sim 1$ as $\varepsilon\to \infty$.}

\begin{proof} 
By (5.7),  we have
$$
\begin{array}{rcl}
m_\infty&\le& m_\varepsilon+\varepsilon^{-\sigma} (\tau_2(w_\varepsilon))^{\frac{N}{2+\alpha}}\varepsilon^{-\tau q}\int_{\mathbb R^N}G(\varepsilon^\tau w_\varepsilon)\\
&\le & m_\varepsilon+\frac{A}{q}\varepsilon^{-\sigma}(\tau_2(w_\varepsilon)^{\frac{N}{2+\alpha}}(\int_{\mathbb R^N}|w_\varepsilon|^q+o_\varepsilon(1)).
\end{array}
$$
Therefore, it follows from (5.6) and Corollary 5.5 that 
$$
\|w_\varepsilon\|_q^q\ge \frac{m_\infty-m_\varepsilon}{(\tau_2(w_\varepsilon))^{\frac{N}{2+\alpha}}}\cdot \frac{q}{A}\varepsilon^{\sigma}+o_\varepsilon(1)\ge \frac{Cq}{A(\tau_2(w_\varepsilon))^{\frac{N}{2+\alpha}}}\ge C>0,
$$
which together with the boundedness of $\{w_\varepsilon\}$ implies the desired  conclusion. 
\end{proof}

\noindent{\bf Lemma 5.7.} {\it Let $N\ge 5$ and $q\in (2,2^*)$,  then there exists a $\zeta_\varepsilon\in (0,\infty)$ verifying 
$$
\zeta_\varepsilon\sim  \varepsilon^{-\frac{4}{(N-2)(q-2)}}
$$
such that the rescaled ground states 
$
w_\varepsilon(x)=\zeta_\varepsilon^{\frac{N-2}{2}}u_\varepsilon(\zeta_\varepsilon x)
$
converges to  $V_{\rho_0}$ in $H^1(\mathbb R^N)$ as $\varepsilon\to \infty$,  where $V_{\rho_0}$ is  a positive ground sate of the equation $-\Delta V=(I_\alpha\ast |V|^{2_\alpha^*})V^{2_\alpha^*-1}$ with
$$
 \rho_0=\left(\frac{NA(2^*-q)\int_{\mathbb R^N}|V_1|^q}{2^*q\int_{\mathbb R^N}|V_1|^2}\right)^ {\frac{2}{(N-2)(q-2)}}.
 \eqno(5.20)
 $$
In the lower dimension cases $N=4$ and $N=3$, there exists $\xi_\varepsilon\in (0,+\infty)$ with $\xi_\varepsilon\to 0$ such that 
$$
w_\varepsilon-\xi_\varepsilon^{-\frac{N-2}{2}}V_1(\xi^{-1}_\varepsilon\cdot)\to 0
$$
as $\varepsilon\to \infty$ in $D^{1,2}(\mathbb R^N)$ and $L^{2^*}(\mathbb R^N)$. }

\begin{proof} We adopt a proof of \cite[Lemma 5.7]{Ma-2}. Note that $w_\varepsilon$ is a positive radially symmetric function, and by Lemma 5.2, $\{w_\varepsilon\}$ is bounded in $H^1(\mathbb R^N)$. Then there exists $w_\infty\in H^1(\mathbb R^N)$ verifying $-\Delta w=(I_\alpha \ast w^{2_\alpha^*})w^{2_\alpha^*-1}$ such that 
$$
w_\varepsilon \rightharpoonup w_\infty  \quad {\rm weakly \ in} \  H^1(\mathbb R^N), \quad  w_\varepsilon\to w_\infty \quad {\rm in} \ L^p(\mathbb R^N) \quad {\rm for \ any} \ p\in (2,2^*),
\eqno(5.21)
$$
and 
$$
w_\varepsilon(x)\to w_\infty(x) \quad a. \ e. \  {\rm on} \ R^N,  \qquad  w_\varepsilon \to w_\infty \quad {\rm in} \   L^2_{loc}(\mathbb R^N).
\eqno(5.22)
$$
Observe that
$$
J_\infty(w_\varepsilon)=J_\varepsilon(w_\varepsilon)+\varepsilon^{-\sigma-\tau q}\int_{\mathbb R^N}G(\varepsilon^\tau w_\varepsilon)-\frac{\varepsilon^{-\sigma}}{2}\int_{\mathbb R^N}|w_\varepsilon|^2=m_\varepsilon+o_\varepsilon(1)=m_\infty+o_\varepsilon(1),
$$
and 
$$
J'_\infty(w_\varepsilon)w=J'_\varepsilon(w_\varepsilon)w+\varepsilon^{-\sigma-\tau(q-1)}\int_{\mathbb R^N}g(\varepsilon^\tau w_\varepsilon)w-\varepsilon^{-\sigma}\int_{\mathbb R^N}w_\varepsilon w=o_\varepsilon(1).
$$
Therefore, $\{w_\varepsilon\}$ is a PS sequence of $J_\varepsilon$ at level  $m_\infty=\frac{2+\alpha}{2(N+\alpha)}S_\alpha^{\frac{N+\alpha}{2+\alpha}}$.

By Lemma 3.4,  it is standard to show \cite{Tintarev-1} that  there exists $\zeta^{(j)}_\varepsilon\in (0,+\infty)$, $w^{(j)}\in D^{1,2}(\mathbb R^N)$ with $j=1,2,\cdots, k$, $k$  a non-negative integer, such that
$$
w_\varepsilon=w_\infty+\sum_{j=1}^k(\zeta^{(j)}_\varepsilon)^{-\frac{N-2}{2}}w^{(j)}((\zeta^{(j)}_\lambda)^{-1} x)+\tilde w_\varepsilon,
\eqno(5.23)
$$
where $\tilde w_\varepsilon \to 0$ in $L^{2^*}(\mathbb R^N)$ and $w^{(j)}$ are nontrivial solutions of the limit equation $-\Delta w=(I_\alpha\ast w^{2_\alpha^*})w^{2_\alpha^*-1}$. Moreover, we have
$$
\limsup_{\varepsilon\to \infty}\|w_\varepsilon\|^2_{D^1(\mathbb R^N)}\ge \|w_\infty\|^2_{D^1(\mathbb R^N)}+\sum_{j=1}^k\|w^{(j)}\|^2_{D^1(\mathbb R^N)}
\eqno(5.24)
$$
and 
$$
m_\infty=J_\infty(w_\infty)+\sum_{j=1}^kJ_0(w^{(j)}).
\eqno(5.25)
$$
Moreover, $J_\infty(w_\infty)\ge 0$ and $J_\infty(w^{(j)})\ge m_\infty$ for all $j=1,2,\cdots, k.$

If $N\ge 5$, then by Lemma 5.6, we have $w_\infty\not=0$ and hence $J_\infty(w_\infty)=m_\infty$ and  $k=0$.  Thus $w_\varepsilon\to w_\infty$ in $L^{2^*}(\mathbb R^N)$. Since $J_\infty'(w_\varepsilon)\to 0$, it follows that
$w_\varepsilon\to w_\infty$ in $D^{1,2}(\mathbb R^N)$.  

Since $w_\varepsilon(x)$ is radial and radially decreasing, for every $x\in\mathbb R^N\setminus \{0\}$, we have 
$$
w^2_\varepsilon(x)\le \frac{1}{|B_{|x|}|}\int_{B_{|x|}}|w_\varepsilon|^2\le \frac{1}{|x|^N}\int_{\mathbb R^N}|w_\varepsilon|^2\le \frac{C}{|x|^N},
$$
then 
$$
w_\varepsilon(x) \le C|x|^{-\frac{N}{2}}, \qquad |x|\ge 1.
\eqno(5.26)
$$

If  $\alpha>N-4$, then we have $2_\alpha^*=\frac{N+\alpha}{N-2}>2$ and hence
$$
|w_\varepsilon|^{2_\alpha^*}|x|^N\le C|x|^{-\frac{N}{2}{2_\alpha^*}+N}=C|x|^{-\frac{N}{2}(2_\alpha^*-2)}\to 0, \qquad {\rm as} \  |x|\to +\infty.
$$
By virtue of  Lemma 3.10, we obtain 
$$
(I_\alpha\ast |w_\varepsilon|^{2_\alpha^*})(x)\le C|x|^{-N+\alpha}, \qquad |x|\ge 1.
$$
and then
$$
(I_\alpha\ast |w_\varepsilon|^{2_\alpha^*})(x)|w_\varepsilon|^{2_\alpha^*-2}(x)\le C|x|^{-\frac{N^2-N\alpha+4\alpha}{2(N-2)}}, \qquad |x|\ge \tilde R.
\eqno(5.27)
$$
As in the proof of \cite[Lemma 5.7]{Ma-2}, for any small $\nu_0>0$, by using a comparison argument, we can find $R>0$ such that for all large $\varepsilon>0$ there holds
$$
w_\varepsilon(x)\lesssim |x|^{-(N-2-\nu_0)},  \qquad {\rm for} \  \ |x|\ge R.
\eqno(5.28)
$$
When $\nu_0>0$ is small enough, the domination is in $L^2(\mathbb R^N)$ for $N\ge 5$, and this shows, by the dominated convergence theorem, that 
$w_\varepsilon\to w_\infty$ in $L^2(\mathbb R^N)$.  Thus, we conclude that $w_\varepsilon\to w_0$ in  $H^1(\mathbb R^N)$. Moreover, by (5.4), we obtain
$$
\|w_\infty\|_2^2=\frac{NA(2^*-q)}{2^*q}\int_{\mathbb R^N}|w_\infty|^q,
$$
from which it follows that $w_\infty=V_{\rho_0}$ with  
$$
 \rho_0=\left(\frac{NA(2^*-q)\int_{\mathbb R^N}|V_1|^q}{2^*q\int_{\mathbb R^N}|V_1|^2}\right)^ {\frac{2}{(N-2)(q-2)}}.
 $$
Therefore, if $\alpha>N-4$, then the first  statement is valid with $\zeta_\varepsilon=\varepsilon^{-\frac{2}{(N-2)(q-2)}}$. 
If $\alpha\le N-4$, then for any $\varepsilon_n\to \infty$, up to a subsequence, we can assume that
$w_{\varepsilon_n}\to V_\rho$ in $D^{1,2}(\mathbb R^N)$ with $\rho\in (0,\rho_0]$. Moreover, $w_{\varepsilon_n}\to V_\rho$ in $L^2(\mathbb R^N)$ if and only if
$\rho=\rho_0$. Arguing as in the proof of \cite[Theorem 2.1]{Ma-2} (see also the proof of Theorem 2.1 in the present paper),  we can show that there exists a $\zeta_\varepsilon\sim \varepsilon^{-\frac{2}{(N-2)(q-2)}}$ such that
$w_\varepsilon(x)=\zeta_\varepsilon^{\frac{N-2}{2}}u_\varepsilon(\zeta_\varepsilon x )$ converges to $V_{\rho_0}$ in $L^2(\mathbb R^N)$, and hence in $H^1(\mathbb R^N)$.

If $N=4$ or $3$. By Fatou's lemma, we have $\|w_\infty\|^2_2\le \liminf_{\varepsilon\to \infty}\|v_\varepsilon\|_2^2<\infty$, therefore,  $w_\infty=0$ and hence $k=1$. Thus, we obtain  
$J_\infty(w^{(1)})=m_\infty$ and hence $w^{(1)}=V_\rho$ for some $\rho\in (0,+\infty)$. Therefore, we conclude that 
$$
w_\varepsilon-\xi_\varepsilon^{-\frac{N-2}{2}}V_1(\xi_\varepsilon^{-1}\cdot )\to 0
$$ 
in $L^{2^*}(\mathbb R^N)$ as $\lambda\to 0$, where 
$\xi_\varepsilon:=\rho\zeta_\varepsilon^{(1)}\in (0,+\infty)$ satisfying $\xi_\varepsilon\to 0$ as $\varepsilon\to \infty$.
Since 
$$
J_\infty'(w_\varepsilon-\xi_\varepsilon^{-\frac{N-2}{2}}V_1(\xi_\varepsilon^{-1}\cdot ))=J'_\infty(v_\varepsilon)+J'_\infty(V_1)+o_\varepsilon(1)=o_\varepsilon(1)
$$ 
as $\varepsilon\to \infty$, it follows that $w_\varepsilon-\xi_\varepsilon^{-\frac{N-2}{2}}V_1(\xi_\varepsilon^{-1}\cdot )\to 0$ in $D^{1,2}(\mathbb R^N)$. 
\end{proof}

\vskip 3mm

\subsection{The second scaling}

In the lower dimension cases $N=4$ and $N=3$, we further perform a  scaling
$$
\tilde w(x)=\xi_\varepsilon^{\frac{N-2}{2}} w(\xi_\varepsilon x), 
\eqno(5.29)
$$
where $\xi_\varepsilon\in (0,+\infty)$ is given in Lemma 5.7.  Then the rescaled equation is as follows
$$
-\Delta \tilde w+\varepsilon^{-\sigma}\xi_\varepsilon^{2}\tilde w=(I_\alpha\ast |\tilde w|^{2_\alpha^*})\tilde w^{2_\alpha^*-1}+\varepsilon^{-\sigma}\xi_\varepsilon^{N-\frac{N-2}{2}q}(\varepsilon^{-\tau}\xi_\varepsilon^{\frac{N-2}{2}})^{q-1}g(\varepsilon^\tau\xi_\varepsilon^{-\frac{N-2}{2}}\tilde w).
\eqno(5.30)
$$
The  corresponding energy functional is given by
$$
\begin{array}{rcl}
\tilde J_\varepsilon(\tilde w):&=&\frac{1}{2}\int_{\mathbb R^N}|\nabla\tilde  w|^2+\varepsilon^{-\sigma}\xi_\varepsilon^2|\tilde w|^2-\frac{1}{22_\alpha^*}\int_{\mathbb R^N}(I_\alpha\ast |\tilde w|^{2_\alpha^*})|\tilde w|^{2_\alpha^*}\\
&\mbox{}& -\varepsilon^{-\sigma}\xi_\varepsilon^{N-\frac{N-2}{2}q}(\varepsilon^{-\tau}\xi_\varepsilon^{\frac{N-2}{2}})^{q}\int_{\mathbb R^N}G(\varepsilon^\tau \xi_\varepsilon^{-\frac{N-2}{2}}\tilde w).
\end{array}
\eqno(5.31)
$$
Clearly, we have $\tilde J_\varepsilon(\tilde w)=J_\varepsilon(w)=I_\varepsilon(u)$. 

\smallskip

Furthermore, we  have the following lemma.

\smallskip

\noindent{\bf Lemma 5.8.}  {\it  Let $u,w,\tilde w\in H^1(\mathbb R^N)$ satisfy  (5.1) and (5.29), then the following statements hold true

(1) $ \ \|\nabla \tilde w\|_2^2= \|\nabla w\|_{2}^{2}=\|\nabla u\|_{2}^{2}, \  \int_{\mathbb R^N}(I_\alpha\ast |\tilde w|^{2_\alpha^*})|\tilde w|^{2_\alpha^*}=\int_{\mathbb R^N}(I_\alpha\ast |w|^{2_\alpha^*})|w|^{2_\alpha^*}=\int_{\mathbb R^N}(I_\alpha\ast |u|^{2_\alpha^*})|u|^{2_\alpha^*},$

(2)  $\xi_\varepsilon^{2}\|\tilde w\|^2_2=\|w\|_2^2=\varepsilon^{1+\sigma}\| u\|_2^2, \   \   \xi_\varepsilon^{N-\frac{N-2}{2}q}\|\tilde w\|^q_q=\|w\|_q^q=\varepsilon^{\sigma} \|u\|_q^q$.
}

\smallskip

Set  $\tilde w_\varepsilon(x)=\xi_\varepsilon^{\frac{N-2}{2}} w_\varepsilon(\xi_\varepsilon x)$, then by Lemma 5.7, we have 
$$
\|\nabla(\tilde w_\varepsilon-V_1)\|_2\to 0, \qquad \|\tilde w_\varepsilon-V_1\|_{2^*}\to 0,   \qquad {\rm as} \  \ \varepsilon\to \infty.
\eqno(5.32)
$$

Note that the corresponding Nehari and Poho\v zaev's identities are as follows
$$
\begin{array}{cl}
&\int_{\mathbb R^N}|\nabla \tilde w_\varepsilon|^2+\varepsilon^{-\sigma}\xi_\varepsilon^{2}\int_{\mathbb R^N}|
\tilde w_\varepsilon|^2\\
&=\int_{\mathbb R^N}(I_\alpha\ast |\tilde w_\varepsilon|^{2_\alpha^*})|\tilde w_\varepsilon|^{2_\alpha^*}+\varepsilon^{-\sigma}\xi_\varepsilon^{N-\frac{N-2}{2}q}(\varepsilon^{-\tau}\xi_\varepsilon^{\frac{N-2}{2}})^{q-1}\int_{\mathbb R^N}g(\varepsilon^\tau\xi_\varepsilon^{-\frac{N-2}{2}}\tilde w)\tilde w.
\end{array}
\eqno(5.33)
$$
and 
$$
\begin{array}{cl}
&\frac{1}{2^*}\int_{\mathbb R^N}|\nabla \tilde w_\varepsilon|^2+\frac{1}{2}
\varepsilon^{-\sigma}\xi_\varepsilon^{2}\int_{\mathbb R^N}|\tilde w_\varepsilon|^2\\
&=\frac{1}{2^*}\int_{\mathbb R^N}(I_\alpha\ast |\tilde w_\varepsilon|^{2_\alpha^*})|\tilde w_\varepsilon|^{2_\alpha^*}
+
\varepsilon^{-\sigma}\xi_\varepsilon^{N-\frac{N-2}{2}q}(\varepsilon^{-\tau}\xi_\varepsilon^{\frac{N-2}{2}})^q\int_{\mathbb R^N}G(\varepsilon^\tau\xi_\varepsilon^{-\frac{N-2}{2}} \tilde w_\varepsilon),
\end{array}
\eqno(5.34)
$$
it follows that 
$$
\begin{array}{rcl}
\frac{1}{N}\varepsilon^{-\sigma}\xi_\varepsilon^{2}\int_{\mathbb R^N}|\tilde w_\varepsilon|^2&=&\varepsilon^{-\sigma-\tau q}\xi_\varepsilon^{N}\int_{\mathbb R^N}[G(\varepsilon^{\tau}\xi_\varepsilon^{-\frac{N-2}{2}}\tilde w_\varepsilon)-\frac{1}{2^*}g(\varepsilon^\tau\xi_\varepsilon^{-\frac{N-2}{2}}\tilde w_\varepsilon)\varepsilon^\tau\xi_\varepsilon^{-\frac{N-2}{2}}\tilde w_\varepsilon]\\
&=&\varepsilon^{-\sigma}\xi_\varepsilon^{N-\frac{N-2}{2}q}[\frac{A(2^*-q)}{2^*q}\int_{\mathbb R^N}|\tilde w_\varepsilon|^q+o_\varepsilon(1)].
\end{array}
$$
Thus, we obtain
$$
\xi_\varepsilon^{\frac{(N-2)(q-2)}{2}}\int_{\mathbb R^N}|\tilde w_\varepsilon|^2=\frac{NA(2^*-q)}{2^*q}\int_{\mathbb R^N}|\tilde w_\varepsilon|^q+o_\varepsilon(1).
\eqno(5.35)
$$

To control the norm $\|\tilde w_\varepsilon\|_2$, we note that  for any $\varepsilon>0$, $\tilde w_\varepsilon>0$ satisfies the linear inequality
$$
-\Delta \tilde w_\varepsilon+\varepsilon^{-\sigma}\xi_\varepsilon^{2}\tilde w_\varepsilon>0.
\eqno(5.36)
$$

\smallskip

\noindent{\bf Lemma 5.9}  {\it There exists a constant $c>0$ such that 
$$
\tilde w_\varepsilon(x)\ge c|x|^{-(N-2)}\exp({-\varepsilon^{-\frac{\sigma}{2}}\xi_\varepsilon}|x|),  \quad |x|\ge 1.
\eqno(5.37)
$$}

The proof of the above lemma  is similar to that of   \cite[Lemma 4.8]{Moroz-1}.  As consequences, we have the following lemma.

\smallskip

\noindent{\bf Lemma 5.10.}  {\it If $N=3$, then $\|\tilde w_\varepsilon\|_2^2\gtrsim \varepsilon^{\frac{\sigma}{2}}\xi_\varepsilon^{-1}$.}

\smallskip

\noindent{\bf Lemma 5.11.}  {\it If $N=4$, then $\|\tilde w_\varepsilon\|_2^2\gtrsim  - \ln(\varepsilon^{-\sigma}\xi_\varepsilon^2).$}

\vskip 3mm

We remark that $\tilde w_\varepsilon$ is only defined for $N=4$ and $N=3$. But the following discussion also apply to the case $N\ge 5$.
To prove our main result, the key point is to show the boundedness of $\|\tilde w_\varepsilon\|_q$. We begin with a technique lemma.

\smallskip

\noindent{\bf Lemma 5.12.}  {\it  Assume $N\ge 3$, $\alpha>N-4$ and $ 2<q<2^*$. Then there exist constants $L_0>0$ and $C_0>0$ such that for any large $\varepsilon>0$ and $|x|\ge L_0\varepsilon^{\sigma/2}\xi_\varepsilon^{-1}$, 
$$
\tilde w_\varepsilon(x)\le C_0\varepsilon^{-\sigma(N-2)/4}\xi_\varepsilon^{(N-2)/2}\exp(-\frac{1}{2}\varepsilon^{-\sigma/2}\xi_\varepsilon |x|).
$$}
\begin{proof} By (5.26) and (5.27),  if $|x|\ge L_0\varepsilon^{\sigma/2}\xi_\varepsilon^{-1}$ with  $L_0>0$ being large enough, we have
$$
\begin{array}{rcl}
(I_\alpha\ast |\tilde w_\varepsilon|^{2_\alpha^*})(x)|\tilde w_\varepsilon|^{2_\alpha^*-2}(x)&=&\xi_\varepsilon^{(N-2)(2_\alpha^*-1)-\alpha}(I_\alpha\ast |w_\varepsilon|^{2_\alpha^*})(\xi_\varepsilon x)|w_\varepsilon|^{2_\alpha^*-2}(\xi_\varepsilon x)\\
&\le &C\xi_\varepsilon^2L_0^{-\frac{N^2-N\alpha+4\alpha}{2(N-2)}}\varepsilon^{-\sigma\cdot \frac{N^2-N\alpha+4\alpha}{4(N-2)}}\\
&\le &\frac{1}{4}\varepsilon^{-\sigma}\xi_\varepsilon^2,
\end{array}
$$
here we have used the fact that 
$$
\frac{N^2-N\alpha+4\alpha}{4(N-2)}>1,
$$
which follows from the inequality $N<N+2<\frac{4(\alpha+2)}{\alpha-N+4}, \ \forall \alpha\in (N-4,N)$. 

By (H1) and (H2), for any $\delta>0,$ there exists $C_\delta>0$ such that
$$
g(s)\le \delta s+C_\delta s^{q-1}, \quad for \ s\ge 0.
$$
Therefore, we have 
$$
(\varepsilon^{-\tau}\xi_\varepsilon^{\frac{N-2}{2}})^{q-1}g(\varepsilon^\tau\xi_\varepsilon^{-\frac{N-2}{2}}\tilde w_\varepsilon(x))
\le \delta (\varepsilon^{-\tau}\xi_\varepsilon^{\frac{N-2}{2}})^{q-2}\tilde w_\varepsilon+C_\delta \tilde w_\varepsilon^{q-1}.
$$
By (5.24) and (5.26), for $|x|\ge L_0\varepsilon^{\sigma/2}\xi_\varepsilon^{-1}$, we get
$$
\begin{array}{rcl}
\varepsilon^{-\sigma}\xi_\varepsilon^{N-\frac{N-2}{2}q}(\varepsilon^{-\tau}\xi_\varepsilon^{\frac{N-2}{2}})^{q-1}g(\varepsilon^\tau\xi_\varepsilon^{-\frac{N-2}{2}}\tilde w_\varepsilon(x))
&\le &\delta \varepsilon^{-\sigma-\tau (q-2)}\xi_\varepsilon^2\tilde w_\varepsilon+C_\delta \varepsilon^{-\sigma}\xi_\varepsilon^{N-\frac{N-2}{2}q}\tilde w_\varepsilon^{q-1}\\
&\le & \frac{1}{4}\varepsilon^{-\sigma}\xi_\varepsilon^2\tilde w_\varepsilon(x).
\end{array}
$$
Therefore, we obtain 
$$
-\Delta\tilde w_\varepsilon(x)+\frac{1}{2}\varepsilon^{-\sigma}\xi_\varepsilon^2\tilde w_\varepsilon(x)\le 0, \quad 
{\rm for \ all} \  |x|\ge L_0\varepsilon^{\sigma/2}\xi_\varepsilon^{-1}.
$$ 
Arguing as in  \cite[Lemma 3.2]{Akahori-2}, we prove that 
$$
\tilde w_\varepsilon(x)\le CL_0^{-\frac{N-2}{2}}e^{L_0/2}\varepsilon^{-\frac{\sigma(N-2)}{4}}\xi_\varepsilon^{\frac{N-2}{2}}e^{-\frac{1}{2}\varepsilon^{-\sigma/2}\xi_\varepsilon |x|},
$$
for all $|x|\ge L_0\varepsilon^{\sigma/2}\xi_\varepsilon^{-1}$. The proof is complete.
\end{proof}

To obtain the optimal decay estimate of $\tilde w_\varepsilon$ at infinity,
we consider the Kelvin transform of $\tilde w_\varepsilon$. For any $w\in H^1(\mathbb R^N)$, we denote by $K[w]$ the  Kelvin transform of $w$, that is,
$$
K[w](x):=|x|^{-(N-2)}w\left(\frac{x}{|x|^2}\right). 
$$
It is easy to see that $\|K[\tilde w_\lambda]\|_{L^\infty(B_1)}\lesssim 1$ implies that 
$$
\tilde w_\varepsilon(x)\lesssim |x|^{-(N-2)},  \qquad |x|\ge 1,
$$
uniformly for large $\varepsilon>0$. Thus, it suffices to show that there exists $\varepsilon_0>$ such that 
$$
\sup_{\varepsilon\in (\varepsilon_0,\infty)}\|K[\tilde w_\varepsilon]\|_{L^{\infty}(B_1)}<\infty.
\eqno(5.38)
$$

It is easy to verify that $K[\tilde w_\varepsilon]$ satisfies 
$$
\begin{array}{rcl}
-\Delta K[\tilde w_\varepsilon]+\frac{\varepsilon^{-\sigma}\xi_\varepsilon^{2}}{|x|^4}K[\tilde w_\varepsilon]&=&\frac{1}{|x|^{4}}(I_\alpha\ast |\tilde w_\varepsilon|^{2_\alpha^*})(\frac{x}{|x|^2})\tilde w_\varepsilon^{2_\alpha^*-2}(\frac{x}{|x|^2})K[\tilde w_\varepsilon]\\
&\mbox{}&+\frac{1}{|x|^4}\varepsilon^{-\sigma-\tau(q-1)}\xi_\varepsilon^{\frac{N+2}{2}}\frac{g(\varepsilon^\tau \xi_\varepsilon^{-\frac{N-2}{2}}\tilde w_\varepsilon(\frac{x}{|x|^2}))}{\tilde w_\varepsilon(\frac{x}{|x|^2})}K[\tilde w_\varepsilon].
\end{array}
\eqno(5.39)
$$
 
We also see from Lemma 5.12 that if $|x|\le \varepsilon^{-\sigma/2}\xi_\varepsilon/L_0$, then 
$$
K[\tilde w_\varepsilon](x)\lesssim \frac{1}{|x|^{N-2}}\varepsilon^{\frac{-\sigma(N-2)}{4}}\xi_\varepsilon^{\frac{N-2}{2}}e^{-\frac{1}{2}\varepsilon^{-\sigma/2}\xi_\varepsilon |x|^{-1}}.
\eqno(5.40)
$$

   Let 
$
\tilde a(x)=\frac{\varepsilon^{-\sigma}\xi_\varepsilon^2}{|x|^4}
$
and 
 $$
 \tilde b(x)=\frac{1}{|x|^{4}}(I_\alpha\ast |\tilde w_\varepsilon|^{2_\alpha^*})(\frac{x}{|x|^2})\tilde w_\varepsilon^{2_\alpha^*-2}(\frac{x}{|x|^2})
+\frac{1}{|x|^4}\varepsilon^{-\sigma-\tau(q-1)}\xi_\varepsilon^{\frac{N+2}{2}}\frac{g(\varepsilon^\tau \xi_\varepsilon^{-\frac{N-2}{2}}\tilde w_\varepsilon(\frac{x}{|x|^2}))}{\tilde w_\varepsilon(\frac{x}{|x|^2})}.
$$
Then (5.39) reads as 
$$
-\Delta K[\tilde w_\varepsilon]+\tilde a(x)K[\tilde w_\varepsilon]=\tilde b(x)K[\tilde w_\varepsilon].
\eqno(5.41)
$$
Since
$g(s)=o(s)$ as $s\to 0$ and $\lim_{s\to +\infty}g(s)s^{1-q}=A>0$, it follows that
$$
g(s)\le C(|s|^{q-1}+|s|), \quad \forall s\ge 0.
$$
 then we have 
$$
\begin{array}{rcl}
\tilde  b(x)&\le &\frac{1}{|x|^{4}}(I_\alpha\ast |\tilde w_\varepsilon|^{2_\alpha^*})(\frac{x}{|x|^2})\tilde w_\varepsilon^{2_\alpha^*-2}(\frac{x}{|x|^2})\\
&
\mbox{}&+\frac{1}{|x|^4}\varepsilon^{-\sigma}\xi_\varepsilon^{N-\frac{N-2}{2}q}|\tilde w_\varepsilon(\frac{x}{|x|^2})|^{q-2}
+\frac{1}{|x|^4}\varepsilon^{-\sigma-\tau(q-2)}\xi_\varepsilon^{2}\\
&=&\frac{1}{|x|^{4}}(I_\alpha\ast |\tilde w_\varepsilon|^{2_\alpha^*})(\frac{x}{|x|^2})\tilde w_\varepsilon^{2_\alpha^*-2}(\frac{x}{|x|^2})\\
&\mbox{}&+\frac{1}{|x|^{\gamma_q}}\varepsilon^{-\sigma}\xi_\varepsilon^{\gamma_q/2}K[\tilde w_\varepsilon]^{q-1}+\frac{1}{|x|^4}\varepsilon^{-\sigma-\tau(q-2)}\xi_\varepsilon^{2}, \end{array}
$$
here and in what follows, for any $q\in (2,2^*)$, we set
$
\gamma_q=2N-(N-2)q.
$
Choose $\eta(x)\in [0,1]$ such that 
$$
\begin{array}{rcl}
\tilde  b(x)
&=&\eta(x)\left\{\frac{1}{|x|^{4}}(I_\alpha\ast |\tilde w_\varepsilon|^{2_\alpha^*})(\frac{x}{|x|^2})\tilde w_\varepsilon^{2_\alpha^*-2}(\frac{x}{|x|^2})\right.\\
&\mbox{}& \qquad \left. +\frac{1}{|x|^{\gamma_q}}\varepsilon^{-\sigma}\xi_\varepsilon^{\gamma_q/2}K[\tilde w_\varepsilon]^{q-1}+\frac{1}{|x|^4}\varepsilon^{-\sigma-\tau(q-2)}\xi_\varepsilon^{2}\right\}.
\end{array}
\eqno(5.42)
$$
Let 
$$
a(x)=(1-\varepsilon^{-\tau(q-2)}\eta(x))\frac{\varepsilon^{-\sigma}\xi_\varepsilon^2}{|x|^4}, 
$$
 $$
b(x)=\frac{\eta(x)}{|x|^{4}}(I_\alpha\ast |\tilde w_\varepsilon|^{2_\alpha^*})(\frac{x}{|x|^2})\tilde w_\varepsilon^{2_\alpha^*-2}(\frac{x}{|x|^2})+\frac{\eta(x)}{|x|^{\gamma_q}}\varepsilon^{-\sigma}\xi_\varepsilon^{\gamma_q/2}K[\tilde w_\varepsilon]^{q-1}.
$$
Then (5.41) and (5.42) yields
$$
-\Delta K[\tilde w_\varepsilon]+ a(x)K[\tilde w_\varepsilon]= b(x)K[\tilde w_\varepsilon].
\eqno(5.43)
$$
Adopting the Moser iteration, 
arguing as in \cite[Proposition 5.14]{Ma-2},  we obtain the following

\smallskip

\noindent{\bf Proposition 5.13.}  {\it Assume   $N\ge 3$, $\alpha>N-4$ and $ 2<q<2^*$.
Then there exists a constant $C>0$ such that for large $\varepsilon>0$, there holds 
$$
\tilde w_\varepsilon(x)\le C(1+|x|)^{-(N-2)}, \qquad x\in \mathbb R^N.
$$}

\noindent{\bf Lemma 5.14.} { If $s>\frac{N}{N-2}$, then $\|\tilde w_\varepsilon\|_s^s\sim 1$ as $\varepsilon\to \infty$. Furthermore, $\tilde w_\varepsilon\to V_1$ in $L^s(\mathbb R^N)$ as $\varepsilon\to \infty$. }

\begin{proof}   Since $\tilde w_\varepsilon\to V_1$ in $L^{2^*}(\mathbb R^N)$, as in \cite[Lemma 4.6]{Moroz-1}, using the embeddings $L^{2^*}(B_1)\hookrightarrow L^s(B_1)$ we prove that
$\liminf_{\varepsilon\to \infty}\|\tilde w_\varepsilon\|_s^s>0$.  

  On the other hand,  by virtue of  Proposition 5.13,  there exists a constant $C>0$ such that for all large $\varepsilon>0$, 
$$
\tilde w_\varepsilon(x)\le \frac{C}{(1+|x|)^{N-2}}, \qquad \forall x\in \mathbb R^N,
$$
which together with the fact that $s>\frac{N}{N-2}$ implies that $\tilde w_\varepsilon$ is bounded in $L^s(\mathbb R^N)$  uniformly for large  $\varepsilon>0$,  and 
 by the dominated convergence theorem $\tilde w_\varepsilon\to V_1$ in $L^s(\mathbb R^N)$ as $\varepsilon\to \infty$. 
\end{proof}
  
  \vskip 3mm

\subsection{The Proof of Theorem 2.2}  We are in the position to prove Theorem 2.2.
\begin{proof}[Proof of Theorem 2.2]
Firstly, we  consider the cases $N=4$ and $N=3$. We  note that a result similar to Lemma 3.2 holds for  $\tilde w_\varepsilon$ and $\tilde J_\varepsilon$. By  Lemma 5.8, we also have  $\tau_2(\tilde w_\varepsilon)=\tau_2(v_\varepsilon)$. Therefore, by (5.35),  we get
$$
\begin{array}{rcl}
m_\infty&\le& \sup_{t\ge 0} \tilde J_\varepsilon((\tilde w_\varepsilon)_t)\\
&\mbox{}&+\varepsilon^{-\sigma}\tau_2( \tilde w_\varepsilon)^{\frac{N}{2}}\left\{\xi_\varepsilon^{N-\frac{N-2}{2}q}(\varepsilon^{-\tau}\xi_\varepsilon^{\frac{N-2}{2}})^q\int_{\mathbb R^N}G(\varepsilon^{\tau}\xi_\varepsilon^{-\frac{N-2}{2}}\tilde w_\varepsilon)-\frac{1}{2}\xi_\varepsilon^{2}\int_{\mathbb R^N}|\tilde w_\varepsilon|^2\right\}\\
&\mbox{}&+\varepsilon^{-\sigma}\tau_2( \tilde w_\varepsilon)^{\frac{N}{2}}\left\{\xi_\varepsilon^{N}\varepsilon^{-\tau q}\int_{\mathbb R^N}G(\varepsilon^{\tau}\xi_\varepsilon^{-\frac{N-2}{2}}\tilde w_\varepsilon)-\frac{1}{2}\xi_\varepsilon^{2}\int_{\mathbb R^N}|\tilde w_\varepsilon|^2\right\}\\
&\le &m_\varepsilon+\varepsilon^{-\sigma}\tau_2(\tilde w_\varepsilon)^{\frac{N}{2}}\frac{A}{q}\xi_\varepsilon^{N-\frac{N-2}{2}q}(\int_{\mathbb R^N}|\tilde w_\varepsilon|^q+o_\varepsilon(1)),
\end{array}
\eqno(5.44)
$$
which implies that
$$
\xi_\varepsilon^{N-\frac{N-2}{2}q}\left(\int_{\mathbb R^N}|\tilde w_\varepsilon|^q+o_\varepsilon(1)\right)\ge \varepsilon^{\sigma}\frac{q}{A\tau_2(w_\varepsilon)^{\frac{N}{2}}}\delta_\varepsilon.
$$
Hence, by Lemma  5.5,  we obtain
$$
\xi_\varepsilon^{N-\frac{N-2}{2}q}\left(\int_{\mathbb R^N}|\tilde w_\varepsilon|^q+o_\varepsilon(1)\right)\gtrsim \varepsilon^{\sigma}\delta_\varepsilon\gtrsim \left\{\begin{array}{rcl} 
 (\ln\varepsilon)^{-\frac{4-q}{q-2}},  \quad  \quad &{\rm if}& \   \ N=4,\\
 \varepsilon^{-\frac{(6-q)^2}{2(q-2)(q-4)}},   \qquad &{\rm if}& \    \  N=3.
 \end{array}\right.
\eqno(5.45)
 $$  
 Therefore, by Lemma 5.14, we have 
$$
\xi_\varepsilon\gtrsim  \left\{\begin{array}{rcl} 
 (\ln\varepsilon)^{-\frac{1}{q-2}},  \quad  \quad &{\rm if}& \   \ N=4,\\
 \varepsilon^{-\frac{6-q}{(q-2)(q-4)}},   \qquad &{\rm if}& \    \  N=3.
 \end{array}\right.
\eqno(5.46)
$$
On the other hand, if  $ N=3$, then by  (5.35) and Lemma 5.10 and Lemma 5.14, we have 
$$
\varepsilon^{\sigma/2}\xi_\varepsilon^{-1+\frac{q-2}{2}}\lesssim \xi_\varepsilon^{\frac{q-2}{2}}\|\tilde w_\lambda\|_2^2\lesssim 1.
$$
Then, noting that  $\sigma=\frac{2^*-q}{q-2}=\frac{6-q}{q-2}$,  for $q\in (4,6)$, we have 
$$
\xi_\varepsilon^{\frac{q-4}{2}}\lesssim \varepsilon^{-\frac{\sigma}{2}}=\varepsilon^{-\frac{6-q}{2(q-2)}},
$$
and hence
$$
\xi_\varepsilon\lesssim \varepsilon^{-\frac{6-q}{(q-2)(q-4)}}.
\eqno(5.47)
$$
If $N=4$, noting  that 
$$
-\ln(\varepsilon^{-\sigma}\xi_\varepsilon^{2})=\sigma\ln\varepsilon-2\ln\xi_\varepsilon\ge \sigma\ln\varepsilon,
$$
  by (5.35) and Lemma 5.11 and Lemma 5.14,  we have 
$$
\xi_\varepsilon^{q-2}\ln\varepsilon\lesssim \xi_\varepsilon^{q-2}\|\tilde w_\lambda\|_2^2\lesssim 1.
$$
Hence,  we  obtain 
$$
\xi_\varepsilon\lesssim   (\ln\varepsilon)^{-\frac{1}{q-2}}.
\eqno(5.48)
$$
Thus, it follows from (5.44), (5.47), (5.48)  and Lemma 5.14 that 
$$
\delta_\varepsilon=m_\infty-m_\varepsilon\lesssim \varepsilon^{-\sigma}\xi_\varepsilon^{N-\frac{N-2}{2}q}\lesssim \left\{\begin{array}{rcl} (\varepsilon\ln\varepsilon)^{-\frac{4-q}{q-2}},   \quad &{\rm if}& \   \ N=4,\\
 \varepsilon^{-\frac{6-q}{2(q-4)}},    \qquad &{\rm if}& \    \  N=3,
 \end{array}\right.
$$
which together with Lemma 5.5 implies that 
$$
\delta_\varepsilon\sim \varepsilon^{-\sigma}\xi_\varepsilon^{N-\frac{N-2}{2}q}\sim \left\{\begin{array}{rcl} (\varepsilon\ln\varepsilon)^{-\frac{4-q}{q-2}},   \quad &{\rm if}& \   \ N=4,\\
 \varepsilon^{-\frac{6-q}{2(q-4)}},  \qquad &{\rm if}& \    \  N=3.
 \end{array}\right.
\eqno(5.49)
$$
By (5.31) and  (5.35),  we  get
$$
\begin{array}{rcl}
m_\varepsilon&=&\frac{1}{2}\int_{\mathbb R^N}|\nabla \tilde w_\varepsilon|^2-\frac{1}{22_\alpha^*}\int_{\mathbb R^N}(I_\alpha\ast |\tilde w_\varepsilon|^{2_\alpha^*})|\tilde w_\varepsilon^{2_\alpha^*}\\
&\mbox{}& -\frac{A(N-2)(q-2)}{4q}\varepsilon^{-\sigma}\xi_\varepsilon^{N-\frac{N-2}{2}q}[\int_{\mathbb R^N}|\tilde w_\varepsilon|^q+o_\varepsilon(1)].
\end{array}
$$
By (5.33) and (5.34), we get
$$
\int_{\mathbb R^N}|\nabla \tilde w_\varepsilon|^2=\int_{\mathbb R^N}(I_\alpha\ast |\tilde w_\varepsilon|^{2_\alpha^*})|\tilde w_\varepsilon|^{2_\alpha^*}+\frac{NA(q-2)}{2q}\varepsilon^{-\sigma}\xi_\varepsilon^{N-\frac{N-2}{2}q}[\int_{\mathbb R^N}|\tilde w_\varepsilon|^q+o_\varepsilon(1)].
$$
Therefore, we have 
$$
m_\varepsilon=\frac{2+\alpha}{2(N+\alpha)}\int_{\mathbb R^N}|\nabla \tilde w_\varepsilon|^2-\frac{A\alpha(N-2)(q-2)}{4q(N+\alpha)}\varepsilon^{-\sigma}\xi_\varepsilon^{N-\frac{N-2}{2}q}[\int_{\mathbb R^N}|\tilde w_\varepsilon|^q+o_\varepsilon(1)].
$$
Similarly, we have
$$
m_\infty=\frac{2+\alpha}{2(N+\alpha)}\int_{\mathbb R^N}|\nabla V_1|^2.
$$
Thus, by virtue of (5.49), we obtain
$$
\begin{array}{rcl}
\|\nabla V_1\|_2^2-\|\nabla \tilde w_\varepsilon\|_2^2
&=&\frac{2(N+\alpha)}{2+\alpha}\delta_\varepsilon-\frac{A\alpha(N-2)(q-2)}{2q(2+\alpha)}\varepsilon^{-\sigma}\xi_\varepsilon^{N-\frac{N-2}{2}q}[\int_{\mathbb R^N}|\tilde w_\varepsilon|^q+o_\varepsilon(1)]\\
&=& \left\{\begin{array}{rcl} O((\varepsilon\ln\varepsilon)^{-\frac{4-q}{q-2}}),   \quad &{\rm if}& \   \ N=4,\\
O( \varepsilon^{-\frac{6-q}{2(q-4)}}),   \qquad &{\rm if}& \    \  N=3.
 \end{array}\right.
 \end{array}
\eqno(5.50)
$$
By (5.31) and (5.34), we have 
$$
\begin{array}{lcl}
m_\varepsilon &=&(\frac{1}{2}-\frac{1}{2^*})\int_{\R^N}|\nabla \tilde w_\varepsilon|^2+(\frac{1}{2^*}-\frac{1}{22_\alpha^*})\int_{\R^N}(I_\alpha\ast |\tilde w_\varepsilon|^{2_\alpha^*})|\tilde w_\varepsilon|^{2_\alpha^*}\\
&=&\frac{1}{N}\int_{\R^N}|\nabla \tilde w_\varepsilon|^2+\frac{\alpha(N-2)}{2N(N+\alpha)}\int_{\R^N}(I_\alpha\ast |\tilde w_\varepsilon|^{2_\alpha^*})|\tilde w_\varepsilon|^{2_\alpha^*}.
\end{array}
$$
Similarly, we also have 
$$
m_\infty=\frac{1}{N}\int_{\R^N}|\nabla V_1|^2+\frac{\alpha(N-2)}{2N(N+\alpha)}\int_{\R^N}(I_\alpha\ast |V_1|^{2_\alpha^*})|V_1|^{2_\alpha^*}.
$$
Then it follows from (5.49) and (5.50)  that 
$$
\begin{array}{rl}
&\int_{\R^N}(I_\alpha\ast |V_1|^{2_\alpha^*})|V_1|^{2_\alpha^*}-\int_{\R^N}(I_\alpha\ast |\tilde w_\varepsilon|^{2_\alpha^*})|\tilde w_\varepsilon|^{2_\alpha^*}\\
&=\frac{2N(N+\alpha)}{\alpha(N-2)}\left[(m_\infty-m_\varepsilon)-\frac{1}{N}(\int_{\R^N}|\nabla V_1|^2-\int_{\R^N}|\nabla \tilde w_\varepsilon|^2)\right]\\
&=\frac{2N(N+\alpha)}{\alpha(N-2)}\delta_\varepsilon-\frac{2(N+\alpha)}{\alpha(N-2)}(\int_{\R^N}|\nabla V_1|^2-\int_{\R^N}|\nabla \tilde w_\varepsilon|^2)\\
&=\left\{\begin{array}{rcl} O((\varepsilon\ln\varepsilon)^{-\frac{4-q}{q-2}}),   \quad &{\rm if}& \   \ N=4,\\
 O(\varepsilon^{-\frac{6-q}{2(q-4)}}),    \qquad &{\rm if}& \    \  N=3.
 \end{array}\right.
\end{array}
\eqno(5.51)
$$
Since 
$
\|\nabla V_1\|_2^2=\int_{\R^N}(I_\alpha\ast |V_1|^{2_\alpha^*})|V_1|^{2_\alpha^*}=S_\alpha^{\frac{N+\alpha}{2+\alpha}},
$
it follows from (5.61) and (5.62) that 
$$
\|\nabla \tilde w_\varepsilon\|_2^2=S_\alpha^{\frac{N+\alpha}{2+\alpha}}+\left\{\begin{array}{rcl} O((\varepsilon\ln\varepsilon)^{-\frac{4-q}{q-2}}),   \quad &{\rm if}& \   \ N=4,\\
 O(\varepsilon^{-\frac{6-q}{2(q-4)}}),   \qquad &{\rm if}& \    \  N=3,
 \end{array}\right.
$$
and 
$$
\int_{\R^N}(I_\alpha\ast |\tilde w_\varepsilon|^{2_\alpha^*})|\tilde w_\varepsilon|^{2_\alpha^*}=S_\alpha^{\frac{N+\alpha}{2+\alpha}}+\left\{\begin{array}{rcl} O((\varepsilon\ln\varepsilon)^{-\frac{4-q}{q-2}}),   \quad &{\rm if}& \   \ N=4,\\
 O(\varepsilon^{-\frac{6-q}{2(q-4)}}),  \qquad &{\rm if}& \    \  N=3.
 \end{array}\right.
$$
Finally, by (5.35), Lemma 5.10 and Lemma 5.11,  we obtain
$$
\|\tilde w_\varepsilon\|_2^2\sim \left\{\begin{array}{rcl}
\ln\varepsilon,   \quad  \quad if \   \ N=4,\\
\varepsilon^{\frac{6-q}{2(q-4)}},  \quad if  \    \  N=3.
\end{array}\right.
$$

For $N\ge 5$, the conclusion follows directly from Lemmas 5.5, 5.6 and 5.8.  More precisely, Replacing $V_1, \tilde w_\varepsilon$ by $V_{\rho_0}$ and $v_\varepsilon$ respectively,   and  setting $\xi_\varepsilon=1$ in the above arguments, we obtain
$$
\|\nabla u_\varepsilon \|_2^2=\|\nabla v_\varepsilon\|_2^2=S_\alpha^{\frac{N+\alpha}{2+\alpha}}+O(\varepsilon^{-\frac{2N-q(N-2)}{(N-2)(q-2)}})
$$
and 
$$
\int_{\mathbb R^N}(I_\alpha\ast |u_\varepsilon|^{2_\alpha^*})|u_\varepsilon|^{2_\alpha^*}=\int_{\mathbb R^N}(I_\alpha\ast |v_\varepsilon|^{2_\alpha^*})|v_\varepsilon|^{2_\alpha^*}=S_\alpha^{\frac{N+\alpha}{2+\alpha}}+O(\varepsilon^{-\frac{2N-q(N-2)}{(N-2)(q-2)}}).
$$
The statements on $u_\varepsilon$ follow from the corresponding results on $v_\varepsilon$ and $\tilde w_\varepsilon$. This completes the proof of Theorem 2.2.
\end{proof}

\vskip 5mm 

\section{  Final remarks }

 Cleatly, if $u_\varepsilon\in H^1(\mathbb R^N)$ is a ground state of $(P_\varepsilon)$, and for some $a>0$ there holds
$$
\|u_\varepsilon\|_2^2=a^2,
\eqno(6.1)
$$ 
then $u_\varepsilon$ is a normalized solution of (1.6) with $\lambda=-\varepsilon$.  We  denote this normalized solution by  a pair $(u_a, \lambda_a)$ with $\lambda_a=-\varepsilon$, or $u_a$ for simplicity. 

\smallskip

As direct consequences of  Theorem 2.1 and Theorem 2.2, we have the following results.
\smallskip

\noindent{\bf Proposition  6.1.}   {\it Assume that (H1) and (H2) hold. If $p=\frac{N+\alpha}{N}$, $q\in (2,2+\frac{4}{N})$, then for any large  $a>0$, the problem (1.6) has at least one positive  normalized solution $u_a\in H^1(\mathbb R^N)$, which is  radially symmetric and radially nonincreasing. 
Moreover,  as $a\to \infty$,  there hold
$$
 \|\nabla u_a \|_2^2\sim 
a^{\frac{2(2N-q(N-2))}{4-N(q-2)}}\to +\infty,
$$
$$
E(u_a)\sim \left\{\begin{array}{rcl}
-a^{\frac{2(N+\alpha)}{N}}\to -\infty,  \qquad      &{\rm if}&   q\in (2,2+\frac{4\alpha}{N(2+\alpha)}),\\
-a^{\frac{2(2N-q(N-2))}{4-N(q-2)}}\to -\infty,   &{\rm if} & q\in (2+\frac{4\alpha}{N(2+\alpha)},2+\frac{4}{N}).
 \end{array}\right. 
$$}

\smallskip

\noindent{\bf Proposition 6.2.}  {\it Assume that (H1) and (H2) hold. If $p=\frac{N+\alpha}{N-2}$, $q\in (2,2^*)$ for $N\ge 4$ and $q\in (4,6)$ for $N=3$, then for any small $a>0$,  the problem  (1.6)  has at least  one  positive normalized solutions  $u_a\in H^1(\mathbb R^N)$,  which are  radially symmetric and radially nonincreasing. Moreover,  as $a\to 0$, there hold
$$
\|\nabla u_a\|_2^2\simeq S_\alpha^{\frac{N+\alpha}{2+\alpha}},  \quad 
E(u_a)\simeq \frac{2+\alpha}{2(N+\alpha)}S_\alpha^{\frac{N+\alpha}{2+\alpha}}.
$$}

\vskip 3mm

\bigskip
{\small
\noindent {\bf Acknowledgements.} 
S.M. was supported by National Natural Science Foundation of China
(Grant Nos.11571187, 11771182)


\vskip 10mm

\begin {thebibliography}{99}
\footnotesize

\bibitem{Akahori-2}
T. Akahori, S. Ibrahim, N. Ikoma, H. Kikuchi and H. Nawa, 
{\it Uniqueness and nondegeneracy of ground states to nonlinear scalar field equations involving the Sobolev critical exponent in their nonlinearities for high frequencies},
 Calc. Var. Partial Differential Equations,  {\bf 58} (2019), Paper No. 120, 32 pp.

\bibitem{Akahori-3} 
T. Akahori, S. Ibrahim, H. Kikuchi and H. Nawa, 
Global dynamics above the ground state energy for the combined power type nonlinear Schr\"odinger equations with energy critical growth at low frequencies,
Memoirs of the AMS,  {\bf 272} (2021), 1331.


\bibitem{Akahori-4} T. Akahori, S. Ibrahim, H. Kikuchi, H. Nawa, {\it Non-existence of ground states and gap of variational values for 3D
Sobolev critical nonlinear scalar field equations},  J. Differential Equations,  {\bf 334} (2022),  25--86.



\bibitem{ASM2012} C. Alves, M. Souto, M. Montenegro, {\it Existence of a ground state solution for a nonlinear scalar field
equation with critical growth},  Calc. Var. Partial Differ. Equ., {\bf  43} (2012), 537--554.

\bibitem{Berestycki-1}  H. Berestycki and P.-L. Lions,  {\it Nonlinear scalar field equations. I. Existence of a ground state}, 
Archive for Rational Mechanics and Analysis,  {\bf  82} (1983), 313--345.

\bibitem{Brezis-1}    H. Brezis and E. Lieb,  {\it A relation between pointwise convergence of functions and convergence of functionals},  Proc. Amer. Math. Soc.,  {\bf 88} (1983),  486--490.


   
 \bibitem{Cazenave-1}  T. Cazenave and P.-L. Lions, {\it Orbital stability of standing waves
for some nonlinear Schr\"odinger equations},  Comm. Math. Phys.,
{\bf 85} (1982), 549--561.   
   
\bibitem{Coles}
M. Coles, S. Gustafson,
 {\it Solitary Waves and Dynamics for Subcritical Perturbations of Energy Critical NLS},
arXiv:1707.07219 (to appear in Differ. Integr. Equ.).


\bibitem{GT}   D. Gilbarg and N. S. Trudinger,  Elliptic Partial Differential Equations of Second Order. Grundlehren der Mathematischen Wissenschaften, Second edn. Springer, Berlin (1983).

\bibitem{Jeanjean-2} L. Jeanjean, {\it Existence of solutions with prescribed norm for semilinear elliptic equations},  Nonlinear Anal., {\bf 28} (1997), 1633--1659.


\bibitem{Jeanjean-3}   L. Jeanjean and   T. Le, {\it Multiple normalized solutions for a Sobolev critical Schr\"odinger equation},  Mathematische Annalen, https://doi.org/10.1007/s00208-021-02228-0.

   
\bibitem{Jeanjean-4}   L. Jeanjean, J. Zhang and X. Zhong, {\it A global branch approach to
normalized solutions for the Schr\"odinger equation}, arXiv:2112.05869v1.


\bibitem{Jeanjean-5}  L. Jeanjean, J. Jendrej, T. Le and N. Visciglia, {\it Orbital stability of ground states for a Sobolev critical schr\"odinger equation},
arXiv:2008.12084v1.
   
\bibitem{Lewin-1} 
M. Lewin and S. R. Nodari, 
{\it The double-power nonlinear Schr\"odingger equation and its generalizations: uniqueness, non-degeneracy and applications}, 
arXiv:2006.02809v1.

\bibitem{Li-2} Xinfu Li, Shiwang Ma and Guang Zhang, {\it Existence and qualitative properties of solutions for
Choquard equations with a local term},  Nonlinear Analysis: Real World Applications, {\bf  45}
(2019), 1--25.   

\bibitem{Li-1} Xinfu Li and Shiwang Ma, {\it Choquard equations with critical nonlinearities},  Commun. Contemp. Math., {\bf 22}
(2019), 1950023.

\bibitem{Li-3} Xinfu Li, {\it Existence and symmetry of normalized ground state to Choquard equation with
local perturbation}, arXiv:2103.07026v1, 2021.

\bibitem{Li-4} Xinfu Li, {\it Standing waves to upper critical Choquard equation with a local perturbation: multiplicity, qualitative properties and stability}, arXiv:2104.09317v1, 19 Apr 2021.

\bibitem{Li-5} Xinfu Li, {\it Nonexistence, existence and symmetry of normalized ground states to Choquard equations with a local perturbation}, arXiv:2103.07026v2 [math.AP] 7 May 2021.

\bibitem{Li-6} Xinfu Li, Jianguang Bao and  Wenguang Tang, {\it Normalized solutions to lower critical Choquard equation with a local perturbation}, 
arXiv:2207.10377v2 [math.AP] 18 Aug 2022.


\bibitem{LiuXQ} X.Q. Liu, J. Q. Liu and Z.-Q. Wang,  {\it Quasilinear elliptic equations with critical growth via perturbation method},  J. Differential Equations,  {\bf 254} (2013), 102--124.

\bibitem{Liu-1} Zeng Liu and V. Moroz, {\it Limit profiles for singularly perturbed Choquard equations with local repulsion}, arXiv:2107.05065v2.


\bibitem{Lieb-Loss 2001} E.H. Lieb and  M. Loss,  Analysis, volume 14 of graduate studies in
mathematics, American Mathematical Society, Providence, RI, (4)
2001.

\bibitem{LLT2017} J. Liu, J.-F. Liao, C.-L. Tang,  {\it Ground state solution for a class of Schrödinger equations involving
general critical growth term},  Nonlinearity, {\bf  30} (2017), 899--911.


\bibitem{Lions-1}   P. H. Lions, {\it The concentration-compactness principle in the calculus of variations: The locally compact cases, Part I and Part II},  Ann. Inst. H. Poincar\'{e} Anal. Non Lin\'{e}aire, {\bf 1} (1984), 223--283.

\bibitem{Ma-1} S. Ma and V. Moroz,  {\it Asymptotic profiles for a nonlinear Schr\"odinger equation with critical combined
powers nonlinearity},  Mathematische Zeitschrift,  (2023) 304:13.

\bibitem{Ma-2} S. Ma and V. Moroz, {\it Asymptotic profiles for Choquard  equations  with  combined
attractive  nonlinearities}, arXiv:2302.13727v1 [math.AP] 27 Feb 2023.


\bibitem{Moroz-1}   V. Moroz and C. B. Muratov,  {\it Asymptotic properties of ground states of scalar field equations with a vanishing
parameter}, J. Eur. Math. Soc.,  {\bf 16} (2014), 1081--1109.


\bibitem {Moroz-2}   V. Moroz  and J. Van Schaftingen, {\it Groundstates of nonlinear Choquard equations: Existence, qualitative properties and decay asymptotics}, 
J. Functional Analysis, {\bf  265}  (2013),  153--184.

\bibitem{MV2017} V. Moroz and J. Van Schaftingen, {\it A guide to the Choquard equation,} J. Fixed Point Theory Appl. 19 (2017) 773–813. 

 \bibitem{Soave-1} N. Soave, {\it Normalized ground states for the NLS equation with combined nonlinearities: The Sobolev critical case}, J. Functional  Analysis, {\bf 279} (2020), 108610.
 
\bibitem{Soave-2}  N. Soave, {\it  Normalized ground state for the NLS equations with combined nonlinearities},  J. Differential Equations, {\bf 269} (2020), 6941--6987.



\bibitem{Tao}
T. Tao, M. Visan and X. Zhang,
{\it The nonlinear Schr\"odinger equation with combined power-type nonlinearities},
Commun. Partial Differ. Equ.,  {\bf 32} (2007), 1281--1343.


\bibitem {Tintarev-1}  K. Tintarev and K-H. Fieseler,  Concentration Compactness: functional-Analytic Grounds and Applications, Imperial College Press, 57 Shelton Street, London, 2007.

\bibitem{VX-1} J. Van Schaftingen and Jiankang Xia, {\it Groundstates for a local nonlinear perturbation of the Choquard
equations with lower critical exponent}, J. Math. Anal. Appl., {\bf 464} (2018), 1184--1202. 


\bibitem{Wei-1}  J. Wei and Y. Wu, {\it Normalized solutions for Schr\"odinger equations with critical Soblev exponent and mixed nonlinearities},  J. Functional   Analysis, {\bf  283} (2022) 109574.

\bibitem{Wei-2}  J. Wei and Y. Wu,    {\it  On some nonlinear Schr\"odinger equations in $\mathbb R^N$},  Proc.  Royal Soc.  Edinburgh, page 1 of 26,
DOI:10.1017/prm.2022.56


\bibitem{Sun-1} Shuai Yao, Juntao Sun and Tsung-fang Wu,  {\it Normalized solutions for the Schr\"odinger equation with combined Hartree type and power nonlinearities}, 
arXiv:2102.10268v1, 2021.

\bibitem{Sun-2}   Shuai Yao, Haibo Chen, 
V. D. R\u adulescu and Juntao Sun,   {\it Normalized solutions for lower critical Choquard equations with critical Sobolev perturbation},    SIAM J. Math. Anal.,  
{\bf 54} (2022),  3696--3723.

\bibitem{ZZ2012}  J.  Zhang,  W. Zou,  {\it A Berestycki-Lions theorem revisited},  Commun. Contemp. Math., {\bf 14} (2012), 1250033. 



\end {thebibliography}
\end {document}